\documentclass[]{article}
\usepackage{graphicx}
\usepackage{graphics}
\usepackage{subfigure}
\usepackage{amssymb}
\usepackage{amsmath}
\usepackage{ulem}
\usepackage{float}
\usepackage{tabu}
\usepackage{multirow}
\usepackage{pifont}
\usepackage{color}
\usepackage{authblk}

\newenvironment{sequation}{\begin{equation}\small}{\end{equation}}

\title{\textbf{The Boundary Element Method of Peridynamics}}

\date{}

\author[a,d]{\sffamily \small Xue Liang}
\author[b]{\sffamily \small Linjuan Wang}
\author[c]{\sffamily \small Jifeng Xu}
\author[a,d]{\sffamily \small Jianxiang Wang \thanks{Corresponding author: jxwang@pku.edu.cn}}

\affil[a]{\scriptsize State Key Laboratory for Turbulence and Complex Systems, Department of Mechanics and Engineering Science, College of Engineering, Peking University, Beijing 100871, P. R. China}
\affil[b]{\scriptsize School of Astronautics, Beihang University, Beijing, 100191, P. R. China}
\affil[c]{\scriptsize Beijing Aeronautical Science and Technology Research Institute, Beijing 100083, P. R. China}
\affil[d]{\scriptsize CAPT-HEDPS, and IFSA Collaborative Innovation Center of MoE, College of Engineering, Peking University, Beijing 100871, P. R. China}

\begin{document}
	
\maketitle

\begin{abstract}

\quad The peridynamic theory brings advantages in dealing with discontinuities, dynamic loading, and non-locality. The integro-differential formulation of peridynamics poses challenges to numerical solutions of complicated and practical problems. Some important issues attract much attention, such as the computation of infinite domains, the treatment of softening of boundaries due to an incomplete horizon, and time accumulation error in dynamic processes. In this work, we develop the \textit{boundary element method of peridynamics} (PD-BEM). The numerical examples demonstrate that the PD-BEM exhibits several features. First, for non-destructive cases, the PD-BEM can be one to two orders of magnitude faster than the meshless particle method of peridynamics (PD-MPM) that directly discretizes the computational domains; second, it eliminates the time accumulation error, and thus conserves the total energy much better than the PD-MPM; third, it does not exhibit spurious boundary softening phenomena. For destructive cases where new boundaries emerge during the loading process, we propose a coupling scheme where the PD-MPM is applied to the cracked region and the PD-BEM is applied to the un-cracked region such that the time of computation can be significantly reduced.
\par\textbf{Keywords: } Peridynamic theory; Boundary element method; Linear elasticity; Coupling method
\end{abstract}

\section{Introduction}\label{atl1}

Peridynamics (PD) is a nonlocal theory in which the governing equation of motion of continuum mechanics is reformulated in an integro-differential form~\cite{tb1}. A material point interacts with others via bonds in the bond-based peridynamic theory or via a description of states in the state-based peridynamic theory~\cite{tb5}. PD facilitates dealing with discontinuous problems~\cite{tb90,tb85,tb86} because of the integral form. The theory can be applied to studies of microscale problems~\cite{tb2,tb3}, and multiscale problems~\cite{tb4,tb6}, and it can degenerate into the classical local theory when the characteristic length (i.e., the size of the horizon) approaches zero. Just like the classical continuum theory, solutions of various practical and complicated problems within the formalism of the peridynamic theory are achieved via numerical computations. Thus, various numerical methods have been developed, which can be divided into three categories.

First, a straightforward way to solve the integro-differential peridynamic governing equation is to discretize it with respect to time and spatial coordinates. In this regard, Silling et al. developed a meshfree particle method (PD-MPM) in which a continuum is discretized into particles, and the equations of motion of the particles are then solved~\cite{tb10}. Thereafter, a new integral strategy is adopted for improving the efficiency of computation~\cite{tb11,tb12}. Nonlinear problems are also explored~\cite{tb14}. This method requires a small time step because of numerical stability. For peridynamic static problems, Kilic et al. adopted an adaptive dynamic relaxation technique to transform the static problem into a dynamic one~\cite{tb15}, based on the fact that the static solution is regarded as the steady-state part of the transient solution whose initial condition is the boundary condition of the static problem, which effectively saves storage space. Littlewood solved the equilibrium equation with an iterative method~\cite{tb13}. The grid refinement technology, where the areas with a high strain energy density or discontinuity are selected to be refined, is introduced in order to improve the accuracy and efficiency~\cite{tb16,tb17,tb18}.

Second, considering that mature numerical methods in the classical continuum mechanics, such as the finite element method, and the finite volume method, are widely applied and readily available, a natural way to solve the PD governing equation is to develop the corresponding numerical methods by following the methodologies of the traditional numerical methods. Macek et al. first developed the finite element method based on the bond-based PD theory~\cite{tb7}. Dorduncu et al. generalized it to the state-based PD~\cite{tb8}. However, discontinuities are a difficult issue for the finite element method. Thus, Chen et al. developed the Galerkin finite element method of the PD theory based on the variational principle~\cite{tb19}.  Subsequently, Zhang et al. adopted a new integral rule for improving the efficiency of computation~\cite{tb20}. Gu et al.~\cite{tb21}, and Chen et al.~\cite{tb22} adopted  multigrid and preconditioned processing to accelerate the solution of the equation, respectively. Soghrati et al. introduced mesh refinement to reduce the cost of calculation~\cite{tb23}. Wang et al. presented a fast collocation method that shows lower computational complexity and less storage~\cite{tb25}. Tian et al. generalized the method to the problem of convection diffusion~\cite{tb24}. Besides the finite element method, other traditional numerical methods are also applied to PD. For example, Tian et al. introduced the finite difference method to PD~\cite{tb26}. Du et al. proposed a finite volume method based on the conservation rate~\cite{tb27}. Kilic et al., combining the finite element method and the PD-MPM~\cite{tb10}, developed a collocation point method with the background grid~\cite{tb28}.

While the above two categories of numerical methods show advantages in dealing with discontinuous problems, four issues have been encountered in the numerical solutions of the peridynamic formulation. The first issue is the weakening of the boundary~\cite{tb73}, which is caused by the incomplete horizon of the material points near the boundaries; the second is the infinite problem~\cite{tb56}, where artificial boundaries have to be introduced, and they pose the difficulty of artificial reflection of waves at these boundaries; the third is imposing boundary conditions of volume constraints~\cite{tb71} that introduce additional constraints compared with classical boundary conditions; the fourth is the computational efficiency.  Therefore, researchers have been trying to couple local numerical methods and nonlocal numerical methods. Thus, these coupling methods are the third category of numerical methods for the peridynamic theory. Specifically, nonlocal numerical methods are used in the regions where discontinuities exist, and local numerical methods are used in other regions. In this way, the advantages of the local numerical methods in dealing with boundary conditions and computational efficiency are combined with those of the nonlocal numerical methods in dealing with discontinuities. In order to achieve this coupling, it is necessary to communicate the information between the region where the local numerical method is implemented and the one where the nonlocal numerical method is implemented. According to the way where the information is transmitted, the coupling methods can be divided into two types. The first is that the coupling is realized with the aid of node force balance. The second is that the coupling is realized with the aid of deformation coordination. Coupling between the classical finite element method and the meshless particle class method in PD has been well achieved~\cite{tb32,tb33,tb34,tb35,tb36,tb87}, and has been extended to  the multiscale finite element and the extended finite element methods~\cite{tb37,tb38}. Considering the connection problem in the coupling zone and some extreme load conditions where the finite element method shows a low precision, computations that couple classical meshless methods with the PD-MPM have been proposed.
For example, the PD-MPM has been coupled with the SPH (Smoothed Particle Hydrodynamics)~\cite{tb39,tb40}, the MLEBF (Meshless Local Exponential Basis Functions)~\cite{tb42}, the FPM (Finite Point Method)~\cite{tb43}, and the molecular dynamics method~\cite{tb45}. Seleson et al. developed a model coupling classical continuum mechanics and peridynamics by force balance in the coupling area~\cite{tb31,tb89}. Han et al. developed a model by combining local stiffness of classical continuum mechanics with nonlocal stiffness of peridynamics in the coupling area~\cite{tb30,tb83}.

In this work, we develop the \textit{boundary element method of peridynamics} (PD-BEM) based on solution of the PD governing equation through boundary discretization. This method has its particular features in dealing with the above four issues. Firstly, the first and third issues are solved  by boundary discretization. Secondly, the boundary discrete numerical method improves the computational efficiency by dimension reduction, and facilitates computation of infinite domains, which addresses the second and the fourth issues.

The outline of the method is described as follows. We derive the boundary integral equation of bond-based peridynamics for static and dynamic problems by using the infinite Green function~\cite{tb58}, and the nonlocal operator theory~\cite{tb59}. The corresponding numerical framework, based on the PD boundary integral equation, is then constructed. For static problems, this method gives the explicit equation, which need not be solved iteratively. For dynamic cases, considering the problems of stability and time  accumulation error, we solve the problem in the Laplace domain by introducing the Laplace transformation and finally obtain the results in the time domain via inversion. This treatment brings two remarkable benefits. The first is that we can independently calculate the dynamic response in any position of the medium at any time, which eliminates time  accumulation error and reduces the computational cost, as shown in Sections \ref{atl42} and \ref{atl43}; the second is that we can easily implement parallel computation, as shown in Section \ref{atl34}. It is worth noting that though the present method is developed for the theory of the linearized bond-based PD, it can be generalized to the state-based peridynamic theory.


\section{The boundary integral equation}\label{atl2}

\subsection{The static boundary integral equation}\label{atl21}

First, we give a brief review of the bond-based peridynamic theory.  Its governing equation is
\begin{sequation}\label{fa1}
\rho \ddot{\mathbf{u}}\left(\mathbf{x},t\right) = \int_{\mathcal{H}_{\mathbf{x}}} \mathbf{f}\left(\mathbf{x},\mathbf{x}',t\right) {\rm d}V_{\mathbf{x}'} + \mathbf{b}\left(\mathbf{x},t\right)
\end{sequation}
where $\rho$ is density; $\mathbf{u}$ is the displacement that is the function of the position vector $\mathbf{x}$ of a material point and time $t$; $\mathcal{H}_{\mathbf{x}}$ is the family of point $\mathbf{x}$, which is a set of all points $\mathbf{x}'$ that interact  with point $\mathbf{x}$;   $\mathbf{f}\left(\mathbf{x},\mathbf{x}',t\right)$ is the pairwise force of point $\mathbf{x}'$ applied to point $\mathbf{x}$ at time $t$; and $\mathbf{b}$ is the body force density.

For linear elastic materials, (\ref{fa1}) is simplified to the following form:
\begin{sequation}\label{fa2}
\rho \ddot{\mathbf{u}}\left(\mathbf{x},t\right) = \int_{\mathcal{H}_{\mathbf{x}}} \mathbf{C} \left(\boldsymbol{\xi}\right) \cdot \left[\mathbf{u} \left(\mathbf{x} + \boldsymbol{\xi},t\right) - \mathbf{u} \left(\mathbf{x},t\right)\right] {\rm d}V_{\boldsymbol{\xi}} + \mathbf{b}\left(\mathbf{x},t\right)
\end{sequation}
where $\boldsymbol{\xi} = \mathbf{x}' - \mathbf{x}$; and $\mathbf{C}\left(\boldsymbol{\xi}\right)$ is the micromodulus tensor between point $\mathbf{x}$ and point $\mathbf{x}'$ and reflects the constitutive relation. Generally, $\mathbf{C}\left(\boldsymbol{\xi}\right)$ is a radial tensor function about $\boldsymbol{\xi}$. (\ref{fa2}) means that there is a spring connection between point $\mathbf{x}$ and each point in the family of points around it.

\begin{figure}[!htb]
	\centerline{\includegraphics[scale=0.3]{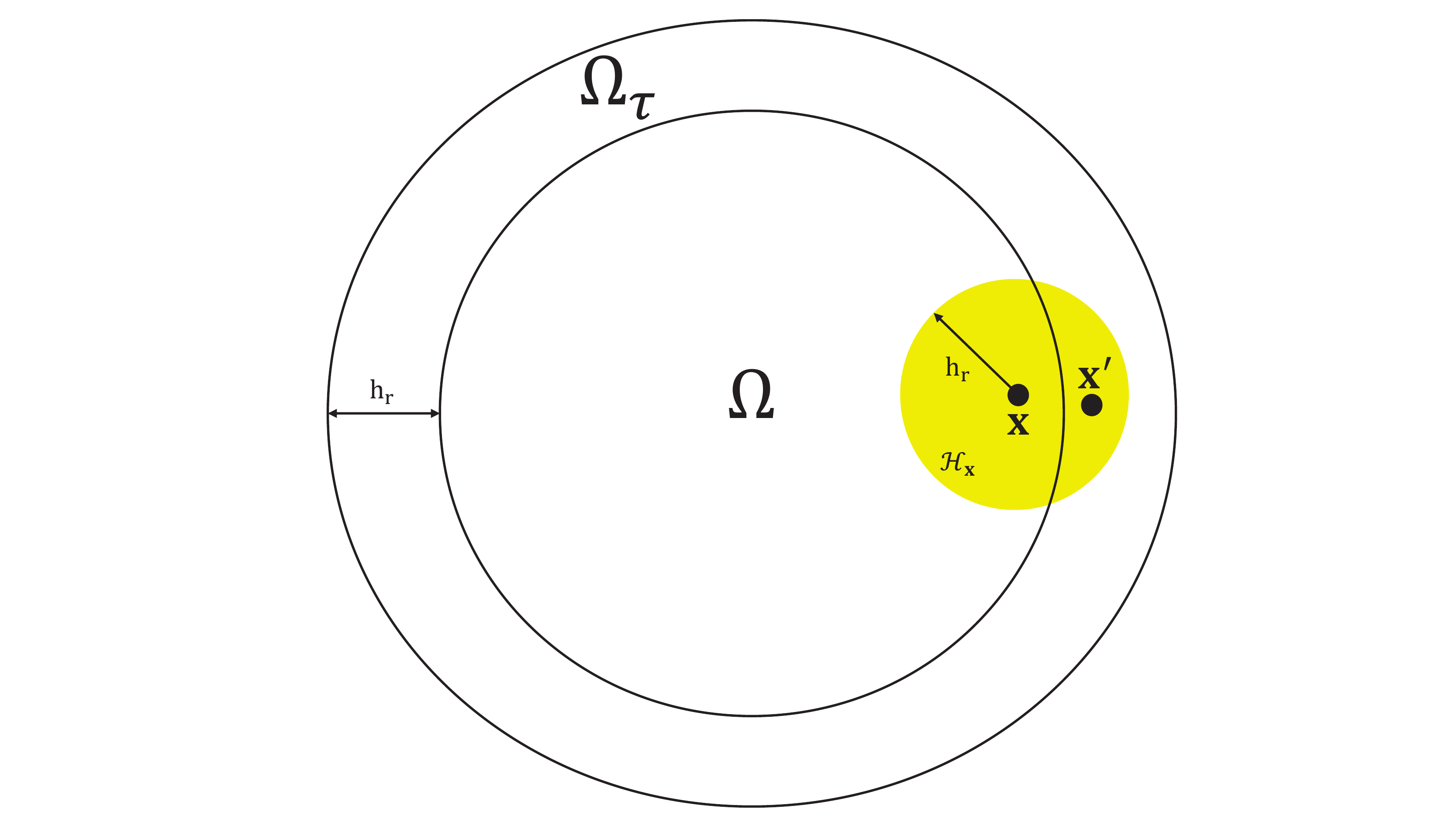}}	
	\caption{The definition of $\Omega_\tau$.\label{sc14}}
\end{figure}

Then, we start the derivation of the static boundary integral equation. Firstly, we derive a nonlocal reciprocal theorem. Secondly, we apply the reciprocal theorem to the state of the fundamental solution and the state of the actual problem, which derive the boundary integral equation. Our derivations here are based on the nonlocal operator theory~\cite{tb59}. Thus, it is necessary to recapitulate the nonlocal operator theory. Du et al. restructured the local differential operator in the form of an integral~\cite{tb59}. The corresponding relation to the classical theory, and the convergence to the classical theory are proved by them. In this theory, the interaction operator $\mathcal{N}$, the divergence operator $\mathcal{D}$ and its conjugate $\mathcal{D}^{*}$ are introduced as follows~\cite{tb59}:
\begin{sequation}\label{fa3}
\mathcal{N}\left(\mathbf{v}\right)\left(\mathbf{x}\right) \equiv - \int_{\Omega\cup\Omega_\tau} \left[\mathbf{v} \left(\mathbf{x},\mathbf{x}'\right) + \mathbf{v} \left(\mathbf{x}',\mathbf{x}\right)\right] \cdot \boldsymbol{\alpha} \left(\mathbf{x},\mathbf{x}'\right) {\rm d}V_{\mathbf{x}'} \qquad {\rm for} \  \mathbf{x} \in \Omega_\tau
\end{sequation}
\begin{sequation}\label{fa4}
\mathcal{D}\left(\mathbf{v}\right)\left(\mathbf{x}\right) \equiv \int_{R^{n}} \left[\mathbf{v} \left(\mathbf{x},\mathbf{x}'\right) + \mathbf{v} \left(\mathbf{x}',\mathbf{x}\right)\right] \cdot \boldsymbol{\alpha} \left(\mathbf{x},\mathbf{x}'\right) {\rm d}V_{\mathbf{x}'} \qquad {\rm for} \  \mathbf{x} \in R^{n}
\end{sequation}
\begin{sequation}\label{fa5}
\mathcal{D}^{*}\left(\mathbf{w}\right) \left(\mathbf{x},\mathbf{x}'\right) \equiv - \left[\mathbf{w}\left(\mathbf{x}'\right) - \mathbf{w}\left(\mathbf{x}\right)\right] \otimes \boldsymbol{\alpha} \left(\mathbf{x},\mathbf{x}'\right) \qquad {\rm for} \  \mathbf{x}, \mathbf{x}' \in R^{n}
\end{sequation}
where  $R^{n}$ denotes the $n$-dimensional Euclidean space; $\mathbf{v}\left(\mathbf{x},\mathbf{x}'\right)$ is a tensor function of two points $\mathbf{x}$ and $\mathbf{x}'$; $\boldsymbol{\alpha} \left(\mathbf{x}',\mathbf{x}\right)$ is the nonlocal weight function, which has the properties $\boldsymbol{\alpha} \left(\mathbf{x},\mathbf{x}'\right) = - \boldsymbol{\alpha} \left(\mathbf{x}',\mathbf{x}\right)$,
$\boldsymbol{\alpha} \left(\mathbf{x};\mathbf{x}'\right) = \mathbf{0} \ $ for $\ \mathbf{x}' \notin \mathcal{H}_{\mathbf{x}}$. $\mathcal{H}_{\mathbf{x}}$ is defined in (\ref{fa1}). $\mathbf{w}\left(\mathbf{x}'\right)$ is a vector function of point $\mathbf{x}'$. $\Omega$ represents the domain imposed by a body force; $\Omega_\tau$ is the volume-constrained boundary that has a nonzero measure. Normally, $\Omega_\tau$, whose definition is $\Omega_\tau = \left\{\mathbf{x}' \in R^{n} \setminus \Omega \ | \ \exists \mathbf{x} \in \Omega, s.t. \boldsymbol{\alpha} \left(\mathbf{x},\mathbf{x}'\right) \ne \mathbf{0} \right\}$, is a banded area whose width is a horizon near the local boundary $\partial\Omega$, and it meets $\Omega \cap \Omega_\tau = \emptyset$. The point in $\Omega$ only interacts with the point in $\Omega\cup\Omega_\tau$. An example of the definition of $\Omega_\tau$ is shown in Figure \ref{sc14},
where $\mathcal{H}_{\mathbf{x}}$ is a circular area with point $\mathbf{x}$ as its center and ${\rm h_r}$ as its radius; $\Omega$ is a  circular domain; $\Omega_\tau$ is the circular ring. The nonlocal operator $\mathcal{N}$ corresponds to the classical surface traction operator, and the nonlocal operator $\mathcal{D}$ corresponds to the classical divergence operator. Thus we can express the nonlocal Gauss theorem as follows~\cite{tb59}:
\begin{sequation}\label{fa6}
\int_{\Omega} \mathcal{D} \left(\mathbf{v}\right) \left(\mathbf{x}\right) {\rm d}V_{\mathbf{x}} = \int_{\Omega_\tau} \mathcal{N} \left(\mathbf{v}\right) \left(\mathbf{x}\right) {\rm d}V_{\mathbf{x}}
\end{sequation}
It should be pointed out that in the paper of Du et al.~\cite{tb59}, $\mathbf{v}$ is a vector, and $\mathbf{w}$ is a scalar. However, we need a second-order tensor for the equations of mechanics; thus, we have generalized the operators to tensors. The generalization can be easily proved in line with the paper of Du et al.~\cite{tb59}, so we do not repeat it here.

To facilitate deriving the boundary integral equation, we introduce two vectors $\boldsymbol{\chi} \left(\mathbf{x}\right)$ and
$ \boldsymbol{\kappa} \left(\mathbf{x},\mathbf{x}'\right)$, as well as a second-order tensor $\boldsymbol{\Psi} \left(\mathbf{x},\mathbf{x}'\right)$ such that $\boldsymbol{\kappa}\left(\mathbf{x},\mathbf{x}'\right) = \boldsymbol{\chi}\left(\mathbf{x}\right) \cdot \boldsymbol{\Psi}\left(\mathbf{x},\mathbf{x}'\right)$. Then, replacing $\mathbf{v}\left(\mathbf{x},\mathbf{x}'\right)$ in (\ref{fa3})-(\ref{fa4}) with $\boldsymbol{\kappa}\left(\mathbf{x},\mathbf{x}'\right)$, we can get
\begin{sequation}\label{faa1}
\begin{split}
\mathcal{D} \left(\boldsymbol{\kappa}\right) (\mathbf{x})
= \int_{R^{n}} \left[\boldsymbol{\kappa}(\mathbf{x},\mathbf{x}') + \boldsymbol{\kappa}(\mathbf{x}',\mathbf{x})\right] \cdot \boldsymbol{\alpha} \left(\mathbf{x},\mathbf{x}'\right) {\rm d}V_{\mathbf{x}'} \\
= \boldsymbol{\chi}(\mathbf{x}) \cdot \mathcal{D} \left(\boldsymbol{\Psi}\right)(\mathbf{x}) + \int_{\Omega\cup\Omega_\tau} \left[ \boldsymbol{\chi}
(\mathbf{x}') - \boldsymbol{\chi} (\mathbf{x})\right] \cdot \boldsymbol{\Psi} (\mathbf{x}',\mathbf{x})\cdot \boldsymbol{\alpha} \left(\mathbf{x},\mathbf{x}'\right) {\rm d}V_{\mathbf{x}'}
\end{split}
\end{sequation}
\begin{sequation}\label{faa2}
\begin{split}
\mathcal{N} \left({\boldsymbol{\kappa}}\right) (\mathbf{x})
= - \int_{\Omega\cup\Omega_\tau} \left[\boldsymbol{\kappa}(\mathbf{x},\mathbf{x}') + \boldsymbol{\kappa}(\mathbf{x}',\mathbf{x})\right] \cdot \boldsymbol{\alpha} \left(\mathbf{x},\mathbf{x}'\right) {\rm d}V_{\mathbf{x}'} \\
= \boldsymbol{\chi}  (\mathbf{x})\cdot \mathcal{N} \left(\boldsymbol{\Psi}\right) (\mathbf{x}) - \int_{\Omega\cup\Omega_\tau} \left[ \boldsymbol{\chi}
(\mathbf{x}') - \boldsymbol{\chi} (\mathbf{x})\right] \cdot \boldsymbol{\Psi} (\mathbf{x}',\mathbf{x})\cdot \boldsymbol{\alpha} \left(\mathbf{x},\mathbf{x}'\right) {\rm d}V_{\mathbf{x}'}
\end{split}
\end{sequation}
The difference between the integrals of (\ref{faa1}) and (\ref{faa2}) with respect to $\mathbf{x}$ is
\begin{sequation}\label{faa3}
\begin{split}
& \int_{\Omega} \mathcal{D} \left(\boldsymbol{\kappa}\right) {\rm d}V_{\mathbf{x}} - \int_{\Omega_\tau} \mathcal{N} \left({\boldsymbol{\kappa}}\right) {\rm d}V_{\mathbf{x}} \\
=& \int_{\Omega} \boldsymbol{\chi} \cdot \mathcal{D} \left( \boldsymbol{\Psi} \right) {\rm d}V_{\mathbf{x}} - \int_{\Omega_\tau} \boldsymbol{\chi} \cdot \mathcal{N} \left( \boldsymbol{\Psi} \right) {\rm d}V_{\mathbf{x}} - \int_{\Omega\cup\Omega_\tau} \int_{\Omega\cup\Omega_\tau}  \mathcal{D}^{*} \left( \boldsymbol{\chi} \right) \left(\mathbf{x},\mathbf{x}' \right)  : \boldsymbol{\Psi}\left(\mathbf{x},\mathbf{x}' \right) {\rm d}V_{\mathbf{x}'} {\rm d}V_{\mathbf{x}}
\end{split}
\end{sequation}
Following (\ref{fa6}), the left hand side of (\ref{faa3}) is equal to zero; thus we can get
\begin{sequation}\label{fa7}
\int_{\Omega} \boldsymbol{\chi} \cdot \mathcal{D} \left(\boldsymbol{\Psi}\right) {\rm d}V_{\mathbf{x}} - \int_{\Omega\cup\Omega_\tau} \int_{\Omega\cup\Omega_\tau} \mathcal{D}^{*} \left(\boldsymbol{\chi}\right) \left(\mathbf{x},\mathbf{x}'\right) : \boldsymbol{\Psi} \left(\mathbf{x},\mathbf{x}'\right) {\rm d}V_{\mathbf{x}'} {\rm d}V_{\mathbf{x}} = \int_{\Omega_\tau} \boldsymbol{\chi} \cdot \mathcal{N} \left(\boldsymbol{\Psi}\right) {\rm d}V_{\mathbf{x}}
\end{sequation}
According to the work of Du et al.~\cite{tb59}, the governing equation and the boundary conditions within the formalism of PD are expressed as
\begin{sequation}\label{fa8}
\begin{cases}
\mathcal{D} \left(\boldsymbol{\Theta} : \mathcal{D}^{*} \left(\mathbf{v}\right) \right) \left(\mathbf{x}\right) + \mathbf{b} \left(\mathbf{x}\right) = \mathbf{0} \qquad & \text{for} \ \mathbf{x} \in \Omega \\
\mathbf{v} \left(\mathbf{x}\right) = \mathbf{g}_{d} \left(\mathbf{x}\right) \qquad & \text{for} \ \mathbf{x} \in \Omega_{\tau_d} \\
\mathcal{N} \left(\boldsymbol{\Theta} : \mathcal{D}^{*} \left(\mathbf{v}\right) \right) \left(\mathbf{x}\right) = \mathbf{g}_{n} \left(\mathbf{x}\right) \qquad & \text{for} \ \mathbf{x} \in \Omega_{\tau_n}
\end{cases}
\end{sequation}
where $\boldsymbol{\Theta}$ is the fourth-order elasticity tensor, which satisfies $\Theta_{ijkl} = \Theta_{jikl} = \Theta_{ijlk} = \Theta_{klij}$ , and $\boldsymbol{\Theta}\left(\mathbf{x},\mathbf{x}'\right) = \boldsymbol{\Theta}\left(\mathbf{x}',\mathbf{x}\right)$. $\mathbf{v}$ is the displacement field. $\mathbf{b}$ is the body force density. $\Omega_{\tau_n}$ represents the force boundary, and $\Omega_{\tau_d}$ represents the displacement boundary. $\Omega_{\tau_n} \cap \Omega_{\tau_d} = \emptyset$, and $\Omega_{\tau_n} \cup \Omega_{\tau_d} = \Omega_\tau$. $\mathbf{g}_d$ is the displacement volume constraint, and $\mathbf{g}_n$ is the force volume constraint. The first equation of (\ref{fa8}) corresponds to the static version of the peridynamic governing equation (\ref{fa1}) where the inertia is neglected.

In order to derive the reciprocal theorem, we introduce two deformation states where the displacements are denoted by $\mathbf{v}_1$ and $\mathbf{v}_2$, the body force densities by $\mathbf{b}$ and $\mathbf{a}$, displacement volume constraints by $\mathbf{g}_{d_1}$ and $\mathbf{g}_{d_2}$, and  force volume constraints by $\mathbf{g}_{n_{1}}$ and $\mathbf{g}_{n_{2}}$, respectively. By applying (\ref{fa7}) to the cases $\boldsymbol{\chi} = \mathbf{v}_1, \boldsymbol{\Psi} = \boldsymbol{\Theta} : \mathcal{D}^{*} \left(\mathbf{v}_2\right)$, and $\boldsymbol{\chi} = \mathbf{v}_2, \boldsymbol{\Psi} = \boldsymbol{\Theta} : \mathcal{D}^{*} \left(\mathbf{v}_1\right)$, we obtain
\begin{sequation}\label{fa11}
\begin{split}
& \int_{\Omega} \mathbf{v}_1 \left(\mathbf{x}\right) \cdot \mathcal{D} \left(\boldsymbol{\Theta} : \mathcal{D}^{*} \left(\mathbf{v}_2\right)\right) \left(\mathbf{x}\right) {\rm d}V_{\mathbf{x}} \\
=& \int_{\Omega\cup\Omega_\tau} \int_{\Omega\cup\Omega_\tau} \mathcal{D}^{*} \left(\mathbf{v}_1\right) \left(\mathbf{x},\mathbf{x}'\right) : \boldsymbol{\Theta} \left(\mathbf{x},\mathbf{x}'\right) : \mathcal{D}^{*} \left(\mathbf{v}_2\right) \left(\mathbf{x},\mathbf{x}'\right) {\rm d}V_{\mathbf{x}'} {\rm d}V_{\mathbf{x}}\\
& + \int_{\Omega_\tau} \mathbf{v}_1 \left(\mathbf{x}\right) \cdot \mathcal{N} \left(\boldsymbol{\Theta} : \mathcal{D}^{*} \left(\mathbf{v}_2\right)\right) \left(\mathbf{x}\right) {\rm d}V_{\mathbf{x}}
\end{split}
\end{sequation}
\begin{sequation}\label{fa12}
\begin{split}
& \int_{\Omega} \mathbf{v}_2 \left(\mathbf{x}\right) \cdot \mathcal{D} \left(\boldsymbol{\Theta} : \mathcal{D}^{*} \left(\mathbf{v}_1\right)\right) \left(\mathbf{x}\right) {\rm d}V_{\mathbf{x}} \\
=& \int_{\Omega\cup\Omega_\tau} \int_{\Omega\cup\Omega_\tau} \mathcal{D}^{*} \left(\mathbf{v}_2\right) \left(\mathbf{x},\mathbf{x}'\right) : \boldsymbol{\Theta} \left(\mathbf{x},\mathbf{x}'\right) : \mathcal{D}^{*} \left(\mathbf{v}_1\right) \left(\mathbf{x},\mathbf{x}'\right) {\rm d}V_{\mathbf{x}'} {\rm d}V_{\mathbf{x}}\\
& + \int_{\Omega_\tau} \mathbf{v}_2 \left(\mathbf{x}\right) \cdot \mathcal{N} \left(\boldsymbol{\Theta} : \mathcal{D}^{*} \left(\mathbf{v}_1\right)\right) \left(\mathbf{x}\right) {\rm d}V_{\mathbf{x}}
\end{split}
\end{sequation}
Subtracting (\ref{fa12}) from (\ref{fa11}), and using $\Theta_{klij} = \Theta_{ijkl}$, we get the reciprocal theorem:
\begin{sequation}\label{fa10}
\begin{split}
& \int_{\Omega} \left[ \mathbf{v}_1 \left(\mathbf{x}\right) \cdot \mathcal{D} \left(\boldsymbol{\Theta} : \mathcal{D}^{*} \left(\mathbf{v}_2\right)\right) \left(\mathbf{x}\right) - \mathbf{v}_2 \left(\mathbf{x}\right) \cdot \mathcal{D} \left(\boldsymbol{\Theta} : \mathcal{D}^{*} \left(\mathbf{v}_1\right)\right) \left(\mathbf{x}\right)\right] {\rm d}V_{\mathbf{x}} \\
=& \int_{\Omega_\tau} \left[ \mathbf{v}_1 \left(\mathbf{x}\right) \cdot \mathcal{N} \left(\boldsymbol{\Theta} : \mathcal{D}^{*} \left(\mathbf{v}_2\right)\right) \left(\mathbf{x}\right) - \mathbf{v}_2 \left(\mathbf{x}\right) \cdot \mathcal{N} \left(\boldsymbol{\Theta} : \mathcal{D}^{*} \left(\mathbf{v}_1\right)\right) \left(\mathbf{x}\right) \right] {\rm d}V_{\mathbf{x}}
\end{split}
\end{sequation}

It is noted that the integration domain on the right hand side of (\ref{fa10}) is the volume constrained boundary that has a nonzero measure~\cite{tb59}. Thus, if we establish the boundary element method based on (\ref{fa10}), then we will lose the advantages of dealing with complex boundaries and dimension reduction, which are the very merits of the classical boundary element method.  Therefore, we introduce additional constraints to the volume constrained boundary~\cite{tb59}, which transform the integral in the volume constrained boundary into the one in the classical boundary for (\ref{fa10}). Firstly, the displacement at $\partial\Omega$ given in (\ref{fa8}) is denoted by $\overline{\mathbf{u}} \left(\mathbf{x}\right)$, that is
\begin{sequation}\label{fa13a}
\overline{\mathbf{u}} \left(\mathbf{x}\right)\equiv  \mathbf{v} \left(\mathbf{x}\right),
\ \ \ {\mathbf{x} \in \partial \Omega}
\end{sequation}
where $\mathbf{v} \left(\mathbf{x}\right)$ denotes the displacement field given by the solution of the problem stipulated in (\ref{fa8}). $\partial \Omega$ is the local boundary. Secondly, once the solution $\mathbf{v} \left(\mathbf{x}\right)$ is given, we can calculate the surface traction according to the definition of the Cauchy stress in PD~\cite{tb79,tb84}.
\begin{sequation}\label{fa13ab}
\overline{\mathbf{T}} \left(\mathbf{v}\right) \left(\mathbf{x}\right)\equiv \dfrac{1}{2} \mathbf{n} \cdot \mathcal{\mathbf{C}} : \left(\left[ \mathbf{F} \left(\mathbf{v}\right) \right]^{\rm{T}} + \mathbf{F} \left(\mathbf{v}\right) - 2\mathbf{I}\right),\ \ \  \ {\mathbf{x} \in \partial \Omega}
\end{sequation}
where $\mathbf{I}$ is the metric tensor. The superscript ${\rm{T}}$ denotes the transpose operation. $\overline{\mathbf{T}}$ is the nonlocal surface traction operator. $\mathcal{\mathbf{C}}$ is the elasticity tensor. $\mathbf{n}$ is the unit vector in the external normal direction of the geometric boundary. $\mathbf{F}$ is the deformation gradient that is introduced for PD~\cite{tb5}:
\begin{sequation}\label{fabb1}
\mathbf{F} \left(\mathbf{v}\right) = \left[\int_{\mathcal{H}_{\mathbf{x}}} \underline{\omega} \left \langle \boldsymbol{\xi} \right \rangle \left(\mathbf{v} \left(\mathbf{x} + \boldsymbol{\xi},t\right) - \mathbf{v} \left(\mathbf{x},t\right)\right) \otimes \boldsymbol{\xi} {\rm d}V_{\boldsymbol{\xi}}\right] \cdot \left[\int_{\mathcal{H}_{\mathbf{x}}} \underline{\omega} \left \langle \boldsymbol{\xi} \right \rangle \boldsymbol{\xi} \otimes \boldsymbol{\xi} {\rm d}V_{\boldsymbol{\xi}}\right]^{-1}
\end{sequation}
where $\underline{\omega} \left \langle \boldsymbol{\xi} \right \rangle$ is the nonlocal weight function; $\otimes$ denotes tensorial product; $-1$ is the inverse of a tensor.

\begin{figure}[!htb]
	\centerline{\includegraphics[scale=0.3]{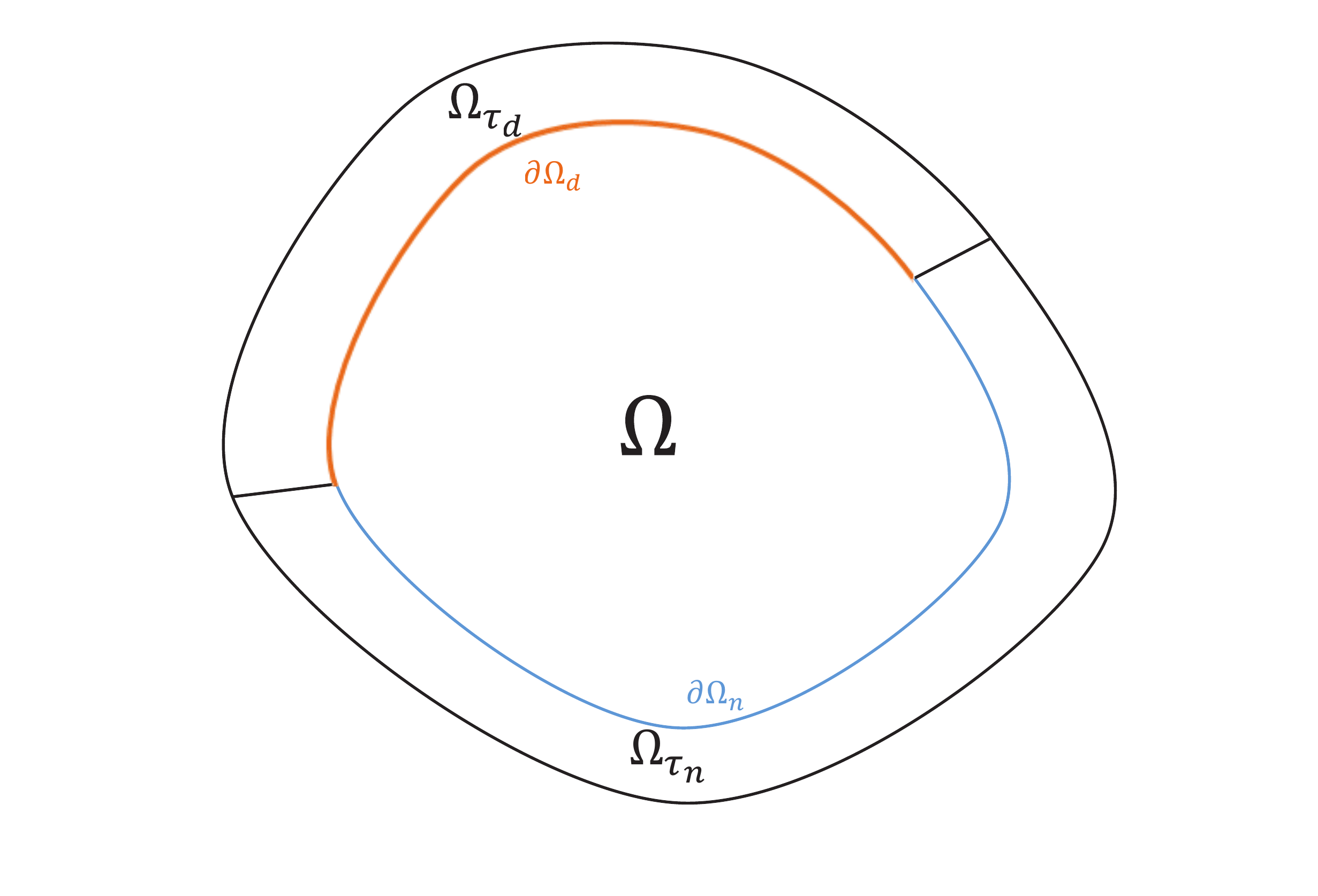}}	
	\caption{The diagram for local boundary.\label{sc11}}
\end{figure}

The boundary $\partial\Omega$ can be divided into two parts: the local displacement boundary $\partial \Omega_d$, which is defined as $\partial \Omega_d=Encl(\Omega_{\tau_{d}})\cap\Omega$, and the local force boundary $\partial \Omega_n$, which is defined as $\partial \Omega_n=Encl(\Omega_{\tau_{n}})\cap\Omega$. They are shown in Figure \ref{sc11}.
Therefore, $\overline{\mathbf{u}} \left(\mathbf{x}\right)$ is given in the local displacement boundary $\partial \Omega_d$ for the classical theory, and $\overline{\mathbf{T}} \left(\mathbf{v}\right) \left(\mathbf{x}\right)$ is given in the local force boundary $\partial \Omega_n$ for the classical theory. It is worth noting that we are developing the PD-BEM, so the nonlocal surface traction operator $\overline{\mathbf{T}} \left(\mathbf{v}\right) \left(\mathbf{x}\right)$ is given based on PD that is a special case of nonlocal theories. We introduce the following equivalency condition to convert the integral in the volume  domain $\Omega_{\tau}$ into the integral on $\partial\Omega$:
\begin{sequation}\label{fa14}
\int_{\Omega_\tau} \mathbf{v}_p \left(\mathbf{x}\right) \cdot \mathcal{N} \left(\boldsymbol{\Theta} : \mathcal{D}^{*} \left(\mathbf{v}\right)\right) \left(\mathbf{x}\right) {\rm d}V_{\mathbf{x}} = \int_{\partial \Omega_{n}} \mathbf{v}_p \left(\mathbf{x}\right) \cdot \overline{\mathbf{T}}\left(\mathbf{v}\right) \left(\mathbf{x}\right) {\rm d}V_{\mathbf{x}}+\int_{\partial \Omega_{d}} \mathbf{v}_p \left(\mathbf{x}\right) \cdot \overline{\mathbf{T}}\left(\mathbf{v}\right) \left(\mathbf{x}\right) {\rm d}V_{\mathbf{x}}
\end{sequation}
where $\mathbf{v}_p$ is a possible deformed state, and $\mathbf{v}$ is the actual deformed state which depends on $\mathbf{g}_{d}$  and $\mathbf{g}_{n}$.  The possible deformed state is any state that meets the following displacement constraint:
\begin{sequation}\label{fa15}
\mathbf{v}_p \left(\mathbf{x}\right) = \mathbf{g}_{d_{p}} \left(\mathbf{x}\right) \qquad  \text{for} \ \mathbf{x} \in \Omega_{\tau_d}
\end{sequation}
where $\mathbf{g}_{d_{p}}\left(\mathbf{x}\right)$ is any given displacement constraint. We note that the operator $\mathcal{N} \left(\boldsymbol{\Theta} : \mathcal{D}^{*} \left(\mathbf{v}\right)\right) \left(\mathbf{x}\right)$ is analogous to the body force density in the classical local theory, so the integral on the left hand side of (\ref{fa14}) represents the work done by the imposed body force on the possible deformation, and the right hand side represents the work done by the surface traction in the classical sense. Thus, the equivalency condition in (\ref{fa14}) can be regarded as the principle of virtual work that relates the volume constraint to the classical boundary conditions.

Substituting (\ref{fa15}) into (\ref{fa14}) yields
\begin{sequation}\label{fa16}
\begin{split}
& \int_{\Omega_{\tau_d}} \mathbf{g}_{d_{p}} \left(\mathbf{x}\right) \cdot \mathcal{N} \left(\boldsymbol{\Theta} : \mathcal{D}^{*} \left(\mathbf{v}\right)\right) \left(\mathbf{x}\right) {\rm d}V_{\mathbf{x}} + \int_{\Omega_{\tau_n}} \mathbf{v}_p \left(\mathbf{x}\right) \cdot \mathbf{g}_{n} \left(\mathbf{x}\right) {\rm d}V_{\mathbf{x}} \\
=& \int_{\partial \Omega_d} \overline{\mathbf{u}}_p \left(\mathbf{x}\right) \cdot T\left(\mathbf{v}\right) \left(\mathbf{x}\right) {\rm d}V_{\mathbf{x}} + \int_{\partial \Omega_n} \mathbf{v}_p \left(\mathbf{x}\right) \cdot \overline{\mathbf{T}} \left(\mathbf{v}\right) \left(\mathbf{x}\right) {\rm d}V_{\mathbf{x}}
\end{split}
\end{sequation}
Therefore, applying $\mathbf{v}_p=\mathbf{v}_1,\mathbf{v}=\mathbf{v}_2$ and $\mathbf{v}_p=\mathbf{v}_2,\mathbf{v}=\mathbf{v}_1$ in (\ref{fa14}), where
$\mathbf{v}_1$ and $\mathbf{v}_2$ represent the displacements of two problems with different body forces and different $\mathbf{g}_{d}$  and $\mathbf{g}_{n}$,
we obtain the following two equations:
\begin{sequation}\label{fa18}
\int_{\Omega_\tau} \mathbf{v}_1 \left(\mathbf{x}\right) \cdot \mathcal{N} \left(\boldsymbol{\Theta} : \mathcal{D}^{*} \left(\mathbf{v}_2\right)\right) \left(\mathbf{x}\right) {\rm d}V_{\mathbf{x}} = \int_{\partial \Omega} \mathbf{v}_1 \left(\mathbf{x}\right) \cdot \mathbf{T}\left(\mathbf{v}_2\right) \left(\mathbf{x}\right) {\rm d}V_{\mathbf{x}}
\end{sequation}
\begin{sequation}\label{fa19}
\int_{\Omega_\tau} \mathbf{v}_2 \left(\mathbf{x}\right) \cdot \mathcal{N} \left(\boldsymbol{\Theta} : \mathcal{D}^{*} \left(\mathbf{v}_1\right)\right) \left(\mathbf{x}\right) {\rm d}V_{\mathbf{x}} = \int_{\partial \Omega} \mathbf{v}_2 \left(\mathbf{x}\right) \cdot \mathbf{T}\left(\mathbf{v}_1\right) \left(\mathbf{x}\right) {\rm d}V_{\mathbf{x}}
\end{sequation}
Substituting (\ref{fa18}) and (\ref{fa19}) into the reciprocal theorem (\ref{fa10}) yields
\begin{sequation}\label{fa20}
\begin{split}
& \int_{\Omega} \mathbf{v}_1 \left(\mathbf{x}\right) \cdot \mathcal{D} \left(\boldsymbol{\Theta} : \mathcal{D}^{*} \left(\mathbf{v}_2\right)\right) \left(\mathbf{x}\right) - \mathbf{v}_2 \left(\mathbf{x}\right) \cdot \mathcal{D} \left(\boldsymbol{\Theta} : \mathcal{D}^{*} \left(\mathbf{v}_1\right)\right) \left(\mathbf{x}\right) {\rm d}V_{\mathbf{x}} \\
=& \int_{\partial \Omega} \mathbf{v}_1 \left(\mathbf{x}\right) \cdot \mathbf{T}\left(\mathbf{v}_2\right) \left(\mathbf{x}\right) - \mathbf{v}_2 \left(\mathbf{x}\right) \cdot \mathbf{T}\left(\mathbf{v}_1\right) \left(\mathbf{x}\right) {\rm d}V_{\mathbf{x}}
\end{split}
\end{sequation}

\begin{figure}[!htb]
	\centerline{\includegraphics[scale=0.3]{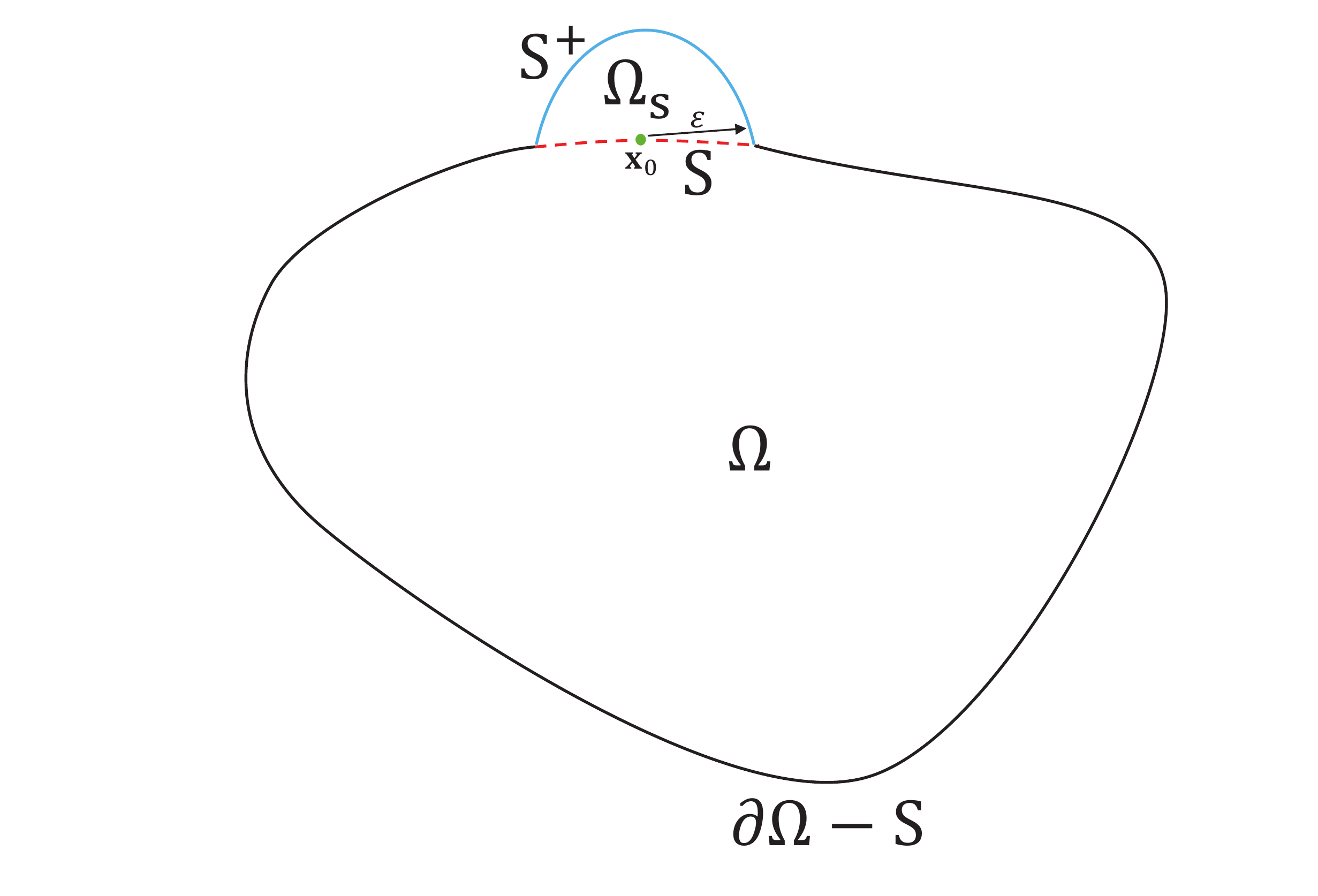}}	
	\caption{The geometric diagram for the boundary limit process.\label{sc10}}
\end{figure}

Next, we use the reciprocal theorem (\ref{fa20}) to derive the static boundary integral equation of PD, following the process in the classical theory~\cite{tb60}. We consider the case that $Problem\ \ 1$ in the reciprocal theorem is the actual problem, and $Problem\ \ 2$ is the fundamental solution $\overline{\mathbf{v}}_k$ of the infinite domain, whose value calculated in the considered finite domain is used in the PD-BEM. Then, (\ref{fa8}) is re-expressed as
\begin{sequation}\label{fa21}
\text{Problem \ 1}: \quad
\begin{cases}
\mathcal{D} \left(\boldsymbol{\Theta} : \mathcal{D}^{*} \left(\mathbf{u}\right) \right) \left(\mathbf{x}\right) + \mathbf{f} \left(\mathbf{x}\right) = \mathbf{0} \qquad \qquad \qquad \qquad \quad  & \text{for} \ \mathbf{x} \in \Omega \qquad \qquad \, \, \\
\mathbf{u} \left(\mathbf{x}\right) = \mathbf{g}_d \left(\mathbf{x}\right) \qquad \qquad \qquad \qquad \quad  & \text{for} \ \mathbf{x} \in \Omega_{\tau_d} \qquad \qquad \, \, \\
\mathcal{N} \left(\boldsymbol{\Theta} : \mathcal{D}^{*} \left(\mathbf{u}\right) \right) \left(\mathbf{x}\right) = \mathbf{g}_n \left(\mathbf{x}\right) \qquad \qquad \qquad \qquad \quad  & \text{for} \ \mathbf{x} \in \Omega_{\tau_n} \qquad \qquad \, \,
\end{cases}
\end{sequation}
\begin{sequation}\label{fa22}
\text{Problem \ 2}: \quad
\begin{cases}
\mathcal{D} \left(\boldsymbol{\Theta} : \mathcal{D}^{*} \left(\mathbf{v}_k\right) \right) \left(\mathbf{x} - \mathbf{x}_0\right) + \delta \left(\mathbf{x} - \mathbf{x}_0\right)\mathbf{e}_k = \mathbf{0} \qquad & \text{for} \ \mathbf{x} \in \Omega \\
\mathbf{v}_k \left(\mathbf{x} - \mathbf{x}_0\right) = \overline{\mathbf{v}}_k \left(\mathbf{x} - \mathbf{x}_0\right) \qquad & \text{for} \ \mathbf{x} \in \Omega_{\tau_d} \mathbf{x}_0 \in \Omega_{\tau_d}\\
\mathcal{N} \left(\boldsymbol{\Theta} : \mathcal{D}^{*} \left(\mathbf{v}_k\right) \right) \left(\mathbf{x} - \mathbf{x}_0\right) = \mathcal{N} \left(\boldsymbol{\Theta} : \mathcal{D}^{*} \left(\overline{\mathbf{v}}_k\right) \right) \left(\mathbf{x} - \mathbf{x}_0\right) \qquad & \text{for} \ \mathbf{x} \in \Omega_{\tau_n}
\end{cases}
\end{sequation}
where $\mathbf{e}_k$ is the coordinate base vector. $\overline{\mathbf{v}}_k$ satisfies
\begin{sequation}\label{faa22}
\mathcal{D} \left(\boldsymbol{\Theta} : \mathcal{D}^{*} \left(\overline{\mathbf{v}}_k\right) \right) \left(\mathbf{x} - \mathbf{x}_0\right) + \delta \left(\mathbf{x} - \mathbf{x}_0\right)\mathbf{e}_k = \mathbf{0}
\end{sequation}
where $\mathbf{x}_0\in \Omega; \mathbf{x}\in R^n$. $\delta \left(\mathbf{x} - \mathbf{x}_0\right)$ is the Dirac function. Applying $\mathbf{v}_1 = \mathbf{u}$ and $\mathbf{v}_2 = \mathbf{v}_k$ to (\ref{fa20}), we get
\begin{sequation}\label{fa23}
\begin{split}
& \int_{\Omega} \mathbf{u} \left(\mathbf{x}\right) \cdot \mathcal{D} \left(\boldsymbol{\Theta} : \mathcal{D}^{*} \left(\mathbf{v}_k\right)\right) \left(\mathbf{x} - \mathbf{x}_0\right) - \mathbf{v}_k \left(\mathbf{x} - \mathbf{x}_0\right) \cdot \mathcal{D} \left(\boldsymbol{\Theta} : \mathcal{D}^{*} \left(\mathbf{u}\right)\right) \left(\mathbf{x}\right) {\rm d}V_{\mathbf{x}} \\
=& \int_{\partial \Omega} \mathbf{u} \left(\mathbf{x}\right) \cdot \mathbf{T}\left(\mathbf{v}_k\right) \left(\mathbf{x} - \mathbf{x}_0\right) - \mathbf{v}_k \left(\mathbf{x} - \mathbf{x}_0\right) \cdot \mathbf{T}\left(\mathbf{u}\right) \left(\mathbf{x}\right) {\rm d}V_{\mathbf{x}}
\end{split}
\end{sequation}
Substituting the first equations of (\ref{fa21}) and (\ref{fa22}) into (\ref{fa23}), we get
\begin{sequation}\label{fa24}
\begin{split}
& \int_{\Omega} \mathbf{u} \left(\mathbf{x}\right) \cdot \left(-\delta \left(\mathbf{x} - \mathbf{x}_0\right)\mathbf{e}_k\right) {\rm d}V_{\mathbf{x}} + \int_{\Omega} \mathbf{v}_k \left(\mathbf{x} - \mathbf{x}_0\right) \cdot \mathbf{f} \left(\mathbf{x}\right) {\rm d}V_{\mathbf{x}} \\
=& \int_{\partial \Omega} \mathbf{u} \left(\mathbf{x}\right) \cdot \mathbf{T}\left(\mathbf{v}_k\right) \left(\mathbf{x} - \mathbf{x}_0\right) - \mathbf{v}_k \left(\mathbf{x} - \mathbf{x}_0\right) \cdot \mathbf{T}\left(\mathbf{u}\right) \left(\mathbf{x}\right) {\rm d}V_{\mathbf{x}}
\end{split}
\end{sequation}
We take the limit process $\mathbf{x}_0 \to \partial \Omega$ for (\ref{fa24}) to get the PD static boundary integral equation. The two-dimensional case is shown graphically in Figure \ref{sc10}, while the three-dimensional case is similar. Following the boundary element for the classical local theory, we transform the limit process for $\mathbf{x}_0$ into the limit process for the boundary geometry. We let $\mathbf{x}_0 \to \partial \Omega$ and divide $\partial \Omega$ into two parts $S$ and $\partial \Omega - S$. Then, we construct a small arc $S^{+}$ with  point $\mathbf{x}_0$ as the center of the circle and $\varepsilon$ as the radius, and use $\partial \Omega - S \cup S^{+}$ to replace $\partial \Omega$ as the new boundary. The geometric relationship is shown in Figure \ref{sc10}.
In this way, the integral on the new boundary is nonsingular. Then, we take the limit process $\varepsilon \to 0$ to make the new boundary to converge to the original one. Applying the new boundary and the corresponding limit process to (\ref{fa24}), we can get
\begin{sequation}\label{fa24a}
\begin{split}
& \lim_{\varepsilon \to 0} \left(\int_{\Omega + \Omega_s} \mathbf{u} \left(\mathbf{x}\right) \cdot \left(-\delta \left(\mathbf{x} - \mathbf{x}_0\right)\mathbf{e}_k\right) {\rm d}V_{\mathbf{x}}\right) + \int_{\Omega + \Omega_s} \mathbf{v}_k \left(\mathbf{x} - \mathbf{x}_0\right) \cdot \mathbf{f} \left(\mathbf{x}\right) {\rm d}V_{\mathbf{x}} \\
=& \lim_{\varepsilon \to 0} \left(\int_{\partial \Omega - S} \mathbf{u} \left(\mathbf{x}\right) \cdot \mathbf{T}\left(\mathbf{v}_k\right) \left(\mathbf{x} - \mathbf{x}_0\right) - \mathbf{v}_k \left(\mathbf{x} - \mathbf{x}_0\right) \cdot \mathbf{T}\left(\mathbf{u}\right) \left(\mathbf{x}\right) {\rm d}V_{\mathbf{x}}\right) \\
& + \lim_{\varepsilon \to 0} \left(\int_{S^{+}} \mathbf{u} \left(\mathbf{x}\right) \cdot \mathbf{T}\left(\mathbf{v}_k\right) \left(\mathbf{x} - \mathbf{x}_0\right) - \mathbf{v}_k \left(\mathbf{x} - \mathbf{x}_0\right) \cdot \mathbf{T}\left(\mathbf{u}\right) \left(\mathbf{x}\right) {\rm d}V_{\mathbf{x}}\right)
\end{split}
\end{sequation}
where $\Omega_{S}$ is the area surrounded by $S$ and $S^{+}$. The limit of the left hand side of (\ref{fa24a}) is $-u_k\left(\mathbf{x}_0\right)$ because of the Dirac function. The  limit of the first term on the right hand side is the Cauchy principal value for the integral in $\partial \Omega$. For the limit of the second term on the right hand side, we need to analyze the singularity of the integrand and the integral interval in order to determine the limit. The integrand has two kinds of singularity: the logarithmic singularity, and the Dirac function. For the Dirac function, the  result of the integral is identically equal to zero because  $\mathbf{x}_0 \not\in S^{+}$. For the logarithmic singularity, we know it is integrable, which renders the limit zero. Thus, the limit
of the second term on the right hand side of (\ref{fa24a}) is zero. Then, (\ref{fa24a}) is simplified as follows:
\newcommand\tbbint{{-\mkern -16mu\int}}
\newcommand\tbint{{\mathchar '26\mkern -14mu\int}}
\newcommand\dbbint{{-\mkern -19mu\int}}
\newcommand\dbint{{\mathchar '26\mkern -18mu\int}}
\newcommand\bbint{{\mathchoice{\dbbint}{\tbbint}{\tbbint}{\tbbint}}}
\begin{sequation}\label{fa24b}
\begin{split}
& -u_k\left(\mathbf{x}_0\right) + \int_{\Omega} \mathbf{v}_k \left(\mathbf{x} - \mathbf{x}_0\right) \cdot \mathbf{f} \left(\mathbf{x}\right) {\rm d}V_{\mathbf{x}} \\
=& \bbint_{\partial \Omega} \mathbf{u} \left(\mathbf{x}\right) \cdot \mathbf{T}\left(\mathbf{v}_k\right) \left(\mathbf{x} - \mathbf{x}_0\right) - \mathbf{v}_k \left(\mathbf{x} - \mathbf{x}_0\right) \cdot \mathbf{T}\left(\mathbf{u}\right) \left(\mathbf{x}\right) {\rm d}V_{\mathbf{x}}
\end{split}
\end{sequation}
where $\bbint$ denotes the Cauchy principal value of the integral. Dividing $\partial \Omega$ into $\partial \Omega_n$ and $\partial \Omega_d$, (\ref{fa24b}) is rewritten as
\begin{sequation}\label{fa25}
\begin{split}
u_k\left(\mathbf{x}_0\right) =& - \bbint_{\partial \Omega_d} \overline{\mathbf{u}} \left(\mathbf{x}\right) \cdot \mathbf{T}\left(\mathbf{v}_k\right) \left(\mathbf{x} - \mathbf{x}_0\right) - \mathbf{v}_k \left(\mathbf{x} - \mathbf{x}_0\right) \cdot \mathbf{T}\left(\mathbf{u}\right) \left(\mathbf{x}\right) {\rm d}V_{\mathbf{x}} \\
& - \bbint_{\partial \Omega_n} \mathbf{u} \left(\mathbf{x}\right) \cdot \mathbf{T}\left(\mathbf{v}_k\right) \left(\mathbf{x} - \mathbf{x}_0\right) - \mathbf{v}_k \left(\mathbf{x} - \mathbf{x}_0\right) \cdot \overline{\mathbf{T}} \left(\mathbf{v}\right) \left(\mathbf{x}\right) {\rm d}V_{\mathbf{x}} \\
& + \int_{\Omega} \mathbf{v}_k \left(\mathbf{x} - \mathbf{x}_0\right) \cdot \mathbf{f} \left(\mathbf{x}\right) {\rm d}V_{\mathbf{x}}
\end{split}
\end{sequation}
(\ref{fa25}) is the PD static boundary integral equation, which has a different form from the one in the classical theory. However, the  reciprocal theorem in (\ref{fa20}) has the same form as the one in the classical theory. In addition, it is ensured that the nonlocal operator and the volume constraint model converge to the classical operator and the classical boundary~\cite{tb5,tb59} when the size of the horizon approaches zero, respectively. Therefore, the solution of (\ref{fa25}) will converge to the classical one when the size of the horizon approaches zero, which will also be demonstrated by the examples in Section \ref{atl4}.

\subsection{The dynamic boundary integral equation}\label{atl22}

For the classic theory, Graff generalizes the reciprocal theorem of elastostatics to elastodynamics~\cite{tb63}. Based on Graff's work, Cole develops the dynamic boundary element method in the time domain~\cite{tb64}. According to this idea, we derive the nonlocal dynamic boundary integral equation in the time domain, which is shown in Appendix \ref{ap4}. Although the time domain method is a straightforward one, it has the problem of time  accumulation error and a high computational cost. A commonly used way in the classical boundary element method is to introduce an integral transformation so that the dynamic problem is solved in the same way as the static one. Because each time point is calculated independently, the integral transformation method does not have time  accumulation error, and can be easily parallelized. Here, we develop the dynamic boundary integral equation of PD with the Laplace transformation. The development is achieved through two steps: One is the separation and derivation of the PD fundamental solution in the Laplace domain, and the other is the derivation of the nonlocal dynamic boundary integral equation in the Laplace domain.

\subsubsection{The PD fundamental solution in Laplace domain}\label{atl221}

The Laplace transformation of the displacement vector $\mathbf{u}\left(\mathbf{x},t\right)$ with respect to time $t$ is denoted by
\begin{sequation}\label{falp}
\tilde{\mathbf{u}}\left(\mathbf{x},s\right)\equiv \mathcal{L}\mathbf{u}\left(\mathbf{x},t\right)=\int_{0}^{\infty} \mathbf{u}\left(\mathbf{x},t\right) e^{-st} {\rm d}t
\end{sequation}
Then, the Fourier transformation of $\tilde{\mathbf{u}}(\mathbf{x},s)$ with respect to $\mathbf{x}$ is denoted by
\begin{sequation}\label{fafr}
\overline{\mathbf{u}}\left(\mathbf{k},s\right)\equiv \mathcal{F} \mathbf{u}\left(\mathbf{x},s\right)=\int \mathbf{u}\left(\mathbf{x},s\right) e^{-2\pi {\rm i} \left(\mathbf{x} \cdot \mathbf{k} \right)} {\rm d}V_{\mathbf{x}}
\end{sequation}
where ${\rm i}=\sqrt{-1}$. We carry out the Laplace transformation on both sides of (\ref{fa2}) to get
\begin{sequation}\label{fa35}
\rho s^2 \tilde{\mathbf{u}}\left(\mathbf{x},s\right) = \int_{\mathcal{H}_{\mathbf{x}}} \mathbf{C}\left(\boldsymbol{\xi}\right) \cdot \left[\tilde{\mathbf{u}}\left(\mathbf{x} + \boldsymbol{\xi},s\right) - \tilde{\mathbf{u}}\left(\mathbf{x},s\right)\right] {\rm d}V_{\boldsymbol{\xi}} + \tilde{\mathbf{F}}\left(\mathbf{x},s\right)
\end{sequation}
where  $\tilde{\mathbf{F}}\left(\mathbf{x},s\right) = \tilde{\mathbf{b}}\left(\mathbf{x},s\right) + \rho s \mathbf{u} \left(\mathbf{x},0\right) + \rho \dot{\mathbf{u}} \left(\mathbf{x},0\right)$. $\tilde{\mathbf{F}}\left(\mathbf{x},s\right)$ is equivalent to the body force density in the Laplace domain. If we let $\tilde{\mathbf{F}}\left(\mathbf{x},s\right) = \delta \left(\mathbf{x}\right) \mathbf{e}$, then we can get the fundamental solution in the Laplace domain when a unit body force in the direction $\mathbf{e}$ is applied in the infinitely extended material by solving (\ref{fa35}). Here, $\mathbf{e}$ is a coordinate base vector. Letting $\tilde{\mathbf{F}}\left(\mathbf{x},s\right) = \delta \left(\mathbf{x}\right) \mathbf{e}$, we introduce the Fourier transformation to change (\ref{fa35}) into an algebraic equation as follows:
\begin{sequation}\label{fa36}
\rho s^2 \overline{\mathbf{u}}_g\left(\mathbf{k},s\right) = - \mathbf{M} \left(\mathbf{k}\right) \cdot \overline{\mathbf{u}}_g\left(\mathbf{k},s\right) + \mathbf{e}
\end{sequation}
where $\overline{\mathbf{u}}_g\left(\mathbf{k},s\right)$ denotes the Green function. $\mathbf{k}$ is the wave vector in the Fourier domain. $\mathbf{M} \left(\mathbf{k}\right)$ has been presented in the literature (Silling~\cite{tb1}, Wang et al.~\cite{tb58}, Weckner et al.~\cite{tb65}). We can get the explicit form of the fundamental solution for (\ref{fa36}) by the inverse Fourier transformation
\begin{sequation}\label{fa37}
\tilde{\mathbf{u}}_g \left(\mathbf{x},s\right) = \mathcal{F}^{-1} \left[\left(\rho s^2 \mathbf{I} + \mathbf{M}\right)^{-1} \cdot \mathbf{e}\right]
\end{sequation}
Further, we can rewrite (\ref{fa37}) in the form of components as
\begin{sequation}\label{fa38}
\left(\tilde{{\rm u}}_g\right)_{ij}\left(\mathbf{x},s\right) = \mathcal{F}^{-1} \left[\frac{\delta_{ij}}{\rho s^2 + M_\perp \left(k\right)} + \frac{k_i k_j}{k^2} \left(\frac{1}{\rho s^2 + M_\Vert \left(k\right)} - \frac{1}{\rho s^2 + M_\perp \left(k\right)}\right)\right]
\end{sequation}
where $\left(\tilde{{\rm u}}_g\right)_{ij}$ is the $j$-th directional component of the Green function when a unit body force in the direction $\mathbf{e}_i$ is applied in the infinitely extended material. $\delta_{ij}$ is the Kronecker delta. $M_\perp \left(k\right)$ reflects the shear stiffness, and $M_\Vert \left(k\right)$ reflects the tensile stiffness. The detailed expressions of $M_\perp \left(k\right)$ and $M_\Vert \left(k\right)$ are given in the literature~\cite{tb58,tb65}.
The general form of the infinite Green function in the Laplace domain is
\begin{sequation}\label{fa39}
\left(\tilde{{\rm u}}_g\right)_{ij} \left(\mathbf{x},s\right) = \left(\tilde{{\rm u}}_g\right)_{A} \left(\mathbf{x},s\right) \delta_{ij} + \frac{x_i x_j}{x^2} \left(\tilde{{\rm u}}_g\right)_{B} \left(\mathbf{x},s\right)
\end{sequation}
For one (I), two (II), and three (III) dimensions, $\left(\tilde{{\rm u}}_g\right)_{A}$ and $\left(\tilde{{\rm u}}_g\right)_{B}$ are
\begin{sequation}\label{fa40}
\left\{
\begin{aligned}
& \left(\tilde{{\rm u}}_g\right)_{A}^{\uppercase\expandafter{\romannumeral1}} \left(\mathbf{x},s\right) = \frac{1}{\pi} \int_{0}^{+\infty} \frac{{\rm cos}\left(kx\right)}{\rho s^2 + M \left(k\right)} {\rm d}k \\
& \left(\tilde{{\rm u}}_g\right)_{B}^{\uppercase\expandafter{\romannumeral1}} \left(\mathbf{x},s\right) = 0
\end{aligned}
\right.
\qquad \qquad \qquad \qquad \qquad \qquad \qquad \qquad \qquad \qquad \qquad \qquad \ \
\end{sequation}
\begin{sequation}\label{fa41}
\left\{
\begin{aligned}
& \left(\tilde{{\rm u}}_g\right)_{A}^{\uppercase\expandafter{\romannumeral2}} \left(\mathbf{x},s\right) = \frac{1}{2 \pi} \int_{0}^{+\infty} \left[\frac{J_0 \left(kx\right)}{\rho s^2 + M_\perp^{\uppercase\expandafter{\romannumeral2}} \left(k\right)} + \frac{J_1 \left(kx\right)}{kx} \left(\frac{1}{\rho s^2 + M_\Vert^{\uppercase\expandafter{\romannumeral2}} \left(k\right)} - \frac{1}{\rho s^2 + M_\perp^{\uppercase\expandafter{\romannumeral2}} \left(k\right)}\right)\right] k {\rm d}k \\
& \left(\tilde{{\rm u}}_g\right)_{B}^{\uppercase\expandafter{\romannumeral2}} \left(\mathbf{x},s\right) = \frac{1}{2 \pi} \int_{0}^{+\infty} \left[J_0 \left(kx\right) - 2 \frac{J_1 \left(kx\right)}{kx}\right] \left[\frac{1}{\rho s^2 + M_\Vert^{\uppercase\expandafter{\romannumeral2}} \left(k\right)} - \frac{1}{\rho s^2 + M_\perp^{\uppercase\expandafter{\romannumeral2}} \left(k\right)}\right] k {\rm d}k
\end{aligned}
\right.
\qquad \qquad \ \ \
\end{sequation}
\begin{sequation}\label{fa42}
\left\{
\begin{aligned}
\left(\tilde{{\rm u}}_g\right)_{A}^{\uppercase\expandafter{\romannumeral3}} \left(\mathbf{x},s\right) =&\frac{1}{2 \pi^2} \int_{0}^{+\infty} \left[\frac{{\rm sin}\left(kx\right)}{\left(kx\right)^3} - \frac{{\rm cos}\left(kx\right)}{\left(kx\right)^2} \right] \left[\frac{1}{\rho s^2 + M_\Vert^{\uppercase\expandafter{\romannumeral3}} \left(k\right)} - \frac{1}{\rho s^2 + M_\perp^{\uppercase\expandafter{\romannumeral3}} \left(k\right)}\right] k^2 {\rm d}k \\
&+ \frac{1}{2 \pi^2} \int_{0}^{+\infty} \frac{{\rm sin}\left(kx\right)}{x} \frac{k}{\rho s^2 + M_\perp^{\uppercase\expandafter{\romannumeral3}} \left(k\right)} {\rm d}k\\
\left(\tilde{{\rm u}}_g\right)_{B}^{\uppercase\expandafter{\romannumeral3}} \left(\mathbf{x},s\right) =& \frac{1}{2 \pi^2} \int_{0}^{+\infty} \left[\frac{{\rm sin}\left(kx\right)}{\left(kx\right)} + \frac{3{\rm cos}\left(kx\right)}{\left(kx\right)^2} - \frac{3{\rm sin}\left(kx\right)}{\left(kx\right)^3}\right] \left[\frac{1}{\rho s^2 + M_\Vert^{\uppercase\expandafter{\romannumeral3}} \left(k\right)} - \frac{1}{\rho s^2 + M_\perp^{\uppercase\expandafter{\romannumeral3}} \left(k\right)}\right] k^2 {\rm d}k
\end{aligned}
\right.
\end{sequation}
where $J_0 \left(x\right)$ and $J_1 \left(x\right)$ are the first- and second-order Bessel functions, respectively.  If the inertia effect disappears, that is, $\rho = 0$, the  above formulas reduce to  the static Green functions~\cite{tb58}.

We note that the integrals in $\left(\tilde{{\rm u}}_g\right)_{A}^{\uppercase\expandafter{\romannumeral1}} \left(\mathbf{x},s\right)$,
$\left(\tilde{{\rm u}}_g\right)_{A}^{\uppercase\expandafter{\romannumeral2}} \left(\mathbf{x},s\right)$ and
$\left(\tilde{{\rm u}}_g\right)_{A}^{\uppercase\expandafter{\romannumeral3}} \left(\mathbf{x},s\right)$ are divergent because the integrands do not approach zero when $k\rightarrow +\infty$. This type of improper integrals has been examined in the literature~\cite{tb58,tb65};
they can be expressed as the parts with the Dirac functions, and convergent integrals
\begin{sequation}\label{fa44}
\left(\tilde{{\rm u}}_g\right)_{A}^{\uppercase\expandafter{\romannumeral1}} \left(\mathbf{x},s\right) = \frac{\delta \left(\mathbf{x}\right)}{\rho s^2 + M \left(\infty\right)} + \frac{1}{\pi} \int_{0}^{+\infty} \left[\frac{{\rm cos}\left(kx\right)}{\rho s^2 + M \left(k\right)} - \frac{{\rm cos}\left(kx\right)}{\rho s^2 + M \left(\infty\right)}\right] {\rm d}k
\end{sequation}
\begin{sequation}\label{fa45}
\begin{aligned}
\left(\tilde{{\rm u}}_g\right)_{A}^{\uppercase\expandafter{\romannumeral2}} \left(\mathbf{x},s\right) =& \frac{\delta^{\uppercase\expandafter{\romannumeral2}} \left(\mathbf{x}\right)}{\rho s^2 + M^{\uppercase\expandafter{\romannumeral2}} \left(\infty\right)} + \frac{1}{2 \pi} \int_{0}^{+\infty} J_0 \left(kx\right)\left[\frac{1}{\rho s^2 + M_\perp^{\uppercase\expandafter{\romannumeral2}} \left(k\right)} - \frac{1}{\rho s^2 + M^{\uppercase\expandafter{\romannumeral2}} \left(\infty\right)}\right] k {\rm d}k \\
& + \frac{1}{2 \pi} \int_{0}^{+\infty} \frac{J_1 \left(kx\right)}{x} \left[\frac{1}{\rho s^2 + M_\Vert^{\uppercase\expandafter{\romannumeral2}} \left(k\right)} - \frac{1}{\rho s^2 + M_\perp^{\uppercase\expandafter{\romannumeral2}} \left(k\right)}\right] {\rm d}k
\end{aligned}
\end{sequation}
\begin{sequation}\label{fa46}
\begin{aligned}
\left(\tilde{{\rm u}}_g\right)_{A}^{\uppercase\expandafter{\romannumeral3}} \left(\mathbf{x},s\right) =& \frac{\delta^{\uppercase\expandafter{\romannumeral3}} \left(\mathbf{x}\right)}{\rho s^2 + M^{\uppercase\expandafter{\romannumeral3}} \left(\infty\right)} + \frac{1}{2 \pi^2} \int_{0}^{+\infty} \frac{{\rm sin}\left(kx\right)}{x} \left[\frac{k}{\rho s^2 + M_\perp^{\uppercase\expandafter{\romannumeral3}} \left(k\right)} - \frac{k}{\rho s^2 + M^{\uppercase\expandafter{\romannumeral3}} \left(\infty\right)}\right] {\rm d}k\\
& + \frac{1}{2 \pi^2} \int_{0}^{+\infty} \left[\frac{{\rm sin}\left(kx\right)}{\left(kx\right)^3} - \frac{{\rm cos}\left(kx\right)}{\left(kx\right)^2} \right] \left[\frac{1}{\rho s^2 + M_\Vert^{\uppercase\expandafter{\romannumeral3}} \left(k\right)} - \frac{1}{\rho s^2 + M_\perp^{\uppercase\expandafter{\romannumeral3}} \left(k\right)}\right] k^2 {\rm d}k \\
&+ \frac{1}{2 \pi^2} \int_{0}^{+\infty} \frac{{\rm sin}\left(kx\right)}{x} \frac{k}{\rho s^2 + M_\perp^{\uppercase\expandafter{\romannumeral3}} \left(k\right)} {\rm d}k
\end{aligned}
\end{sequation}
where $\delta \left(\mathbf{x}\right), \delta^{\uppercase\expandafter{\romannumeral2}} \left(\mathbf{x}\right), \delta^{\uppercase\expandafter{\romannumeral3}} \left(\mathbf{x}\right)$ are one-, two- and three-dimensional Dirac functions, respectively.

\subsubsection{The nonlocal dynamic boundary integral equation in Laplace domain}\label{atl222}

Following the derivation in the classical theory, we here use the weighted residual method to derive the boundary integral equation instead of the reciprocal theorem. First of all, we need to give the description of the problem in the Laplace domain. (\ref{dpf1}) and (\ref{dpf2}) in Appendix \ref{ap4} are the descriptions of the problem in the time domain. By executing the Laplace transformation for (\ref{dpf1}) and (\ref{dpf2}), we get the description in the Laplace domain as follows:
\begin{sequation}\label{fa105}
\begin{cases}
\mathcal{D} \left(\boldsymbol{\Theta} : \mathcal{D}^{*} \left(\tilde{\mathbf{u}}\right) \right) \left(\mathbf{x},s\right) + \tilde{\mathbf{F}} \left(\mathbf{x},s\right) = \rho s^2 \tilde{\mathbf{u}} \left(\mathbf{x},s\right) \qquad & \text{for} \ \mathbf{x} \in \Omega \\
\tilde{\mathbf{u}} \left(\mathbf{x},s\right) = \tilde{\mathbf{g}}_{d} \left(\mathbf{x},s\right) \qquad & \text{for} \ \mathbf{x} \in \Omega_{\tau_d} \\
\mathcal{N} \left(\boldsymbol{\Theta} : \mathcal{D}^{*} \left(\tilde{\mathbf{u}}\right) \right) \left(\mathbf{x},s\right) = \tilde{\mathbf{g}}_{n} \left(\mathbf{x},s\right) \qquad & \text{for} \ \mathbf{x} \in \Omega_{\tau_n} \\
\tilde{\mathbf{F}}\left(\mathbf{x},s\right) = \tilde{\mathbf{b}}\left(\mathbf{x},s\right) + \rho s \mathbf{u} \left(\mathbf{x},t_0\right) + \rho \dot{\mathbf{u}} \left(\mathbf{x},t_0\right) \\
\mathbf{u} \left(\mathbf{x},t_0\right) = \mathbf{u}_0 \left(\mathbf{x}\right) \qquad \dot{\mathbf{u}} \left(\mathbf{x},t_0\right) = \dot{\mathbf{u}}_0 \left(\mathbf{x}\right)
\end{cases}\ \ \
\end{sequation}
where $\tilde{\mathbf{b}}$ is the Laplace transformation of the body force density in the time domain. $\mathbf{u}_0$ and $\dot{\mathbf{u}}_0$ are the initial conditions in the time domain. We still need to give two constraints on the boundary conditions, which is similar to what we did in the time domain. We perform Laplace transformations for the dynamic constraints (\ref{dpf10}) and (\ref{dpf11}) in Appendix \ref{ap4} to get the ones in the Laplace domain as follows:
\begin{sequation}\label{fa106}
\begin{cases}
\tilde{\mathbf{u}} \left(\mathbf{x},s\right) = \mathbf{D} \left(\mathbf{x},s\right) \quad & \text{for} \  \mathbf{x} \in \partial \Omega \\
\mathbf{T} \left(\tilde{\mathbf{u}}\right) \left(\mathbf{x},s\right) = \mathbf{P} \left(\mathbf{x},s\right) \quad & \text{for} \  \mathbf{x} \in \partial \Omega
\ \ \
\end{cases}
\end{sequation}
\begin{sequation}\label{fa107}
\int_{\Omega_\tau} \tilde{\mathbf{v}}_p \left(\mathbf{x},s'\right) \cdot \mathcal{N} \left(\boldsymbol{\Theta} : \mathcal{D}^{*} \left(\tilde{\mathbf{u}}\right)\right) \left(\mathbf{x},s\right) {\rm d}V_{\mathbf{x}} = \int_{\partial \Omega} \tilde{\mathbf{v}}_p \left(\mathbf{x},s'\right) \cdot \mathbf{T} \left(\tilde{\mathbf{u}}\right) \left(\mathbf{x},s\right) {\rm d}V_{\mathbf{x}}
\end{sequation}
$\tilde{\mathbf{v}}_p$ is still a possible deformed state and meets the following condition:
\begin{sequation}\label{fa108}
\begin{cases}
\tilde{\mathbf{v}}_p \left(\mathbf{x},s\right) = \tilde{\mathbf{g}}_{d_{p}} \left(\mathbf{x},s\right) \qquad & \text{for} \ \mathbf{x} \in \Omega_{\tau_d}\\
\tilde{\mathbf{v}}_p \left(\mathbf{x},s\right) = \tilde{\mathbf{u}}_p \left(\mathbf{x},s\right) \qquad & \text{for} \ \mathbf{x} \in \partial \Omega_d
\end{cases}
\end{sequation}
where $\tilde{\mathbf{g}}_{d_{p}}, \tilde{\mathbf{u}}_p$ are the given displacement constraints on the nonlocal boundary and on the local boundary, respectively.  We rewrite the first formulas of (\ref{fa105}) and (\ref{fa106}) in the form of weighted residual as follows:
\begin{sequation}\label{fa109a}
\int_{\Omega} \left(-\rho s^2 \tilde{\mathbf{u}} \left(\mathbf{x},s\right) + \mathcal{D} \left(\boldsymbol{\Theta} : \mathcal{D}^{*} \left(\tilde{\mathbf{u}}\right) \right) \left(\mathbf{x},s\right) + \tilde{\mathbf{F}} \left(\mathbf{x},s\right)\right) \cdot \tilde{\mathbf{u}}_g \left(\mathbf{x},s\right) {\rm d}V_{\mathbf{x}} = 0
\end{sequation}
\begin{sequation}\label{fa109b}
\int_{\partial \Omega_n} \left(\mathbf{P} \left(\mathbf{x},s\right) - \mathbf{T} \left(\tilde{\mathbf{u}}\right)
\left(\mathbf{x},s\right)\right) \cdot \tilde{\mathbf{u}}_g \left(\mathbf{x},s\right) {\rm d}V_{\mathbf{x}} - \int_{\partial \Omega_d} \left(\mathbf{D} \left(\mathbf{x},s\right) - \tilde{\mathbf{u}} \left(\mathbf{x},s\right)\right) \cdot \mathbf{T} \left(\tilde{\mathbf{u}}_g\right) \left(\mathbf{x},s\right) {\rm d}V_{\mathbf{x}} = 0
\end{sequation}
where  $\tilde{\mathbf{u}}_g \left(\mathbf{x},s\right)$ is the fundamental solution in the Laplace domain for the infinite body; $\mathbf{T} \left(\tilde{\mathbf{u}}_g\right) \left(\mathbf{x},s\right)$ is the surface force that corresponds to the fundamental solution $\tilde{\mathbf{u}}_g \left(\mathbf{x},s\right)$. In the following part, $\mathbf{T} \left(\tilde{\mathbf{u}}\right) \left(\mathbf{x},s\right)$ and $\mathbf{T} \left(\tilde{\mathbf{u}}_g\right) \left(\mathbf{x},s\right)$ are simply denoted as $\tilde{\mathbf{T}} \left(\mathbf{x},s\right)$ and $\tilde{\mathbf{T}}_g \left(\mathbf{x},s\right)$, respectively. Merging (\ref{fa109a}) and (\ref{fa109b}), we get
\begin{sequation}\label{fa109}
\begin{aligned}
& \int_{\Omega} \left(-\rho s^2 \tilde{\mathbf{u}} \left(\mathbf{x},s\right) + \mathcal{D} \left(\boldsymbol{\Theta} : \mathcal{D}^{*} \left(\tilde{\mathbf{u}}\right) \right) \left(\mathbf{x},s\right) + \tilde{\mathbf{F}} \left(\mathbf{x},s\right)\right) \cdot \tilde{\mathbf{u}}_g \left(\mathbf{x},s\right) {\rm d}V_{\mathbf{x}} \\
=& \int_{\partial \Omega_n} \left(\mathbf{P} \left(\mathbf{x},s\right) - \tilde{\mathbf{T}} \left(\mathbf{x},s\right)\right) \cdot \tilde{\mathbf{u}}_g \left(\mathbf{x},s\right) {\rm d}V_{\mathbf{x}} - \int_{\partial \Omega_d} \left(\mathbf{D} \left(\mathbf{x},s\right) - \tilde{\mathbf{u}} \left(\mathbf{x},s\right)\right) \cdot \tilde{\mathbf{T}}_g \left(\mathbf{x},s\right) {\rm d}V_{\mathbf{x}}
\end{aligned}
\end{sequation}
Next, we construct two static states $\tilde{\mathbf{u}} \left(\mathbf{x},s\right)$ and $\tilde{\mathbf{u}}_g \left(\mathbf{x},s\right)$. According to (\ref{fa8}), (\ref{fa13a}) and (\ref{fa13ab}), $\tilde{\mathbf{u}} \left(\mathbf{x},s\right)$ and $\tilde{\mathbf{u}}_g \left(\mathbf{x},s\right)$ satisfy the following two sets of governing equations and boundary/initial conditions:
\begin{sequation}\label{fa110}
\begin{cases}
\mathcal{D} \left(\boldsymbol{\Theta} : \mathcal{D}^{*} \left(\tilde{\mathbf{u}}\right) \right) \left(\mathbf{x},s\right) + \mathbf{b} \left(\tilde{\mathbf{u}}\right) \left(\mathbf{x},s\right) = \mathbf{0} \qquad & \text{for} \ \mathbf{x} \in \Omega \\
\tilde{\mathbf{u}} \left(\mathbf{x},s\right) = \tilde{\mathbf{u}} \left(\mathbf{x},s\right) \qquad & \text{for} \ \mathbf{x} \in \Omega_{\tau_d} \\
\mathcal{N} \left(\boldsymbol{\Theta} : \mathcal{D}^{*} \left(\tilde{\mathbf{u}}\right) \right) \left(\mathbf{x},s\right) = \phi \left(\tilde{\mathbf{u}}\right) \left(\mathbf{x},s\right) \qquad & \text{for} \ \mathbf{x} \in \Omega_{\tau_n}
\end{cases}
\end{sequation}
\begin{sequation}\label{fa111}
\begin{cases}
\mathcal{D} \left(\boldsymbol{\Theta} : \mathcal{D}^{*} \left(\tilde{\mathbf{u}}_g\right) \right) \left(\mathbf{x},s\right) + \mathbf{b} \left(\tilde{\mathbf{u}}_g\right) \left(\mathbf{x},s\right) = \mathbf{0} \qquad & \text{for} \ \mathbf{x} \in \Omega \\
\tilde{\mathbf{u}}_g \left(\mathbf{x},s\right) = \tilde{\mathbf{u}}_g \left(\mathbf{x},s\right) \qquad & \text{for} \ \mathbf{x} \in \Omega_{\tau_d} \\
\mathcal{N} \left(\boldsymbol{\Theta} : \mathcal{D}^{*} \left(\tilde{\mathbf{u}}_g\right) \right) \left(\mathbf{x},s\right) = \phi \left(\tilde{\mathbf{u}}_g\right) \left(\mathbf{x},s\right) \qquad & \text{for} \ \mathbf{x} \in \Omega_{\tau_n}
\end{cases}
\end{sequation}
where $\mathbf{b} \left(\tilde{\mathbf{u}}\right)$ and $\mathbf{b} \left(\tilde{\mathbf{u}}_g\right)$ are the body force densities that correspond to the displacement $\tilde{\mathbf{u}}$ and the displacement $\tilde{\mathbf{u}}_g$, respectively. $\phi \left(\tilde{\mathbf{u}}\right)$ and $\phi \left(\tilde{\mathbf{u}}_g\right)$ are the boundary constraints that correspond to the displacement $\tilde{\mathbf{u}}$ and the displacement $\tilde{\mathbf{u}}_g$, respectively. These two states satisfy the static reciprocal theorem (\ref{fa20})
\begin{sequation}\label{fa112}
\begin{aligned}
& \int_{\Omega} \left(\mathcal{D} \left(\boldsymbol{\Theta} : \mathcal{D}^{*} \left(\tilde{\mathbf{u}}_g\right)\right) \left(\mathbf{x},s\right) \cdot \tilde{\mathbf{u}} \left(\mathbf{x},s\right) - \mathcal{D} \left(\boldsymbol{\Theta} : \mathcal{D}^{*} \left(\tilde{\mathbf{u}}\right)\right) \left(\mathbf{x},s\right) \cdot \tilde{\mathbf{u}}_g \left(\mathbf{x},s\right)\right) {\rm d}V_{\mathbf{x}} \\
=& \int_{\partial \Omega} \left(\tilde{\mathbf{T}}_g \left(\mathbf{x},s\right) \cdot \tilde{\mathbf{u}} \left(\mathbf{x},s\right) - \tilde{\mathbf{T}} \left(\mathbf{x},s\right) \cdot \tilde{\mathbf{u}}_g \left(\mathbf{x},s\right)\right) {\rm d}V_{\mathbf{x}}
\end{aligned}
\end{sequation}
Substituting (\ref{fa112}) into (\ref{fa109}), we obtain
\begin{sequation}\label{fa114}
\begin{aligned}
& \int_{\Omega} \left(\mathcal{D} \left(\boldsymbol{\Theta} : \mathcal{D}^{*} \left(\tilde{\mathbf{u}}_g\right)\right) \left(\mathbf{x},s\right) - \rho s^2 \tilde{\mathbf{u}}_g \left(\mathbf{x},s\right)\right) \cdot \tilde{\mathbf{u}} \left(\mathbf{x},s\right) {\rm d}V_{\mathbf{x}} + \int_{\Omega} \tilde{\mathbf{F}} \left(\mathbf{x},s\right) \cdot \tilde{\mathbf{u}}_g \left(\mathbf{x},s\right) {\rm d}V_{\mathbf{x}} \\
=& \int_{\partial \Omega} \left(\tilde{\mathbf{T}}_g \left(\mathbf{x},s\right) \cdot \tilde{\mathbf{u}} \left(\mathbf{x},s\right) - \tilde{\mathbf{T}} \left(\mathbf{x},s\right) \cdot \tilde{\mathbf{u}}_g \left(\mathbf{x},s\right)\right) {\rm d}V_{\mathbf{x}} + \int_{\partial \Omega_n} \left(\mathbf{P} \left(\mathbf{x},s\right) - \tilde{\mathbf{T}} \left(\mathbf{x},s\right)\right) \cdot \tilde{\mathbf{u}}_g \left(\mathbf{x},s\right) {\rm d}V_{\mathbf{x}} \\
& - \int_{\partial \Omega_d} \left(\mathbf{D} \left(\mathbf{x},s\right) - \tilde{\mathbf{u}} \left(\mathbf{x},s\right)\right) \cdot \tilde{\mathbf{T}}_g \left(\mathbf{x},s\right) {\rm d}V_{\mathbf{x}}
\end{aligned}
\end{sequation}
The equilibrium equation corresponding to the fundamental solution $\tilde{\mathbf{u}}_g^k$ in the Laplace domain when a unit body force in the direction $\mathbf{e}_k$ is applied in the infinitely extended material is
\begin{sequation}\label{fa115}
- \mathcal{D} \left(\boldsymbol{\Theta} : \mathcal{D}^{*} \left(\tilde{\mathbf{u}}_g^k\right)\right) \left(\mathbf{x} - \mathbf{x}_0,s\right) + \rho s^2 \tilde{\mathbf{u}}_g^k \left(\mathbf{x} - \mathbf{x}_0,s\right) = \delta \left(\mathbf{x} - \mathbf{x}_0\right) \mathbf{e}_k
\end{sequation}
Substituting (\ref{fa115}) into (\ref{fa114}) and dividing the integral domain $\partial \Omega$ in the first integral on the right hand side of (\ref{fa114}) into two parts $\partial \Omega_d$ and $\partial \Omega_n$, we get
\begin{sequation}\label{fa116a}
\begin{aligned}
& \int_{\Omega} \delta \left(\mathbf{x} - \mathbf{x}_0\right) \tilde{\mathbf{u}}_k \left(\mathbf{x},s\right) {\rm d}V_{\mathbf{x}} - \int_{\Omega} \tilde{\mathbf{F}} \left(\mathbf{x},s\right) \cdot \tilde{\mathbf{u}}_g^k \left(\mathbf{x} - \mathbf{x}_0,s\right) {\rm d}V_{\mathbf{x}} \\
=& \int_{\partial \Omega_d} \tilde{\mathbf{T}} \left(\mathbf{x},s\right) \cdot \tilde{\mathbf{u}}_g^k \left(\mathbf{x} - \mathbf{x}_0,s\right) {\rm d}V_{\mathbf{x}} - \int_{\partial \Omega_d} \tilde{\mathbf{u}} \left(\mathbf{x},s\right) \cdot \tilde{\mathbf{T}}_g^k \left(\mathbf{x} - \mathbf{x}_0,s\right) {\rm d}V_{\mathbf{x}} \\
& + \int_{\partial \Omega_n} \tilde{\mathbf{T}} \left(\mathbf{x},s\right) \cdot \tilde{\mathbf{u}}_g^k \left(\mathbf{x} - \mathbf{x}_0,s\right) {\rm d}V_{\mathbf{x}} - \int_{\partial \Omega_n} \tilde{\mathbf{T}}_g^k \left(\mathbf{x} - \mathbf{x}_0,s\right) \cdot \tilde{\mathbf{u}} \left(\mathbf{x},s\right) {\rm d}V_{\mathbf{x}}
\end{aligned}
\end{sequation}
where $\mathbf{T} \left(\tilde{\mathbf{u}}_g^k\right) \left(\mathbf{x}\right)$ is simply denoted as $\tilde{\mathbf{T}}_g^k \left(\mathbf{x}\right)$. Taking the limit process $\mathbf{x}_0 \to \partial \Omega_d \cup \partial \Omega_n$ for (\ref{fa116a}) similarly to (\ref{fa24}) - (\ref{fa24b}), we can get the boundary integral equation in the Laplace domain as follows:
\begin{sequation}\label{fa117}
\begin{aligned}
\tilde{{\rm u}}_k \left(\mathbf{x},s\right) =& \bbint_{\partial \Omega_d} \left(\tilde{\mathbf{T}} \left(\mathbf{x}',s\right) \cdot \tilde{\mathbf{u}}_g^k \left(\mathbf{x}'-\mathbf{x},s\right) - \tilde{\mathbf{T}}_g^k \left(\mathbf{x}'-\mathbf{x},s\right) \cdot \tilde{\mathbf{u}} \left(\mathbf{x}',s\right)\right) {\rm d}V_{\mathbf{x}'} \\
& + \bbint_{\partial \Omega_n} \left(\tilde{\mathbf{u}}_g^k \left(\mathbf{x}'-\mathbf{x},s\right) \cdot \tilde{\mathbf{T}} \left(\mathbf{x}',s\right) - \tilde{\mathbf{u}} \left(\mathbf{x}',s\right) \cdot \tilde{\mathbf{T}}_g^k \left(\mathbf{x}'-\mathbf{x},s\right)\right) {\rm d}V_{\mathbf{x}'} \\
& + \int_{\Omega} \tilde{\mathbf{F}} \left(\mathbf{x}',s\right) \cdot \tilde{\mathbf{u}}_g^k \left(\mathbf{x}'-\mathbf{x},s\right) {\rm d}V_{\mathbf{x}'} 
\end{aligned}
\end{sequation}

So far, we have derived the PD boundary integral equations (\ref{fa25}) and (\ref{fa117}) for the dynamic and static problems.
Eqs. (\ref{fa25}) and (\ref{fa117}) can be transformed into the following form for $\mathbf{x} \in \Omega \cup \partial \Omega$:
\begin{sequation}\label{faa26}
\begin{split}
{\rm u}_k \left(\mathbf{x}\right) =& - \int_{\partial \Omega_d} \overline{\mathbf{u}} \left(\mathbf{x}'\right) \cdot \mathbf{T}\left(\mathbf{v}_k\right) \left(\mathbf{x}' - \mathbf{x}\right) - \mathbf{v}_k \left(\mathbf{x}' - \mathbf{x}\right) \cdot \mathbf{T}\left(\mathbf{u}\right) \left(\mathbf{x}'\right) {\rm d}V_{\mathbf{x}'} \\
& - \int_{\partial \Omega_n} \mathbf{u} \left(\mathbf{x}'\right) \cdot \mathbf{T}\left(\mathbf{v}_k\right) \left(\mathbf{x}' - \mathbf{x}\right) - \mathbf{v}_k \left(\mathbf{x}' - \mathbf{x}\right) \cdot \overline{\mathbf{T}} \left(\mathbf{v}\right) \left(\mathbf{x}'\right) {\rm d}V_{\mathbf{x}'} \\
& + \int_{\Omega} \mathbf{v}_k \left(\mathbf{x}' - \mathbf{x}\right) \cdot \mathbf{f} \left(\mathbf{x}'\right) {\rm d}V_{\mathbf{x}'}
\end{split}
\end{sequation}
\begin{sequation}\label{faa117}
\begin{aligned}
\tilde{{\rm u}}_k \left(\mathbf{x},s\right) =& \int_{\partial \Omega_d} \left(\tilde{\mathbf{T}} \left(\mathbf{x}',s\right) \cdot \tilde{\mathbf{u}}_g^k \left(\mathbf{x}'-\mathbf{x},s\right) - \tilde{\mathbf{T}}_g^k \left(\mathbf{x}'-\mathbf{x},s\right) \cdot \tilde{\mathbf{u}} \left(\mathbf{x}',s\right)\right) {\rm d}V_{\mathbf{x}'} \\
& + \int_{\partial \Omega_n} \left(\tilde{\mathbf{u}}_g^k \left(\mathbf{x}'-\mathbf{x},s\right) \cdot \tilde{\mathbf{T}} \left(\mathbf{x}',s\right) - \tilde{\mathbf{u}} \left(\mathbf{x}',s\right) \cdot \tilde{\mathbf{T}}_g^k \left(\mathbf{x}'-\mathbf{x},s\right)\right) {\rm d}V_{\mathbf{x}'} \\
& + \int_{\Omega} \tilde{\mathbf{F}} \left(\mathbf{x}',s\right) \cdot \tilde{\mathbf{u}}_g^k \left(\mathbf{x}'-\mathbf{x},s\right) {\rm d}V_{\mathbf{x}'} 
\end{aligned}
\end{sequation}
Eqs. (\ref{faa26}) and (\ref{faa117}) are the PD Somigliana formulas for the static and the dynamic problems, respectively. We see that the PD Somigliana formulas (\ref{faa26}) and (\ref{faa117}) and the PD boundary integral equations (\ref{fa25}) and (\ref{fa117}) have the same form, which is due to the volume constraint model. Introducing the volume constraint boundary leads to equal status for the domain and boundary. Generally, for a given problem, the solution strategy is that we get the boundary quantities $\mathbf{u} \left(\mathbf{x}\right)$ and $\mathbf{T}\left(\mathbf{u}\right) \left(\mathbf{x}\right)$ from (\ref{fa25}), and the boundary quantities $\tilde{\mathbf{u}} \left(\mathbf{x},s\right)$ and $\tilde{\mathbf{T}} \left(\mathbf{x},s\right)$ from (\ref{fa117}) by the PD-BEM; and then, we can get the displacement distribution in the domain by means of the PD Somigliana formulas (\ref{faa26}) and (\ref{faa117}). The details will be presented in the following sections.

\section{Numerical framework}\label{atl3}

In this section, we will discuss the numerical treatment of the PD-BEM for  static and  dynamic problems.

\subsection{Numerical treatment for static problems}\label{atl35}

For static problems, the steps for the solution follow the classical boundary element method~\cite{tb60}. The difference is about the calculation of the integral in the discretized boundary element formulation.

\subsubsection{The discretized boundary element formulation}\label{atl32}

\begin{figure}[!htb]
	\centerline{\includegraphics[scale=0.3]{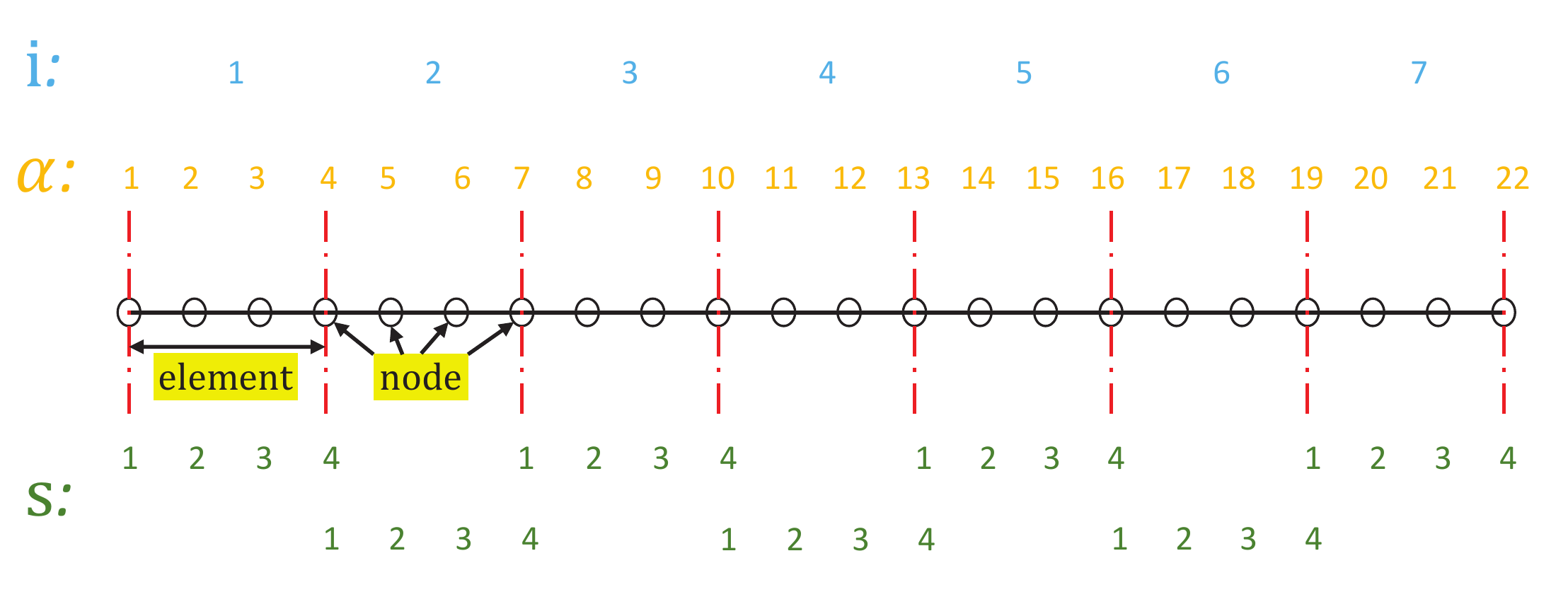}}	
	\caption{The relationship of numbers of elements and nodes.\label{sc8}}
\end{figure}

Firstly, we give the discrete form of grid interpolation as follows:
\begin{sequation}\label{fa28}
\left\{
\begin{aligned}
& {\rm u}_j \left(\mathbf{x}'\right) = \sum_{i=1}^{m} \sum_{s=1}^{N_e^i} \overline{N}_s^i \left(\mathbf{x}'\right) u_j \left(\mathbf{x}_s^i\right) \\
& {\rm T}_j \left(\mathbf{x}'\right) = \sum_{i=1}^{m} \sum_{s=1}^{N_e^i} \overline{N}_s^i \left(\mathbf{x}'\right) T_j \left(\mathbf{x}_s^i\right) \\
& \partial \Omega = \sum_{i=1}^{m} \partial \Omega_i
\end{aligned}
\right.
\end{sequation}
Here, the geometric boundary is divided into $m$ elements; $i$ is the element number; $N_e^i$ is the total number of nodes for element $i$; $\overline{N}_s^i$ is the shape function of node $s$ for element $i$; $s$ is the node number inside the element; $\mathbf{x}_s^i$ is the coordinate of node $s$ of element $i$. The surface force $\mathbf{T}\left(\mathbf{u}\right)$ is simply written as $\mathbf{T}$. $\mathbf{x}' \in \partial \Omega$. $u_{j}$ and $T_{j}$ are the components of $\mathbf{u}$ and $\mathbf{T}\left(\mathbf{u}\right)$, respectively. Secondly, based on (\ref{fa28}), we can discretize (\ref{fa25}) into the following form:
\begin{sequation}\label{fa29}
\begin{split}
& {\rm u}_k \left(\mathbf{x}_\alpha\right) + \sum_{i=1}^{m} \sum_{s=1}^{N_e^i} \sum_{j=1}^{dm} \bbint_{\partial \Omega_i} {\rm T}_{kj} \left(\mathbf{x}' - \mathbf{x}_\alpha\right) \overline{N}_s^i \left(\mathbf{x}'\right) {\rm u}_j \left(\mathbf{x}_s^i\right) {\rm d}V_{\mathbf{x}'} \\
=& \sum_{j=1}^{dm} \int_{\Omega} {\rm v}_{kj} \left(\mathbf{x}' - \mathbf{x}_\alpha\right) {\rm f}_j \left(\mathbf{x}'\right) {\rm d}V_{\mathbf{x}'} + \sum_{i=1}^{m} \sum_{s=1}^{N_e^i} \sum_{j=1}^{dm} \bbint_{\partial \Omega_i} {\rm v}_{kj} \left(\mathbf{x}' - \mathbf{x}_\alpha\right) \overline{N}_s^i \left(\mathbf{x}'\right) {\rm T}_j \left(\mathbf{x}_s^i\right) {\rm d}V_{\mathbf{x}'}
\end{split}
\end{sequation}
{where $dm$ is the dimension of space; ${\rm u}_{kj}$ is the $j$-th directional component of the  Green function when a unit body force in the direction $\mathbf{e}_k$ is applied in the infinitely extended material; ${\rm T}_{kj}$ is the $j$-th directional component of the surface force that corresponds to the  Green function; $\mathbf{x}_\alpha$ is the node coordinate. $\alpha$ denotes the node marking number for the boundary, and the total number of nodes is $n$. It is clear that there is a correspondence between the index $\alpha$ and the index pair $\left(i, s\right)$. The relationships between these numbers are shown in Figure \ref{sc8}.}

\begin{figure}[!htb]
	\centering
	\subfigure[]{
		\label{sc9b}
		\includegraphics[scale=0.20]{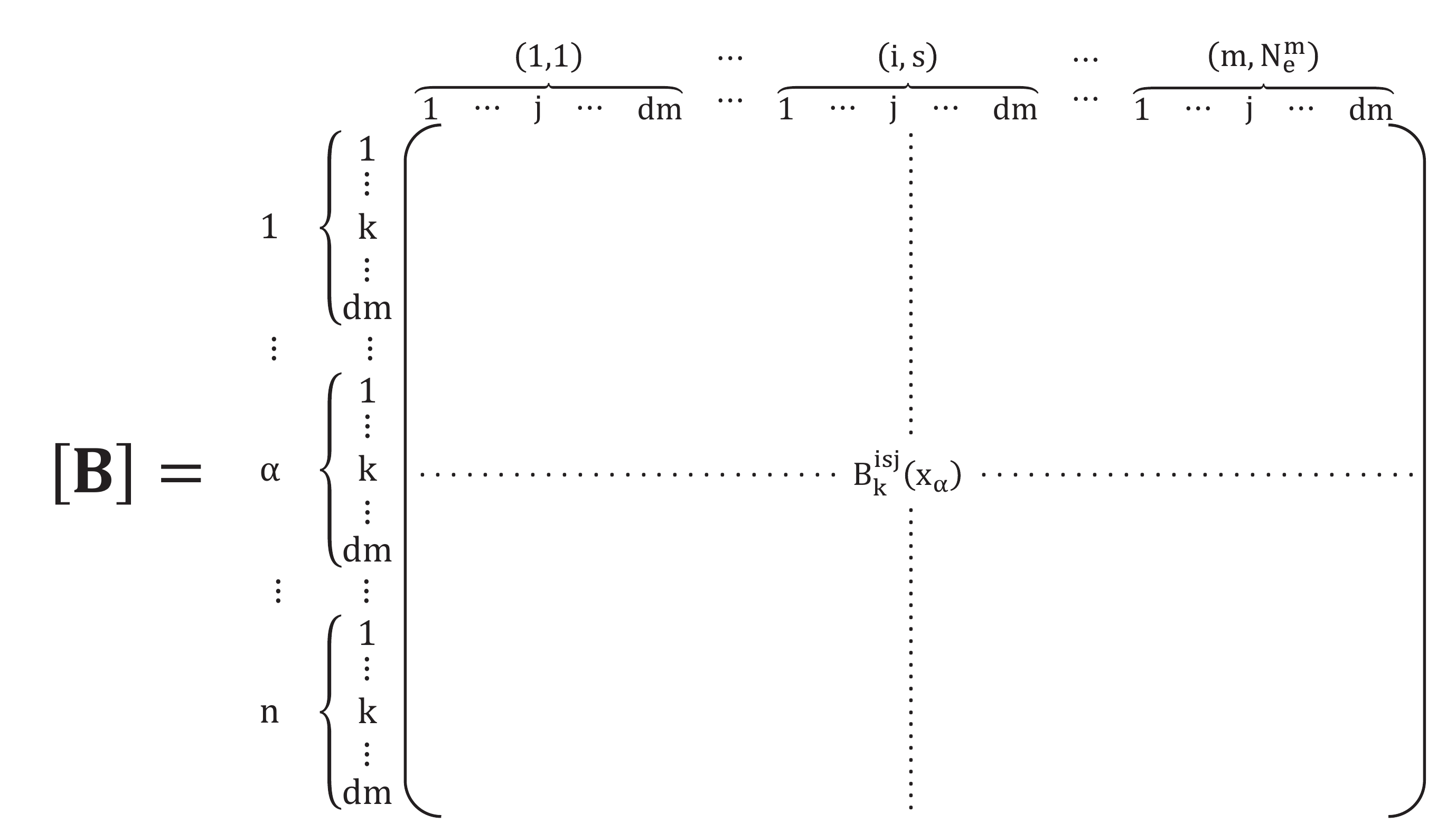}}
	\hspace{0.01in}
	\subfigure[]{
		\label{sc9c}
		\includegraphics[scale=0.20]{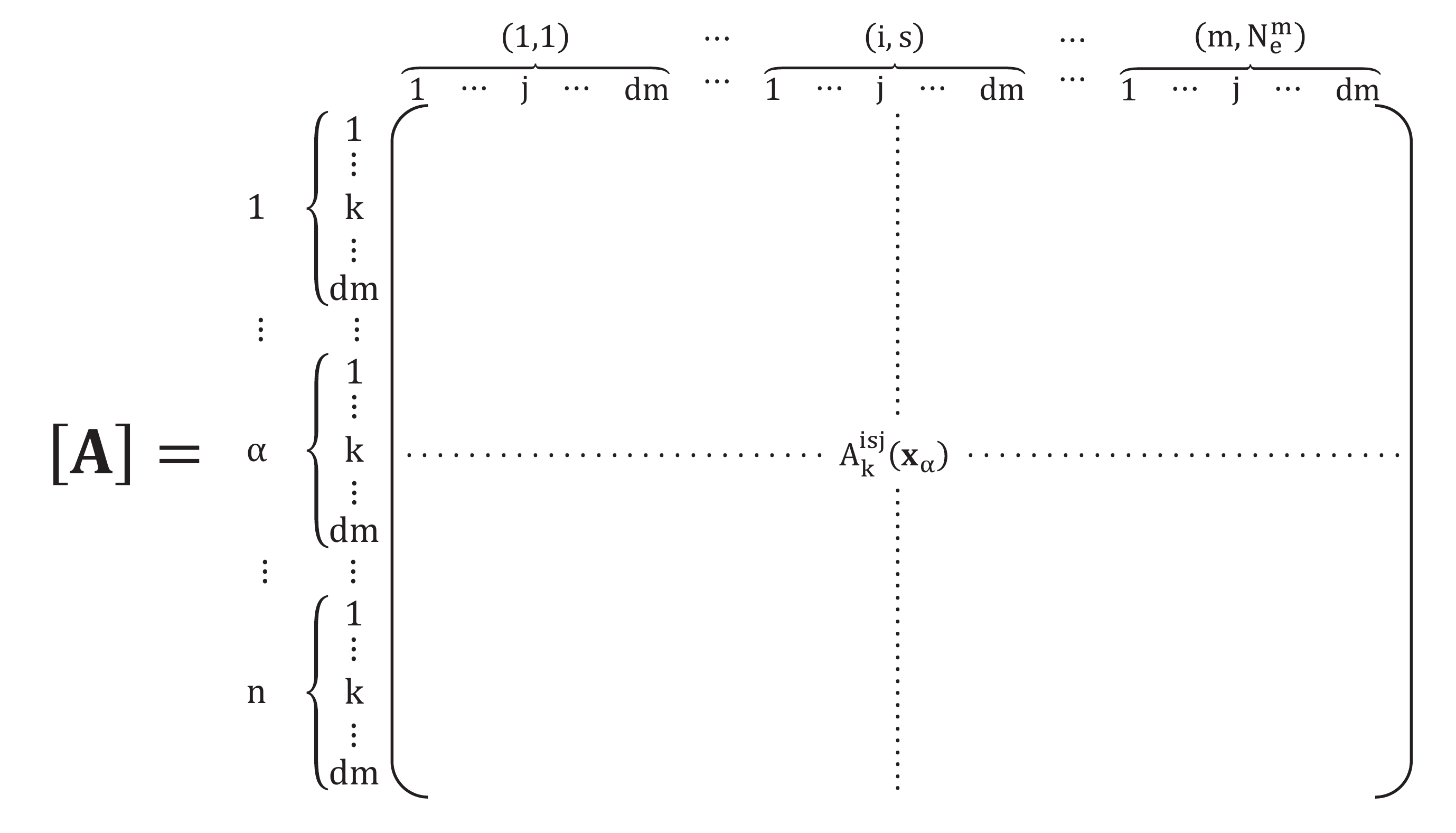}}
	\vfill
	\subfigure[]{
		\label{sc9d}
		\includegraphics[width=1.5in,height=1.20in]{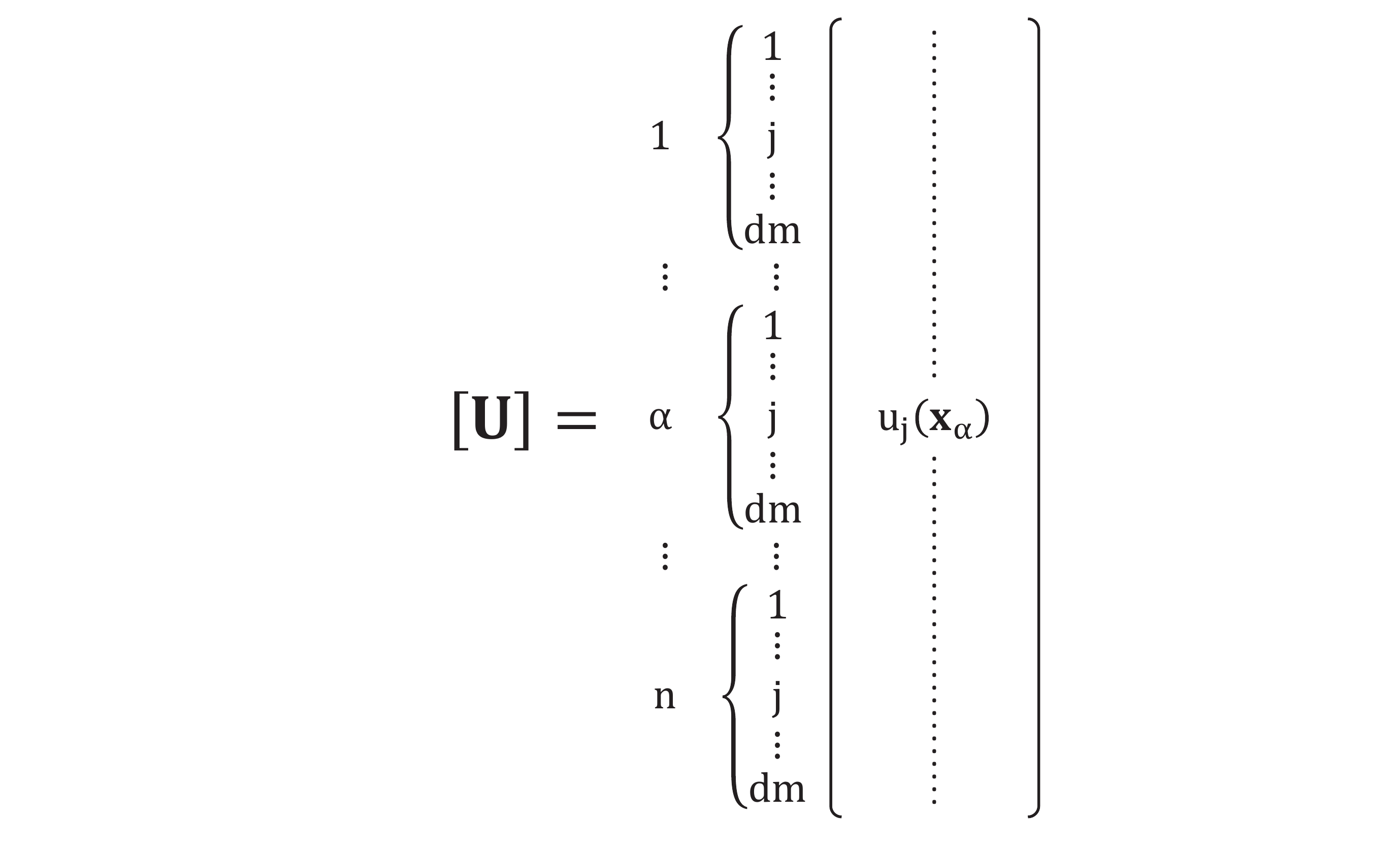}}
	\hspace{0.01in}
	\subfigure[]{
		\label{sc9e}
		\includegraphics[width=1.5in,height=1.20in]{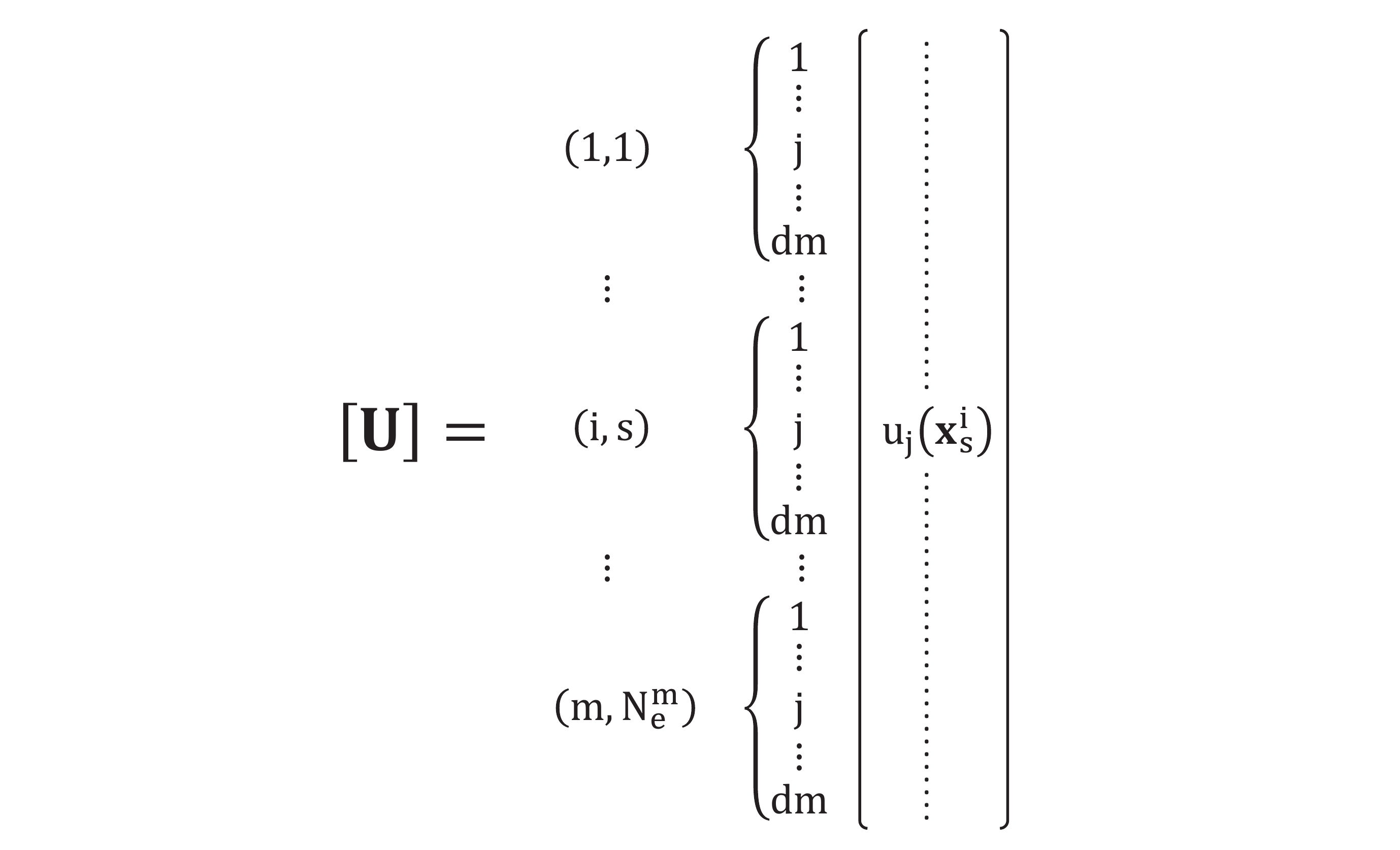}}	
	\vfill
	\subfigure[]{
		\label{sc9f}
		\includegraphics[width=1.5in,height=1.20in]{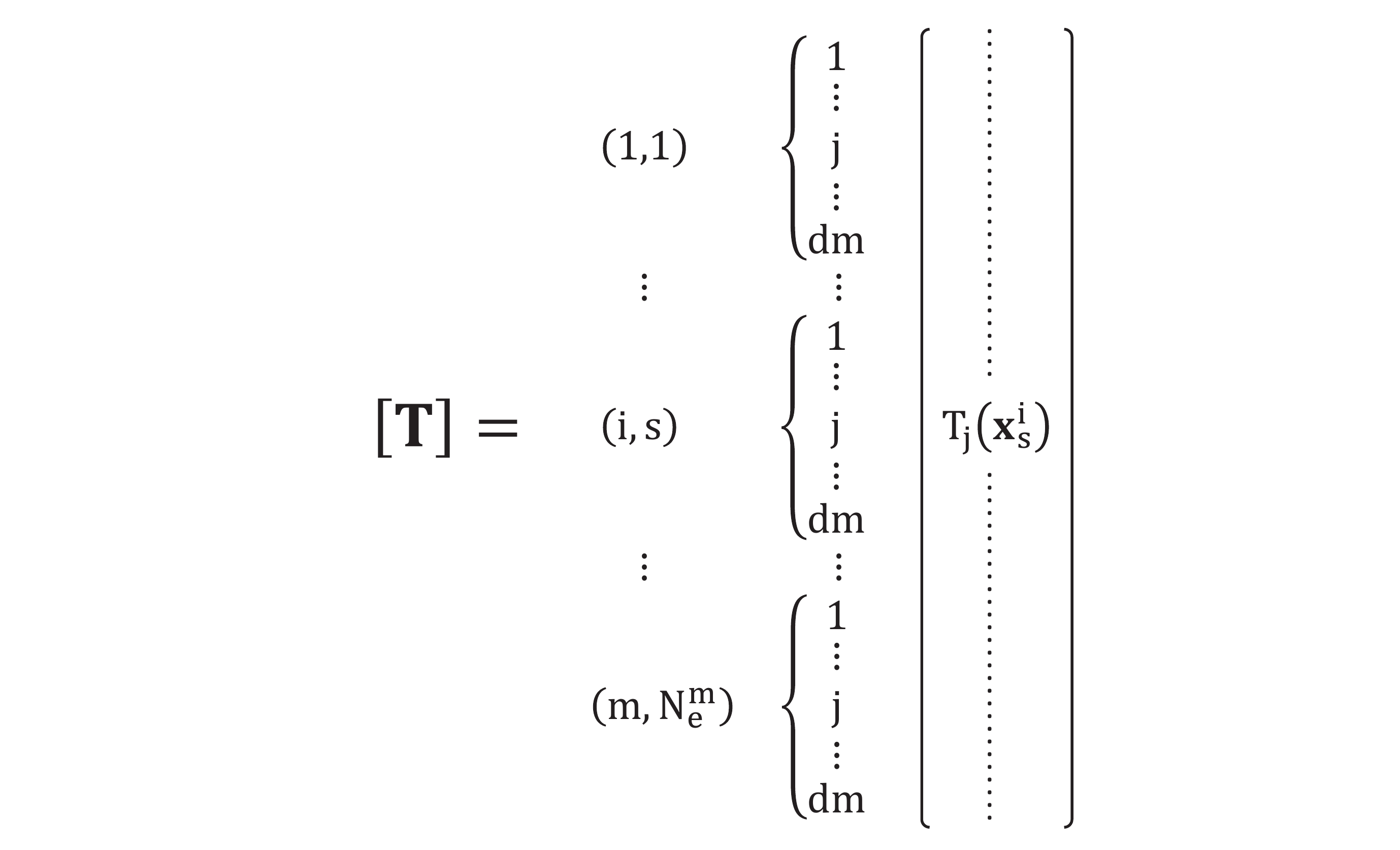}}	
	\hspace{0.01in}
	\subfigure[]{
		\label{sc9g}
		\includegraphics[width=1.5in,height=1.20in]{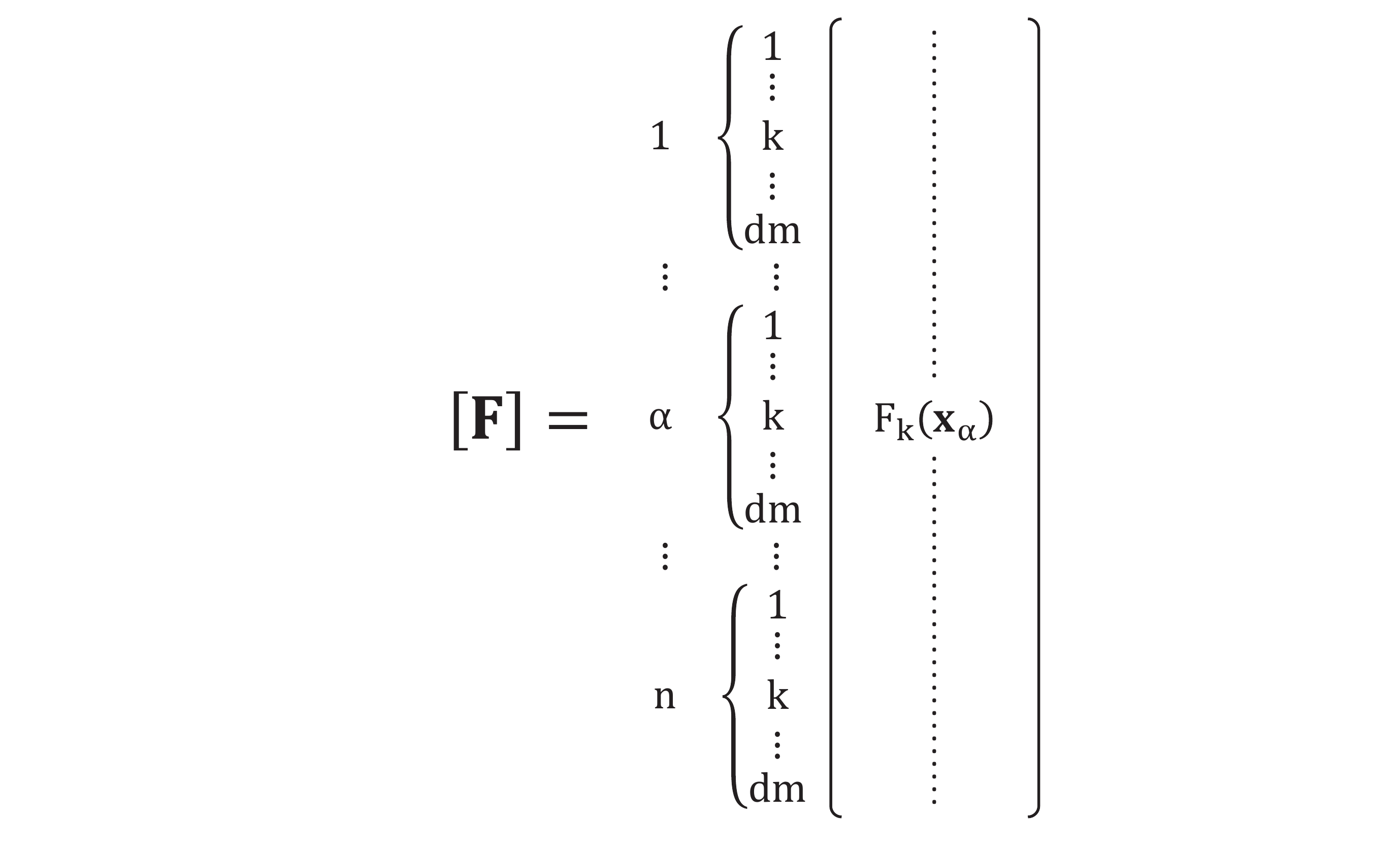}}	
	\caption{The locations of entries in the matrices in (\ref{fa31}) and (\ref{fa32}): ${\rm B}_{k}^{isj} \left(\mathbf{x}_\alpha\right)$(a); ${\rm A}_{k}^{isj} \left(\mathbf{x}_\alpha\right)$(b); ${\rm u}_j \left(\mathbf{x}_\alpha\right)$(c); ${\rm u}_j \left(\mathbf{x}_s^i\right)$(d); ${\rm T}_j \left(\mathbf{x}_s^i\right)$(e); ${\rm F}_k \left(\mathbf{x}_\alpha\right)$(f).}
	\label{sc9}	
\end{figure}

The infinite Green function can be expressed as~\cite{tb58}
\begin{sequation}\label{fa27}
{\rm u}_{ij} \left(\mathbf{x}\right) = {\rm u}_A \left(x\right) \delta_{ij} + \frac{{\rm x}_i {\rm x}_j}{x^2} {\rm u}_B \left(x\right)
\end{sequation}
where $\delta_{ij}$ is the Kronecker delta; $u_{ij}$ is the component $j$ of the Green function when a unit body force in the direction $\mathbf{e}_i$ is applied in the infinitely extended material; $u_A \left(x\right)$ and $u_B \left(x\right)$ are related to the dimension, and they are different for static problems and dynamic problems. We introduce the following notation to simplify (\ref{fa25}):
\begin{sequation}\label{fa30}
\left\{
\begin{aligned}
& {\rm A}_{k}^{isj} \left(\mathbf{x}_\alpha\right) \equiv \bbint_{\partial \Omega_i} {\rm T}_{kj} \left(\mathbf{x}' - \mathbf{x}_\alpha\right) \overline{N}_s^i \left(\mathbf{x}'\right) {\rm d}V_{\mathbf{x}'} \\
& {\rm B}_{k}^{isj} \left(\mathbf{x}_\alpha\right) \equiv \bbint_{\partial \Omega_i} {\rm u}_{kj} \left(\mathbf{x}' - \mathbf{x}_\alpha\right) \overline{N}_s^i \left(\mathbf{x}'\right) {\rm d}V_{\mathbf{x}'} \\
& {\rm F}_k \left(\mathbf{x}_\alpha\right) \equiv \sum_{j=1}^{dm} \int_{\Omega} {\rm u}_{kj} \left(\mathbf{x}' - \mathbf{x}_\alpha\right) {\rm f}_j \left(\mathbf{x}'\right) {\rm d}V_{\mathbf{x}'}
\end{aligned}
\right.
\end{sequation}
(\ref{fa29}) is simplified as follows:
\begin{sequation}\label{fa31}
{\rm u}_k \left(\mathbf{x}_\alpha\right) + \sum_{i=1}^{m} \sum_{s=1}^{N_e^i} \sum_{j=1}^{dm} {\rm A}_{k}^{isj} \left(\mathbf{x}_\alpha\right) {\rm u}_j \left(\mathbf{x}_s^i\right) = \sum_{i=1}^{m} \sum_{s=1}^{N_e^i} \sum_{j=1}^{dm} {\rm B}_{k}^{isj} \left(\mathbf{x}_\alpha\right) {\rm T}_j \left(\mathbf{x}_s^i\right) + {\rm F}_k \left(\mathbf{x}_\alpha\right),
\end{sequation}
which can be arranged into a matrix form
\begin{sequation}\label{fa32}
\left[\mathbf{I}\right] \left[\mathbf{U}\right] + \left[\mathbf{A}\right] \left[\mathbf{U}\right] = \left[\mathbf{B}\right] \left[\mathbf{T}\right] + \left[\mathbf{F}\right]
\end{sequation}
where $\left[\mathbf{U}\right]$ is a vector that is made up of the displacement components of all nodes, and $\left[\mathbf{T}\right]$ is a vector that is made up of the surface force components of all nodes. $\mathbf{x}_s^i \in \left\{\mathbf{x}_\alpha | \alpha = 1, \ldots, n\right\}$. The correspondence between the index $\alpha$ and the index pair $\left(i, s\right)$ is shown in Figure \ref{sc8}. $\alpha_s^i$ represents  $\alpha$ at $\left(i, s\right)$. ${\rm u}_j \left(\mathbf{x}_\alpha\right)$ and ${\rm u}_j \left(\mathbf{x}_s^i\right)$ correspond to $\left[\mathbf{U}_q\right]$, and $q = \left(\alpha-1\right)*dm + j, q = \left(\alpha_s^i-1\right)*dm + j$; ${\rm T}_j \left(\mathbf{x}_s^i\right)$ corresponds to $\left[\mathbf{T}_q\right]$, and $q = \left(\alpha_s^i-1\right)*dm + j$; ${\rm F}_k \left(\mathbf{x}_\alpha\right)$ corresponds to $\left[\mathbf{F}_q\right]$, and $q = \left(\alpha-1\right)*dm + k$; ${\rm A}_{k}^{isj} \left(\mathbf{x}_\alpha\right)$ corresponds to $\left[\mathbf{A}_{qp}\right]$, and $q = \left(\alpha-1\right)*dm + k, p = \left(\alpha_s^i-1\right)*dm + j$; ${\rm B}_{k}^{isj} \left(\mathbf{x}_\alpha\right)$ corresponds to $\left[\mathbf{B}_{qp}\right]$, and $q = \left(\alpha-1\right)*dm + k, p = \left(\alpha_s^i-1\right)*dm + j$; $\mathbf{I}$ is identity matrix. The location of the entry corresponding to (\ref{fa31}) in the matrix corresponding to (\ref{fa32}) is shown in Figure \ref{sc9}.
Merging two items on the left hand side of (\ref{fa32}), we can get the following matrix form:
\begin{sequation}\label{faa32}
\left[\mathbf{C}\right] \left[\mathbf{U}\right] = \left[\mathbf{B}\right] \left[\mathbf{T}\right] + \left[\mathbf{F}\right]
\end{sequation}
where $\left[\mathbf{C}\right] = \left[\mathbf{I}\right] + \left[\mathbf{A}\right]$.

\subsubsection{The calculation of the discretized boundary element formulation}\label{atl33}

Now, we tend to the calculations of the three integrals in (\ref{fa30}). Substituting (\ref{fa27}), (\ref{fa13ab}) and the definition of the deformation gradient (\ref{fabb1}) into (\ref{fa30}), we can get
\begin{sequation}\label{faaa22}
\left\{
\begin{aligned}
& {\rm A}_{k}^{isj} \left(\mathbf{x}_\alpha\right) = \sum_{l=1}^{dm} \sum_{p=1}^{dm} \sum_{q=1}^{dm} \sum_{w=1}^{dm} \bbint_{\partial \Omega_i} \left[{\rm P}_{kpw}\left(\mathbf{x}' - \mathbf{x}_\alpha\right) {\rm K}^{-1}_{wq} + {\rm P}_{kqw}\left(\mathbf{x}' - \mathbf{x}_\alpha\right) {\rm K}^{-1}_{wp}\right] {\rm n}_{l} {\rm C}_{ljpq} \left(\mathbf{x}' - \mathbf{x}_\alpha\right) \overline{N}_s^i \left(\mathbf{x}'\right) {\rm d}V_{\mathbf{x}'} \\
& {\rm B}_{k}^{isj} \left(\mathbf{x}_\alpha\right) = \bbint_{\partial \Omega_i} {\rm u}_A \left(\left| \mathbf{x}' - \mathbf{x}_\alpha \right|\right) \delta_{kj} \overline{N}_s^i \left(\mathbf{x}'\right) {\rm d}V_{\mathbf{x}'} + \bbint_{\partial \Omega_i} \frac{\left({\rm x}' - {\rm x}_\alpha\right)_k \left({\rm x}' - {\rm x}_\alpha\right)_j}{\left| \mathbf{x}' - \mathbf{x}_\alpha \right|^2} {\rm u}_B \left(\left| \mathbf{x}' - \mathbf{x}_\alpha \right|\right) \overline{N}_s^i \left(\mathbf{x}'\right) {\rm d}V_{\mathbf{x}'}\\
& {\rm F}_k \left(\mathbf{x}_\alpha\right) = \int_{\Omega} {\rm u}_A \left(\left| \mathbf{x}' - \mathbf{x}_\alpha \right|\right) {\rm f}_k \left(\mathbf{x}'\right) {\rm d}V_{\mathbf{x}'} + \sum_{j=1}^{dm} \int_{\Omega} \frac{\left(x' - x_\alpha\right)_k \left({\rm x}' - {\rm x}_\alpha\right)_j}{\left| \mathbf{x}' - \mathbf{x}_\alpha \right|^2} {\rm u}_B \left(\left| \mathbf{x}' - \mathbf{x}_\alpha \right|\right) {\rm f}_j \left(\mathbf{x}'\right) {\rm d}V_{\mathbf{x}'}
\end{aligned}
\right.
\end{sequation}
where $\left| \mathbf{x}' - \mathbf{x}_\alpha \right|$ is the length of the vector $\mathbf{x}' - \mathbf{x}_\alpha$. $\left({\rm x}' - {\rm x}_\alpha\right)_i$ is the $i$ directional component of the vector $\mathbf{x}' - \mathbf{x}_\alpha$. The definition of ${\rm P}_{kpw}\left(\mathbf{x}' - \mathbf{x}_\alpha\right)$ and ${\rm K}_{wq}$ for the first formula of (\ref{faaa22}) are
\begin{sequation}\label{faaa23}
\begin{aligned}
& {\rm P}_{kpw}\left(\mathbf{x}' - \mathbf{x}_\alpha\right) \\
=& \int_{\mathcal{H}_{\mathbf{x}}} \underline{\omega} \left \langle \boldsymbol{\xi} \right \rangle \left({\rm u}_A \left(\left| \mathbf{x}' - \mathbf{x}_\alpha + \boldsymbol{\xi} \right|\right) \delta_{kp} + \frac{\left({\rm x}' - {\rm x}_\alpha + {\rm \xi}\right)_k \left({\rm x}' - {\rm x}_\alpha + {\rm \xi}\right)_p}{\left| \mathbf{x}' - \mathbf{x}_\alpha + \boldsymbol{\xi} \right|^2} {\rm u}_B \left(\left| \mathbf{x}' - \mathbf{x}_\alpha + \boldsymbol{\xi} \right|\right)\right) {\rm \xi}_w {\rm d}V_{\boldsymbol{\xi}} \\
& - \int_{\mathcal{H}_{\mathbf{x}}} \underline{\omega} \left \langle \boldsymbol{\xi} \right \rangle \left({\rm {u}}_A \left(\left| \mathbf{x}' - \mathbf{x}_\alpha \right|\right) \delta_{kp} + \frac{\left({\rm x}' - {\rm x}_\alpha\right)_k \left({\rm x}' - {\rm x}_\alpha\right)_p}{\left| \mathbf{x}' - \mathbf{x}_\alpha \right|^2} {\rm {u}}_B \left(\left| \mathbf{x}' - \mathbf{x}_\alpha \right|\right)\right) {\rm \xi}_w {\rm d}V_{\boldsymbol{\xi}}
\end{aligned}
\end{sequation}
\begin{sequation}\label{faaa24}
{\rm K}_{wq} = \int_{\mathcal{H}_{\mathbf{x}}} \underline{\omega} \left \langle \boldsymbol{\xi} \right \rangle {\rm \xi}_w {\rm \xi}_q {\rm d}V_{\boldsymbol{\xi}}
\end{sequation}
where $\underline{\omega} \left \langle \boldsymbol{\xi} \right \rangle$ is the nonlocal weight function. Two types of nonlocal weight functions $\underline{\omega} \left \langle \boldsymbol{\xi} \right \rangle$, namely, the constant distribution and the Gauss distribution, are widely used in the peridynamic theory
\begin{sequation}\label{faaa24a}
\underline{\omega} \left \langle \boldsymbol{\xi} \right \rangle =
\begin{cases}
1  & \left|\boldsymbol{\xi}\right| \leq {\rm h_r} \\
0  & \left|\boldsymbol{\xi}\right| > {\rm h_r}
\end{cases}
\end{sequation}
\begin{sequation}\label{faaa24b}
\underline{\omega} \left \langle \boldsymbol{\xi} \right \rangle = \exp\left(-\dfrac{\left|\boldsymbol{\xi}\right|^2}{\left({\rm h_r}\right)^2}\right)
\end{sequation}
where ${\rm h_r}$ is the characteristic length. We see that only the integral involving  ${\rm u}_A$ part of the infinite Green function needs spacial treatment because of its singularity, whereas the other integrals without singularities are easy to handle with the ordinary Gaussian integration. The integrals involving ${\rm u}_A$ are divided into the following two boundary and volume categories:
\begin{sequation}\label{faaa25}
\bbint_{\partial \Omega} {\rm u}_A \left(\left| \mathbf{x}' - \mathbf{x}_0 \right|\right) \mathbf{f} \left(\mathbf{x}\right) {\rm d}V_{\mathbf{x}}
\end{sequation}
\begin{sequation}\label{faaa26}
\int_{\Omega} {\rm u}_A \left(\left| \mathbf{x}' - \mathbf{x}_0 \right|\right) \mathbf{f} \left(\mathbf{x}\right) {\rm d}V_{\mathbf{x}}
\end{sequation}
where $\mathbf{f} \left(\mathbf{x}\right) \in C^{\infty}$ and $\mathbf{x}_0 \in \partial \Omega$. Substituting (\ref{faaa23}) into (\ref{faaa22})$_{1}$ results in the following integrals:
\begin{sequation}\label{faaa27}
\bbint_{\partial \Omega_i} \int_{\mathcal{H}_{\mathbf{x}}} \underline{\omega} \left \langle \boldsymbol{\xi} \right \rangle {\rm u}_A \left(\left| \mathbf{x}' - \mathbf{x}_\alpha + \boldsymbol{\xi} \right|\right) \delta_{kp} {\rm \xi}_w {\rm K}^{-1}_{wq} {\rm n}_{l} {\rm C}_{ljpq} \left(\mathbf{x}' - \mathbf{x}_\alpha\right) \overline{N}_s^i \left(\mathbf{x}'\right) {\rm d}V_{\boldsymbol{\xi}} {\rm d}V_{\mathbf{x}'}
\end{sequation}
\begin{sequation}\label{faaa28}
\bbint_{\partial \Omega_i} \int_{\mathcal{H}_{\mathbf{x}}} \underline{\omega} \left \langle \boldsymbol{\xi} \right \rangle {\rm u}_A \left(\left| \mathbf{x}' - \mathbf{x}_\alpha \right|\right) \delta_{kp} {\rm \xi}_w {\rm K}^{-1}_{wq} {\rm n}_{l} {\rm C}_{ljpq} \left(\mathbf{x}' - \mathbf{x}_\alpha\right) \overline{N}_s^i \left(\mathbf{x}'\right) {\rm d}V_{\boldsymbol{\xi}} {\rm d}V_{\mathbf{x}'}
\end{sequation}
They have the same forms as (\ref{faaa26}) and (\ref{faaa25}), respectively. The first integral on the right hand side of (\ref{faaa22})$_{2}$ is
\begin{sequation}\label{faaa29}
\bbint_{\partial \Omega_i} {\rm u}_A \left(\left| \mathbf{x}' - \mathbf{x}_\alpha \right|\right) \delta_{kj} \overline{N}_s^i \left(\mathbf{x}'\right) {\rm d}V_{\mathbf{x}'},
\end{sequation}
which has the same form as (\ref{faaa25}). The first integral on the right hand side of (\ref{faaa22})$_{3}$ is
\begin{sequation}\label{faaa30}
\int_{\Omega} {\rm u}_A \left(\left| \mathbf{x}' - \mathbf{x}_\alpha \right|\right) {\rm f}_k \left(\mathbf{x}'\right) {\rm d}V_{\mathbf{x}'},
\end{sequation}
which has the same form as (\ref{faaa26}). According to Wang et al.~\cite{tb58}, $u_A$ can be decomposed into three parts in the two-dimension problem: the classical part, the nonlocal and convergent part and the divergent part. They are given as
\begin{sequation}\label{faaa31}
\left({\rm u}_A\right)^{classical} \left(x\right) = \dfrac{-8\ln x}{9\pi E}
\end{sequation}
\begin{sequation}\label{faaa32}
\begin{aligned}
\left({\rm u}_A\right)^{nonlocal} \left(x\right) =& \dfrac{1}{2 \pi} \int_{0}^{+\infty} J_0 \left(kx\right) k \left(\dfrac{1}{M_\perp^{\uppercase\expandafter{\romannumeral2}} \left(k\right)} - \dfrac{1}{M^{\uppercase\expandafter{\romannumeral2}} \left(\infty\right)} - \dfrac{8}{3Ek^2}\right) {\rm d}k \\
& + \dfrac{1}{2 \pi} \int_{0}^{+\infty} \dfrac{J_1 \left(kx\right)}{x} \left(\dfrac{1}{M_\Vert^{\uppercase\expandafter{\romannumeral2}} \left(k\right)} - \dfrac{1}{M_\perp^{\uppercase\expandafter{\romannumeral2}} \left(k\right)} + \dfrac{16}{9Ek^2} \right) {\rm d}k
\end{aligned}
\end{sequation}
\begin{sequation}\label{faaa33}
\left({\rm u}_A\right)^{divergent} \left(x\right) = \frac{\delta^{\uppercase\expandafter{\romannumeral2}} \left(\mathbf{x}\right)}{M^{\uppercase\expandafter{\romannumeral2}} \left(\infty\right)}
\end{sequation}
Therefore, the integrals in (\ref{faaa25}) and (\ref{faaa26}) are divided into three parts, based on (\ref{faaa31}), (\ref{faaa32}) and (\ref{faaa33}). The part involving (\ref{faaa32}) is nonsingular, so the ordinary Gaussian integration is sufficient to deal with it. The part involving (\ref{faaa31}) corresponds to the local boundary element,  so it can be calculated by following the classical boundary element method~\cite{tb60}. Next, we will focus on the integrals that involve (\ref{faaa33}). They are actually the convolutions of the Dirac function  in $\partial \Omega$ and $\Omega$
\begin{sequation}\label{faaa34}
\bbint_{\partial \Omega} \delta^{\uppercase\expandafter{\romannumeral2}} \left(\mathbf{x} - \mathbf{x}_0\right) \mathbf{f} \left(\mathbf{x}\right) {\rm d}V_{\mathbf{x}}
\end{sequation}
\begin{sequation}\label{faaa35}
\int_{\Omega} \delta^{\uppercase\expandafter{\romannumeral2}} \left(\mathbf{x} - \mathbf{x}_0\right) \mathbf{f} \left(\mathbf{x}\right) {\rm d}V_{\mathbf{x}}
\end{sequation}
According to the definition of the Cauchy principal value, (\ref{faaa34}) is expressed as
\begin{sequation}\label{fa331}
\bbint_{\partial \Omega} \delta \left(\mathbf{x} - \mathbf{x}_0\right) \mathbf{f} \left(\mathbf{x}\right) {\rm d}V_{\mathbf{x}} = \lim_{\varepsilon \to 0} \int_{\partial \Omega \setminus \left\{U_0 \left(\mathbf{x}_0,\varepsilon\right) \cap \partial \Omega \right\}} \delta \left(\mathbf{x} - \mathbf{x}_0\right) \mathbf{f} \left(\mathbf{x}\right) {\rm d}V_{\mathbf{x}} = 0
\end{sequation}
Using convolution invariance of the Dirac function, the result for (\ref{faaa35}) is $\mathbf{f} \left(\mathbf{x}_0\right)$.
In addition, for matrix $\mathbf{A}$, we need not compute its singular integrals (\ref{faaa27}) and (\ref{faaa28}) directly, because we can remove its singularity by applying a rigid body displacement~\cite{tb60}. For the three dimensional problem, we can merge the classical part and the nonlocal convergent part of  ${\rm u}_A$~\cite{tb58} as follows:
\begin{sequation}\label{fa335}
\begin{aligned}
\left({\rm u}_A\right)^{nonsingular} \left(x\right) =& \dfrac{1}{2 \pi^2} \int_{0}^{+\infty} \left(\dfrac{{\rm sin}\left(kx\right)}{\left(kx\right)^3} - \dfrac{{\rm cos}\left(kx\right)}{\left(kx\right)^2} \right) \left(\dfrac{1}{ M_\Vert^{\uppercase\expandafter{\romannumeral3}} \left(k\right)} - \dfrac{1}{M_\perp^{\uppercase\expandafter{\romannumeral3}} \left(k\right)}\right) k^2 {\rm d}k \\
& + \dfrac{1}{2 \pi^2} \int_{0}^{+\infty} \dfrac{{\rm sin}\left(kx\right)}{x} \left( \dfrac{k}{M_\perp^{\uppercase\expandafter{\romannumeral3}} \left(k\right)} - \dfrac{k}{M^{\uppercase\expandafter{\romannumeral3}} \left(\infty\right)} \right) {\rm d}k
\end{aligned}
\end{sequation}
The divergent part of ${\rm u}_A$ is
\begin{sequation}\label{fa336}
\left({\rm u}_A\right)^{divergent} \left(x\right) = \frac{\delta^{\uppercase\expandafter{\romannumeral3}} \left(\mathbf{x}\right)}{M^{\uppercase\expandafter{\romannumeral3}} \left(\infty\right)}
\end{sequation}
Thus, the integrals for (\ref{faaa25}) and (\ref{faaa26}) are divided into two parts. According to Wang et al.~\cite{tb58}, (\ref{fa335}) is smooth, bounded, and nonsingular, and then the ordinary Gaussian integration is sufficient to deal with it. The integral that is related to the divergent part (\ref{fa336}) is the same as the two-dimension problem, and it is zero.

\subsection{Numerical treatment for dynamic problems}\label{atl34}

Since we have introduced the Laplace transformation to deal with the dynamic problem, solving the dynamic problem will be divided into two steps. The first step is to calculate the deformation in the Laplace domain, and the second step is to carry out the transformation inversion to get the dynamic response in the time domain. For the first step, it is similar to the preceding method of the static problem in Section \ref{atl35}. For the second step, the specific Laplace transformation inversion methods have been introduced in detail in the literature~\cite{tb67}. We choose the category that is usually used in the boundary element method. The details of the method can be found in the literature~\cite{tb69}, with an overview given in Appendix \ref{ap9}. The Laplace transformation inversion~\cite{tb69} is expressed as
\begin{sequation}\label{faaa36}
f \left(t\right) \approx \dfrac{2{\rm e}^{at}}{T} \left[\dfrac{Re\left\{F\left(a\right)\right\}}{2} + \sum_{k = 1}^{N} Re\left\{F\left(a + \dfrac{k\pi {\rm i}}{T}\right)\right\}{\rm cos} \dfrac{k\pi t}{T} \right]
\end{sequation}
where $f \left(t\right)$ and $F \left(s\right)$ are the original function in the time domain and the image function in the Laplace domain, respectively. $T$ is the half of the total time that needs to be calculated. $N$ is the number of terms for the numerical truncation. $a$ is the parameter of convergence control, usually $aT=10$. ${\rm i}=\sqrt{-1}$. $Re$ represents the real part. It should be noted that theoretically, the infinite summation is needed for calculating the inversion. In fact, we only need to control the parameters $a$ and $N$ to make the summation convergent with finite terms.

According to the numerical approximation formula of the Laplace transformation inversion, we only need to compute the value of $F \left(s\right)$ in the point $s_k$, and the calculation for different $F \left(s_k\right)$ is independent, which is the basis of parallel computation. According to (\ref{faaa36}), $s_k = a + \dfrac{k\pi {\rm i}}{T}$. Comparison of (\ref{fa25}) and (\ref{fa117}) shows that the static boundary integral equation and the dynamic boundary integral equation have the same form, if we regard the variable $s$ in the Laplace domain as a parameter. Therefore, we only need to take $s$ as $s_k$ in (\ref{fa117}) and solve $\tilde{\mathbf{u}} \left(\mathbf{x},s\right)$ and $\tilde{\mathbf{T}} \left(\mathbf{x},s\right)$ in (\ref{fa117}) with the static PD-BEM. This method has the advantages of fast calculation and simple realization.

For the dynamic problem, we also need to deal with the integrals for  $\left(\tilde{{\rm u}}_g\right)_{A}$  of the infinite Green function in (\ref{fa39}), as in (\ref{faaa25}) and (\ref{faaa26}). In the same way as that in Section \ref{atl33}, the expression of $\left(\tilde{{\rm u}}_g\right)_{A}$ in (\ref{fa45}) can be decomposed into two parts: the classical part, and the nonlocal convergent part
\begin{sequation}\label{faaa37}
\left(\tilde{{\rm u}}_g\right)_{A}^{divergent} \left(x,s\right) = \frac{\delta^{\uppercase\expandafter{\romannumeral2}} \left(\mathbf{x}\right)}{\rho s^2 + M^{\uppercase\expandafter{\romannumeral2}} \left(\infty\right)}
\end{sequation}
\begin{sequation}\label{faaa38}
\begin{aligned}
\left(\tilde{{\rm u}}_g\right)_{A}^{nonlocal} \left(x,s\right) =& \frac{1}{2 \pi} \int_{0}^{+\infty} J_0 \left(kx\right)\left(\frac{1}{\rho s^2 + M_\perp^{\uppercase\expandafter{\romannumeral2}} \left(k\right)} - \frac{1}{\rho s^2 + M^{\uppercase\expandafter{\romannumeral2}} \left(\infty\right)}\right) k {\rm d}k \\
& + \frac{1}{2 \pi} \int_{0}^{+\infty} \frac{J_1 \left(kx\right)}{x} \left(\frac{1}{\rho s^2 + M_\Vert^{\uppercase\expandafter{\romannumeral2}} \left(k\right)} - \frac{1}{\rho s^2 + M_\perp^{\uppercase\expandafter{\romannumeral2}} \left(k\right)}\right) {\rm d}k
\end{aligned}
\end{sequation}
The treatment for the divergent part is the same as that in Section \ref{atl33}. The result will also be zero. We only use the points which are not on the negative real axis for the complex variable $s$, that is to say, $s^2$ is not a negative real number; then the nonlocal  part in (\ref{faaa38}) is nonsingular.

\section{Numerical examples}\label{atl4}

In this section, we first give three numerical examples to verify our numerical method from the aspects of dynamics and statics. Secondly, we give an additional example to show the advantages of our method, in comparison with the PD-MPM in this case. Thirdly, we compute the fracture process of the Brazilian disk by coupling the PD-BEM with the PD-MPM to improve the efficiency. The computations are performed on a PC with a CPU of Intel Core i7-6700 3.40GHZ, RAM of 16.0GB, and Mononuclear, except the last example which is computed by the SMP (Symmetric Multiprocessor) parallel computing.  The model of the computing node is NeXtScale nx360 M5, and the parameter of the computing node is 2*Intel Xeon E5-2697A V4, 128G.

\subsection{A two dimensional problem }\label{atl41}

\begin{figure}[!htb]
	\centerline{\includegraphics[scale=0.3]{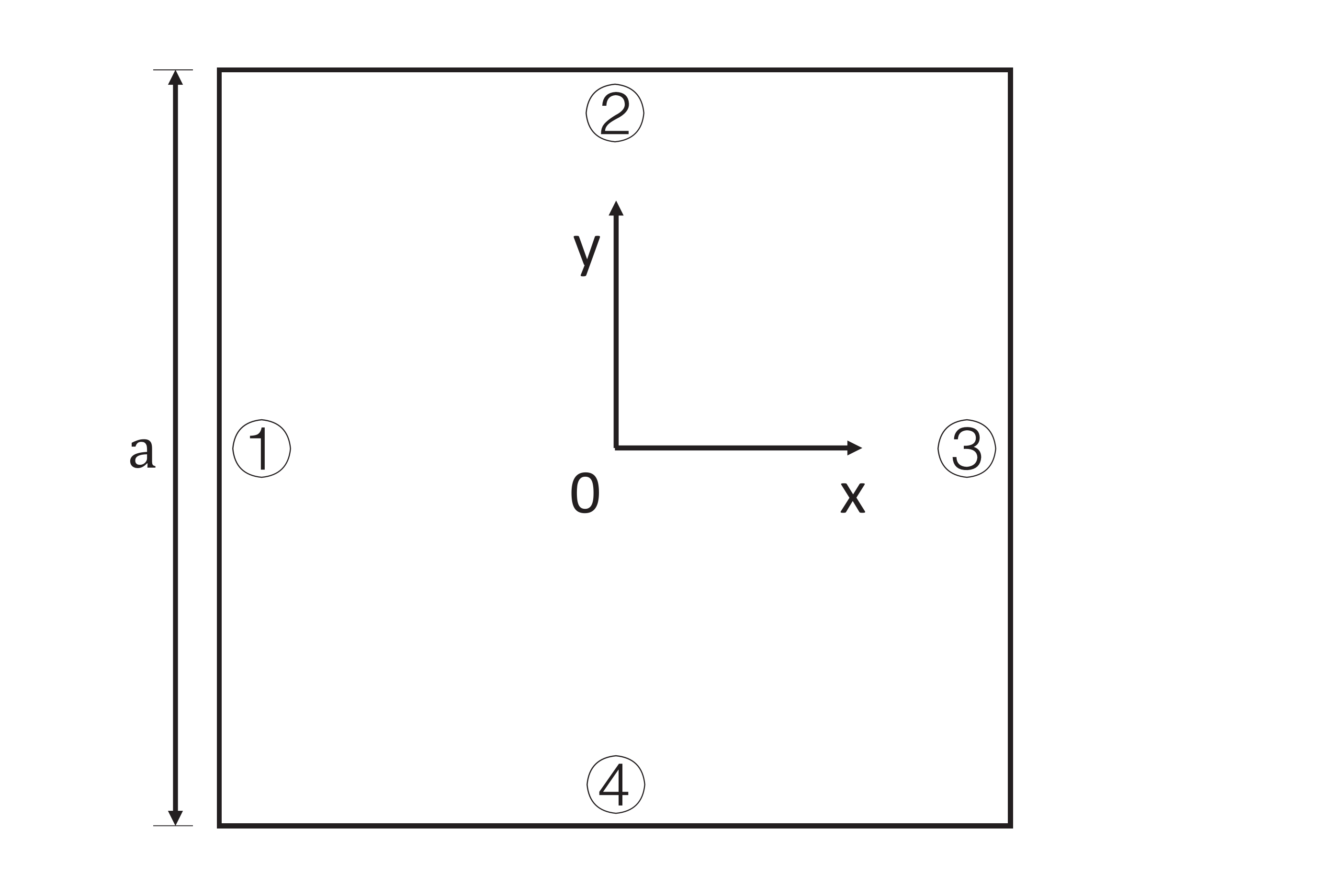}}	
	\caption{A square domain where the numbers in the circle denote the sides.\label{sc1}}
\end{figure}

We apply the PD-BEM  to a square domain, shown in Figure \ref{sc1}, which is subjected to three loading conditions. The first two are biaxial tension and self weight loading. The major purpose of these two cases are to verify the method and demonstrate its convergence to the classical theory when the size of the horizon tends to zero. The third loading case is uniaxial tension, the results of which demonstrate the efficiency and accuracy of the PD-BEM.  The used parameters are listed in Table \ref{tab2}. The boundary and loading conditions are listed in Table \ref{tab3}.
\begin{center}
	\begin{table}[!htb]
		\caption{The geometric and  material parameters for the two-dimensional problem of the square domain.\label{tab2}}
		\centering
		\begin{tabular}{ccccc}
			\hline
			Loading conditions&${\rm a}\left({\rm m}\right)$&$E \left({\rm Pa}\right)$&$\rho \left({\rm kg/m^2}\right)$&$\nu$\\
			\hline
			1,2&$3$&$1$&$1$&$1/3$\\
			3&$1.0$&$1$&$1$&$1/3$\\
			\hline
		\end{tabular}
	\end{table}
\end{center}
\begin{center}
	\begin{table}[!htb]
		\centering
		\caption{The boundary and loading conditions for the two-dimensional problem of the square domain.\label{tab3}}
		\begin{tabular}{ccc}
			\hline
			Loading&\multirow{2}*{}&\multirow{2}*{Boundary conditions}\\
			conditions&&\\
			\hline
			\multirow{2}*{1}&\multirow{2}*{Biaxial tension}&\ding{172},\ding{174} --- ${\rm t}_{n}=0.03 {\rm N/m}$, ${\rm t}_{\tau}=0.0 {\rm N/m}$, ${\rm b}_{x}=0.0 N/m^2$\\
			&&\ding{173},\ding{175} --- ${\rm t}_{n}=0.03 {\rm N/m}$, ${\rm t}_{\tau}=0.0 {\rm N/m}$, ${\rm b}_{y}=0.0 N/m^2$\\
			\multirow{2}*{2}&\multirow{2}*{Self weight}&\ding{172},\ding{174} --- ${\rm u}_{n}=0.0 {\rm m}$, ${\rm u}_{\tau}=0.0 {\rm m}$, ${\rm b}_{x}=0.0 {\rm N/m^2}$\\
			&&\ding{173},\ding{175} --- ${\rm t}_{n}=0.0 {\rm N/m}$, ${\rm t}_{\tau}=0.0 {\rm N/m}$, ${\rm b}_{y}=-1.0 {\rm N/m^2}$\\
			\multirow{2}*{3}&\multirow{2}*{Uniaxial tension}&\ding{172},\ding{174} --- ${\rm t}_{n}=0.0 {\rm N/m}$, ${\rm t}_{\tau}=0.0 {\rm N/m}$, ${\rm b}_{x}=0.0 {\rm N/m^2}$\\
			&&\ding{173},\ding{175} --- ${\rm t}_{n}=1.0 {\rm N/m}$, ${\rm t}_{\tau}=0.0 {\rm N/m}$, ${\rm b}_{y}=0.0 {\rm N/m^2}$\\
			\hline
		\end{tabular}
	\end{table}
\end{center}

\begin{figure}[!htb]
	\centering
	\subfigure[]{
		\label{fig13a}
		\includegraphics[scale=0.16]{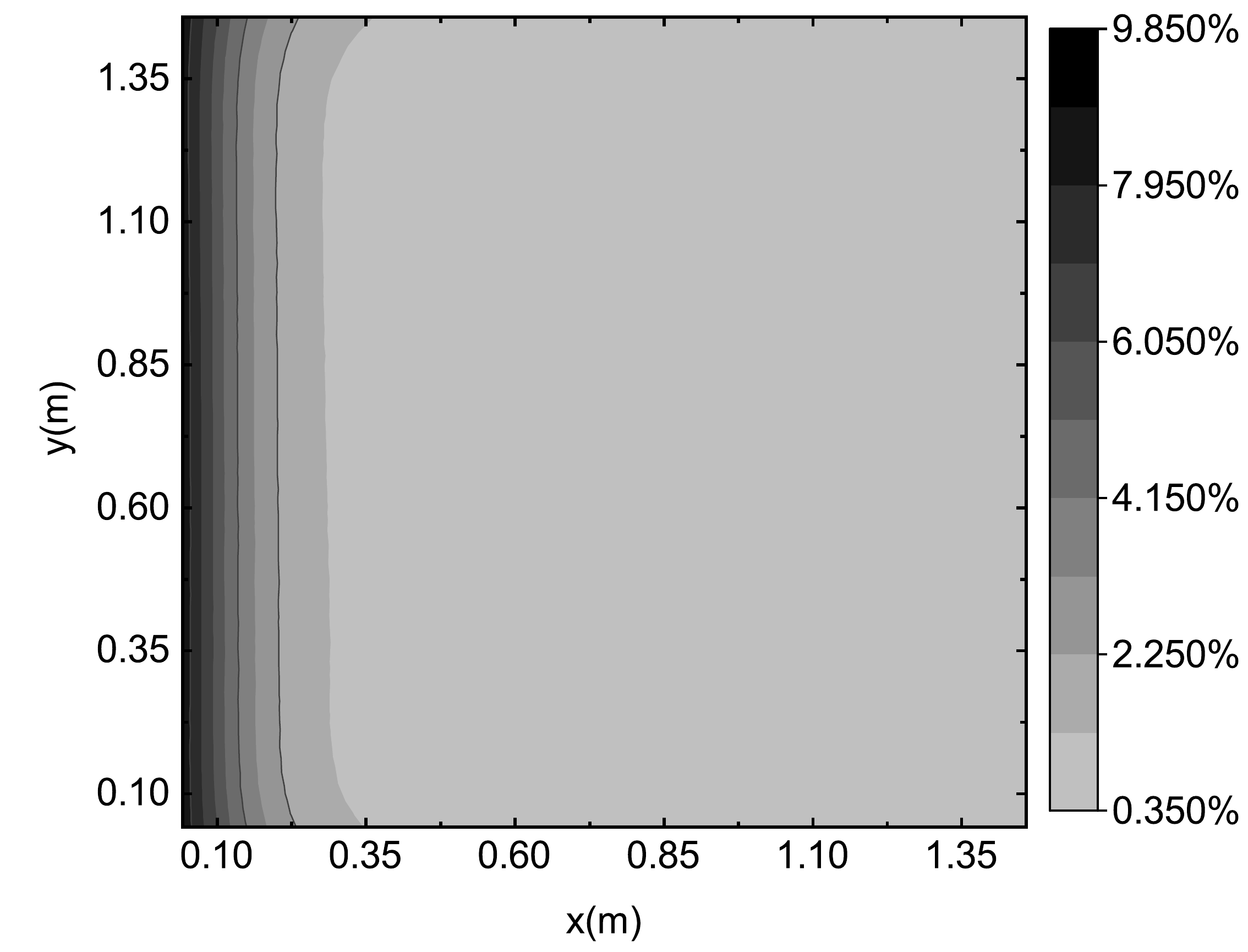}}
	\hspace{0.01in}
	\subfigure[]{
		\label{fig13b}
		\includegraphics[scale=0.16]{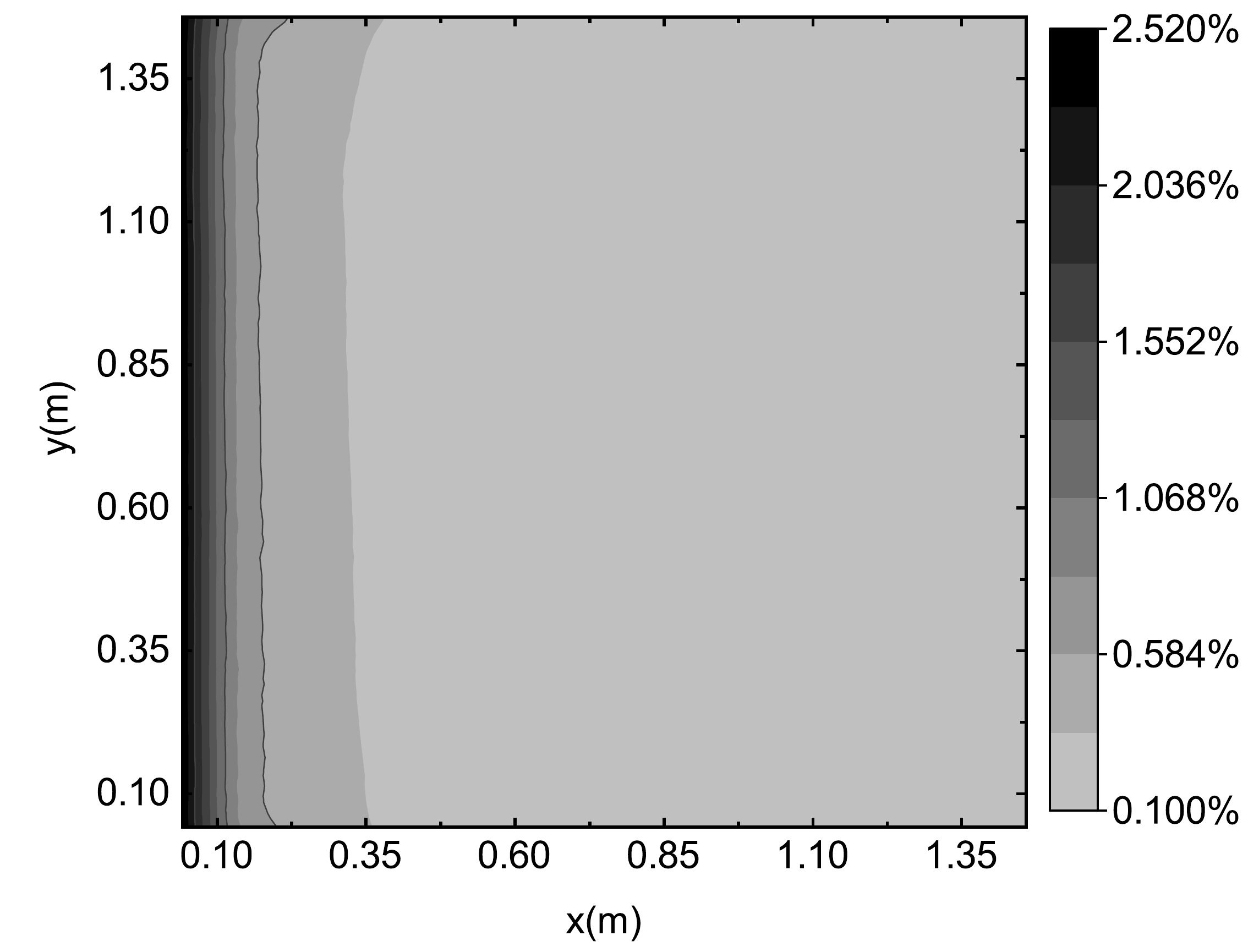}}
	\hspace{0.01in}
	\subfigure[]{
		\label{fig13c}
		\includegraphics[scale=0.16]{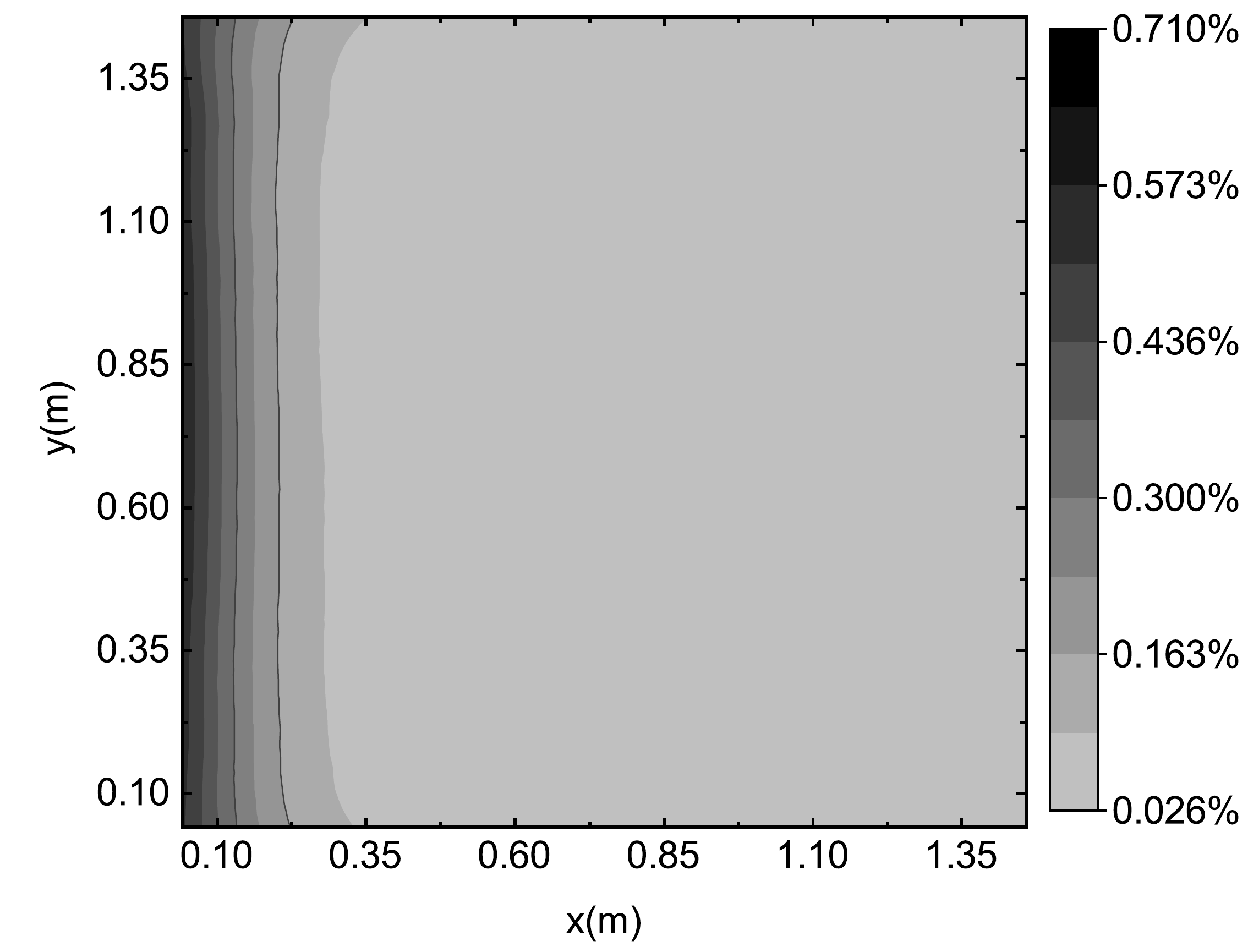}}
	\caption{The relative difference between the PD-BEM results the classical BEM results under the first loading condition with the constant value kernel: ${\rm h_r} = 1/100 {\rm m}$(a); ${\rm h_r} = 1/400 {\rm m}$(b); ${\rm h_r} = 1/1600 {\rm m}$(c).}
	\label{fig13}
\end{figure}

\begin{figure}[!htb]
	\centering
	\subfigure[]{
		\label{fig14a}
		\includegraphics[scale=0.16]{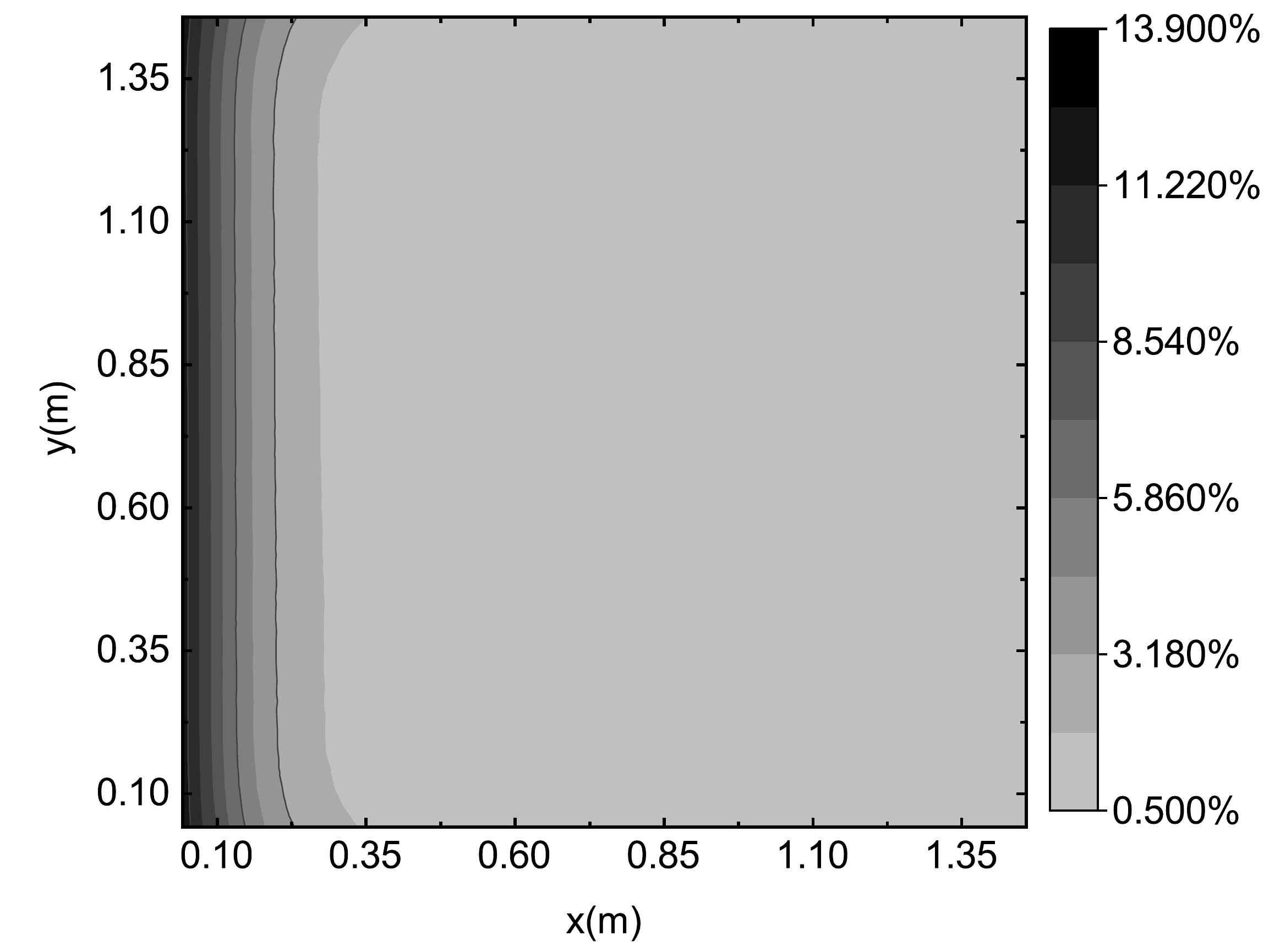}}
	\hspace{0.01in}
	\subfigure[]{
		\label{fig14b}
		\includegraphics[scale=0.16]{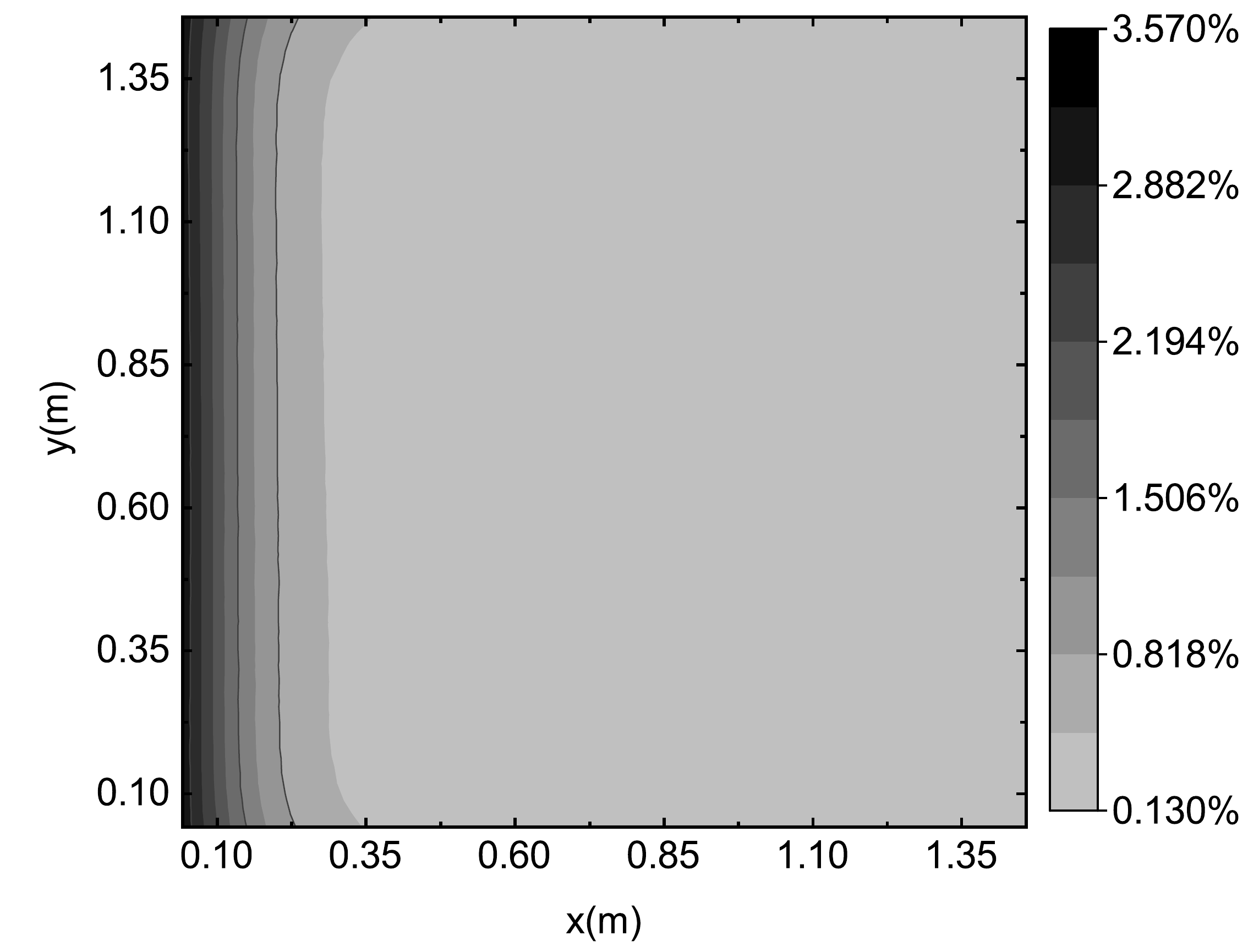}}
	\hspace{0.01in}
	\subfigure[]{
		\label{fig14c}
		\includegraphics[scale=0.16]{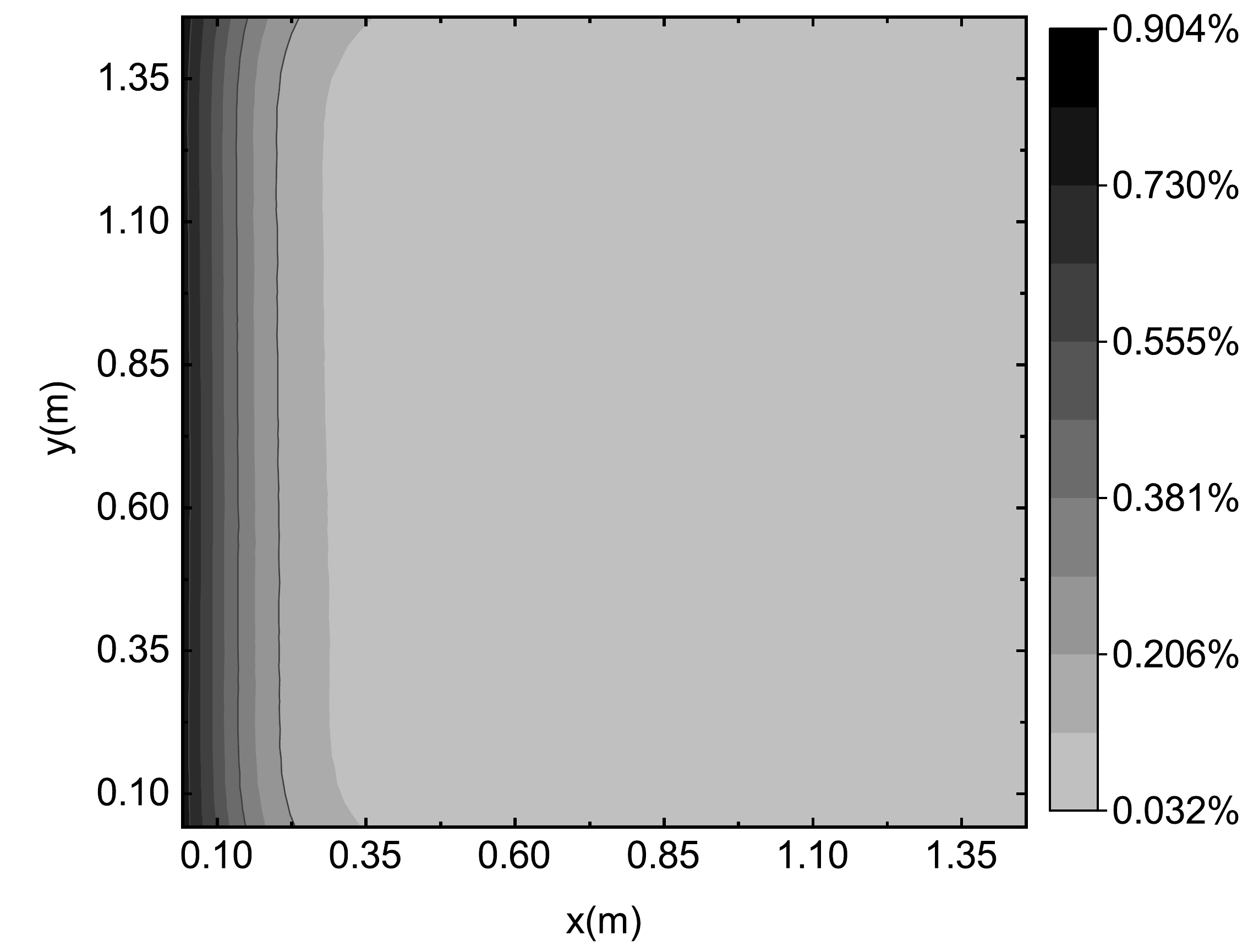}}
	\caption{The relative difference between the PD-BEM results the classical BEM results under the first loading condition with
		the Gauss kernel: ${\rm h_r} = 1/100 {\rm m}$(a); ${\rm h_r} = 1/400 {\rm m}$(b); ${\rm h_r} = 1/1600 {\rm m}$(c).}
	\label{fig14}
\end{figure}

\begin{figure}[!htb]
	\centering
	\subfigure[]{
		\label{fig15d}
		\includegraphics[scale=0.16]{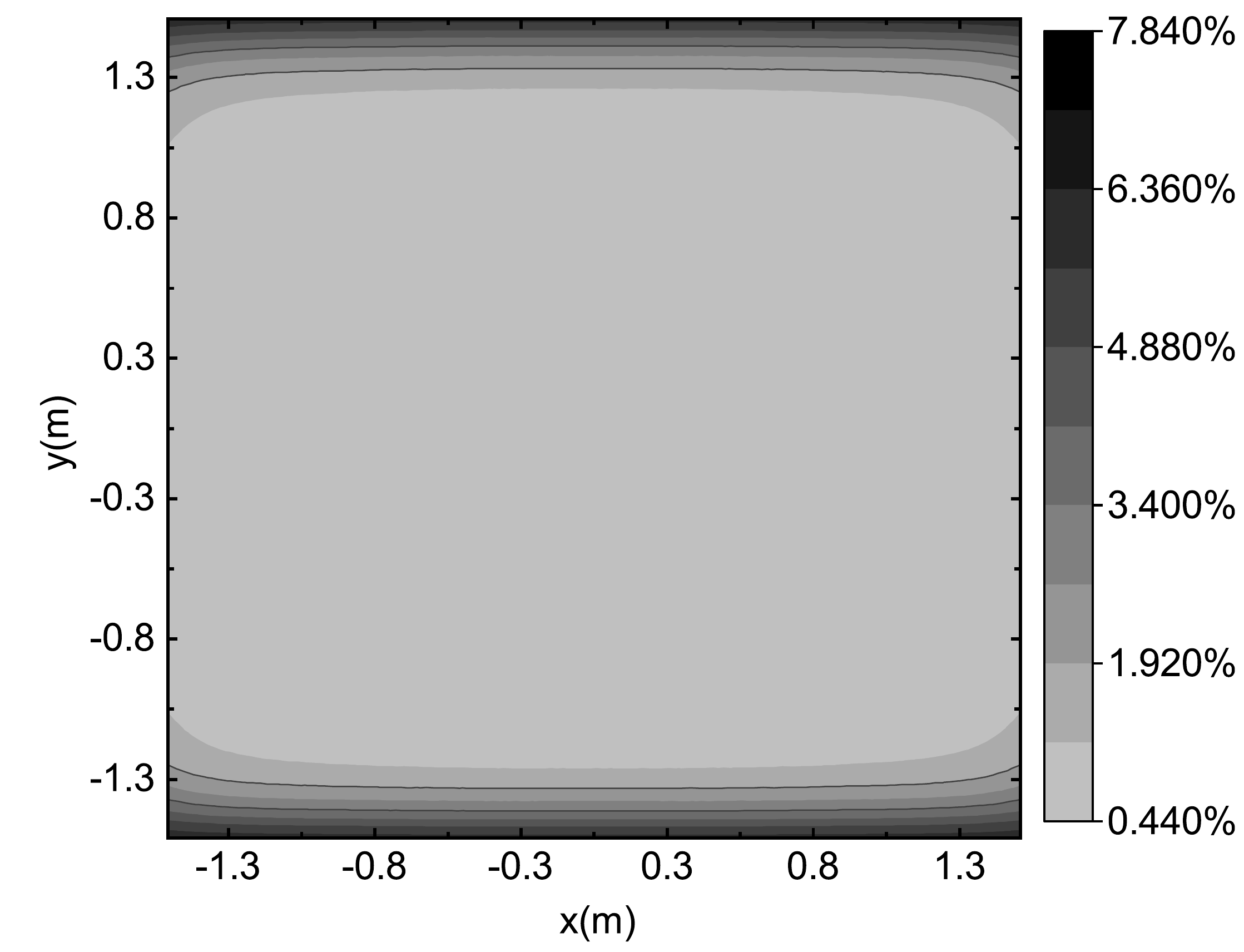}}
	\hspace{0.01in}
	\subfigure[]{
		\label{fig15e}
		\includegraphics[scale=0.16]{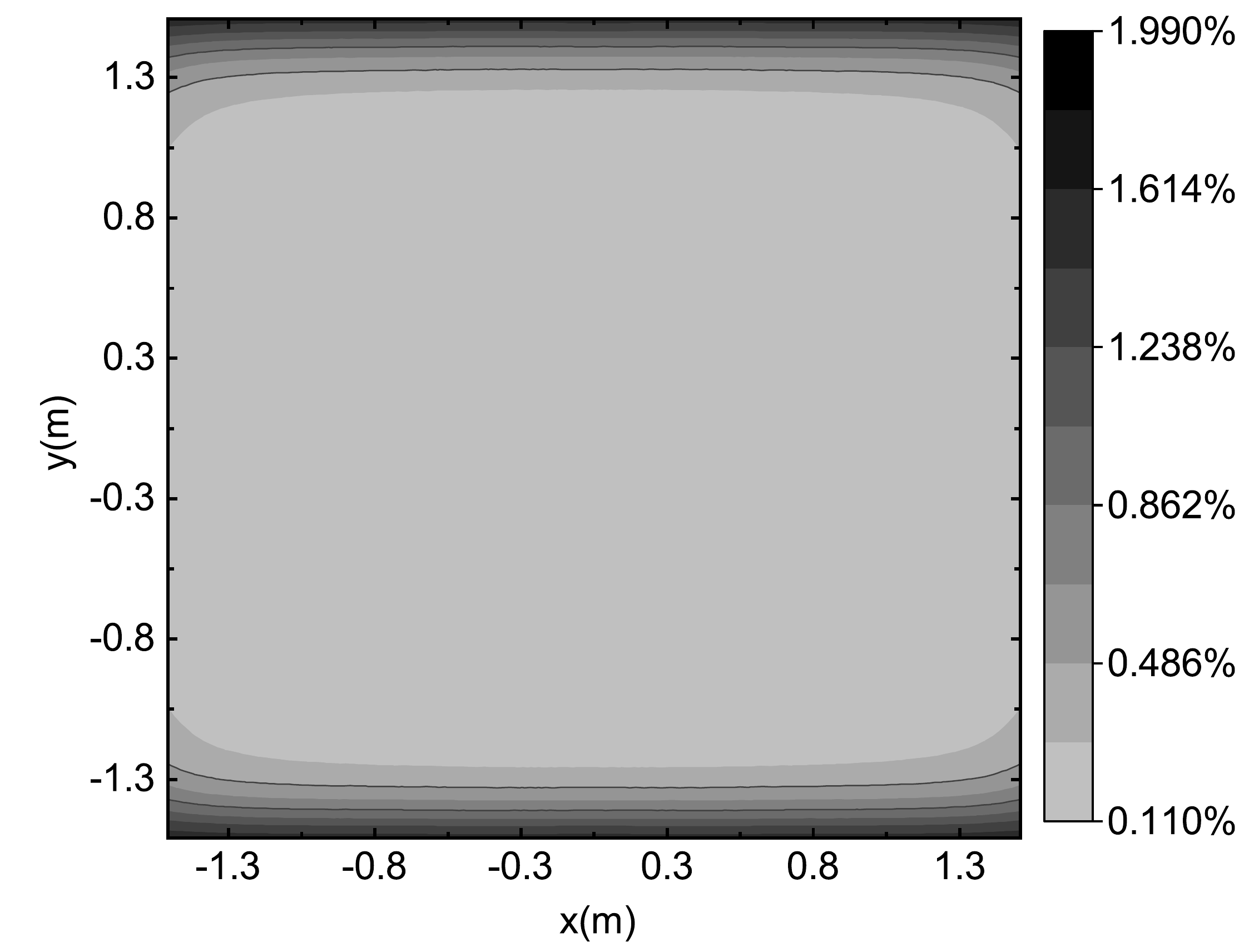}}
	\hspace{0.01in}
	\subfigure[]{
		\label{fig15f}
		\includegraphics[scale=0.16]{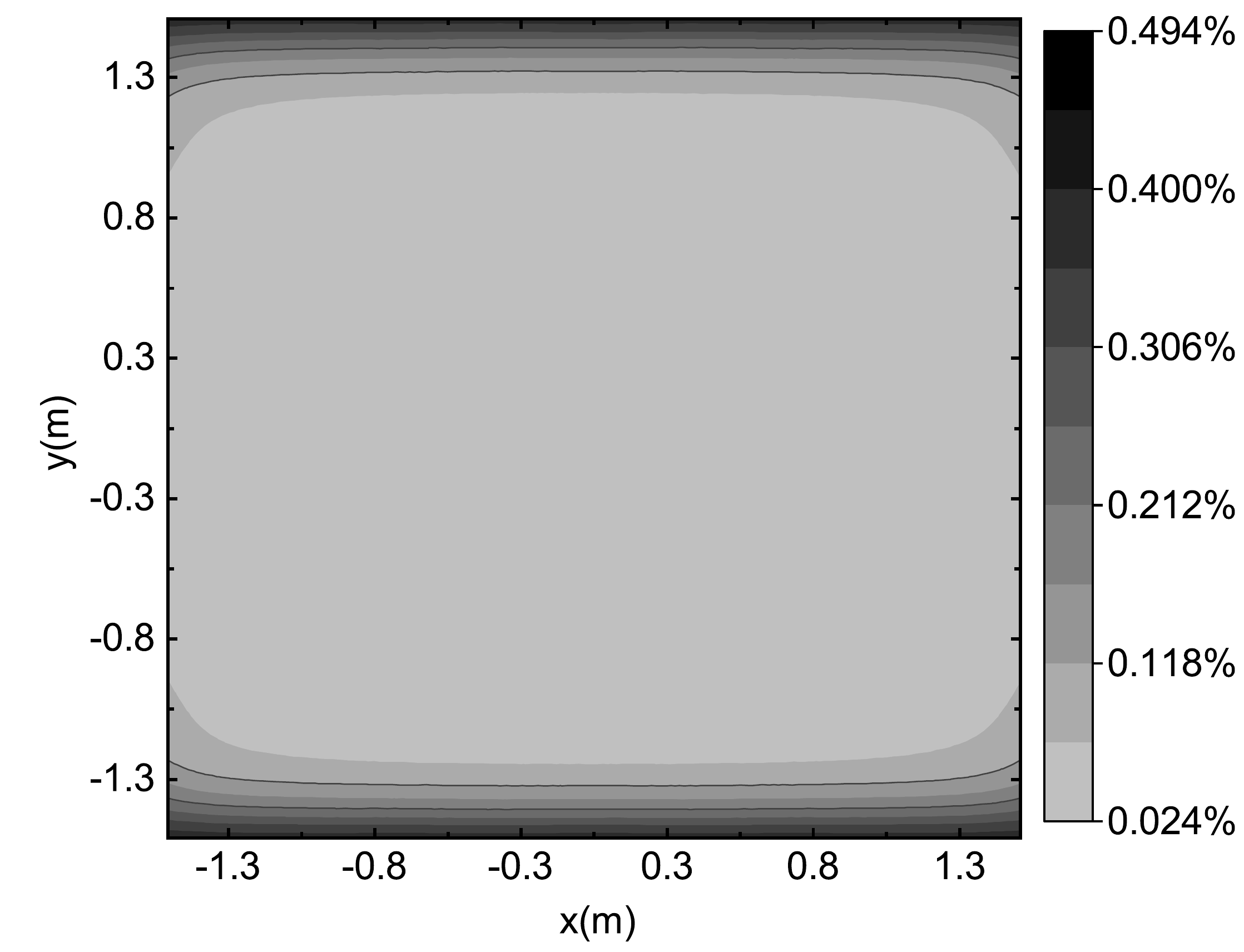}}
	\caption{The relative difference between the PD-BEM results the classical BEM results under the second loading condition with
		the constant value kernel: ${\rm h_r} = 1/100 {\rm m}$(a); ${\rm h_r} = 1/400 {\rm m}$(b); ${\rm h_r} = 1/1600 {\rm m}$(c).}
	\label{fig15}
\end{figure}

\begin{figure}[!htb]
	\centering
	\subfigure[]{
		\label{fig16d}
		\includegraphics[scale=0.16]{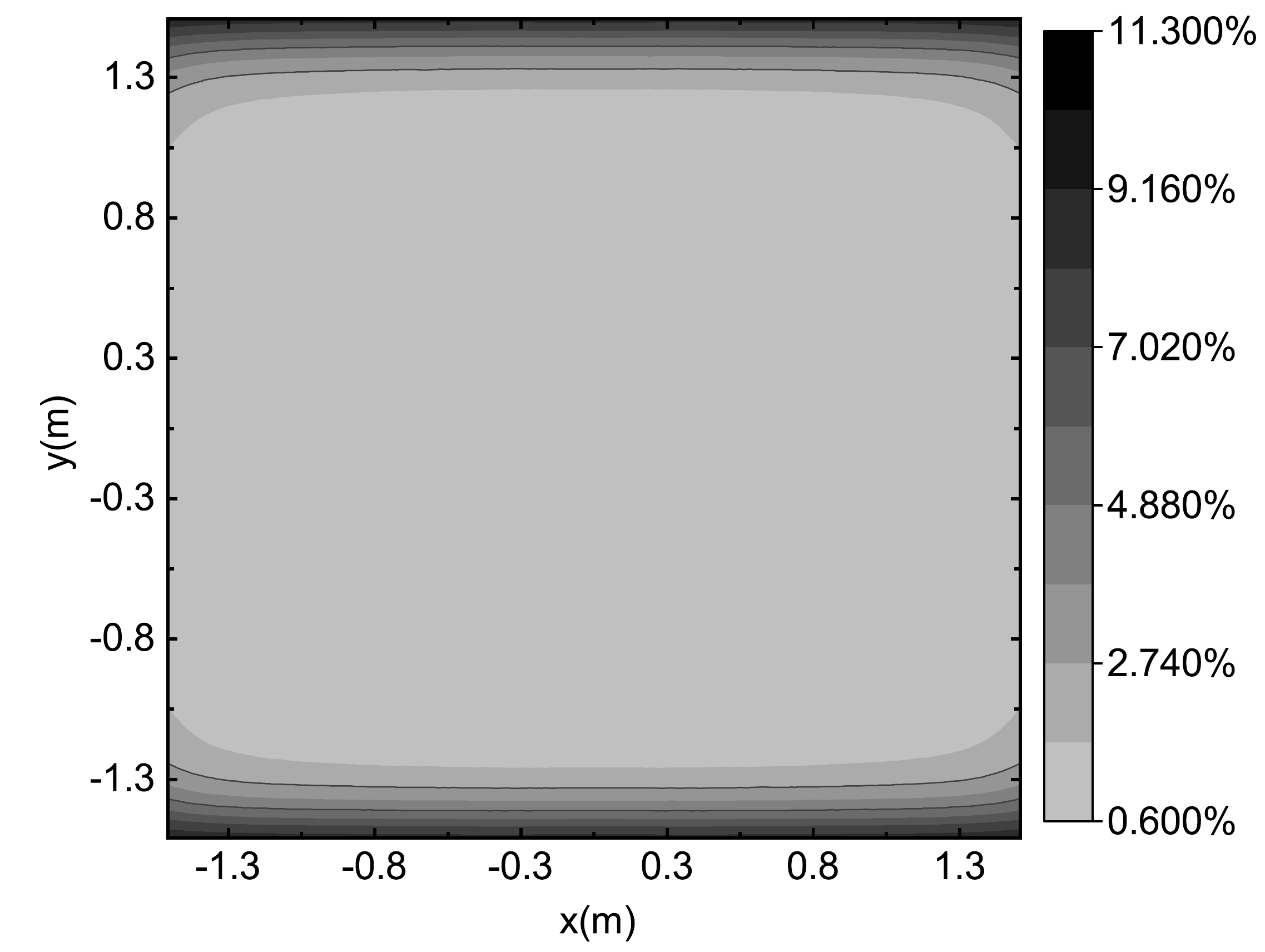}}
	\hspace{0.01in}
	\subfigure[]{
		\label{fig16e}
		\includegraphics[scale=0.16]{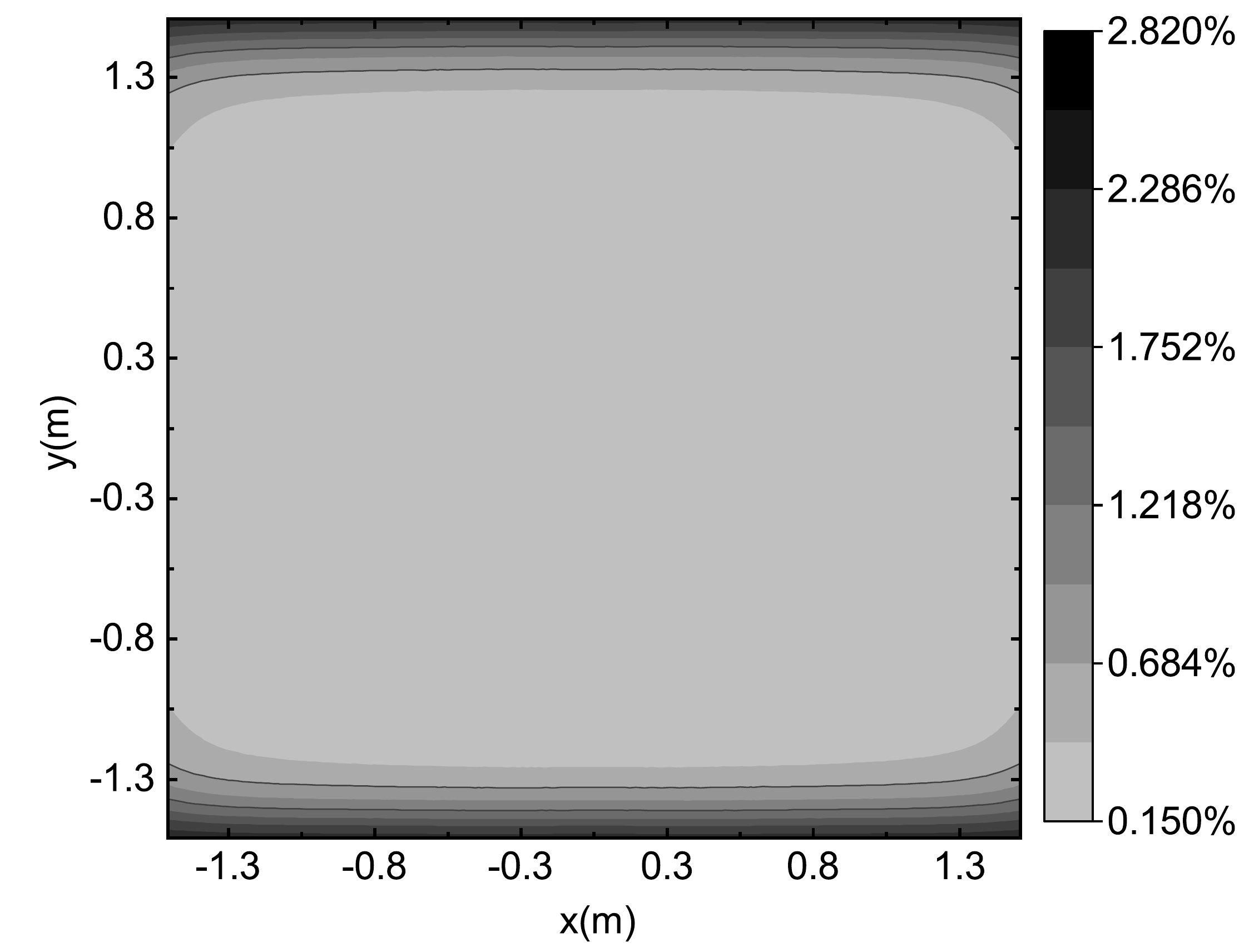}}
	\hspace{0.01in}
	\subfigure[]{
		\label{fig16f}
		\includegraphics[scale=0.16]{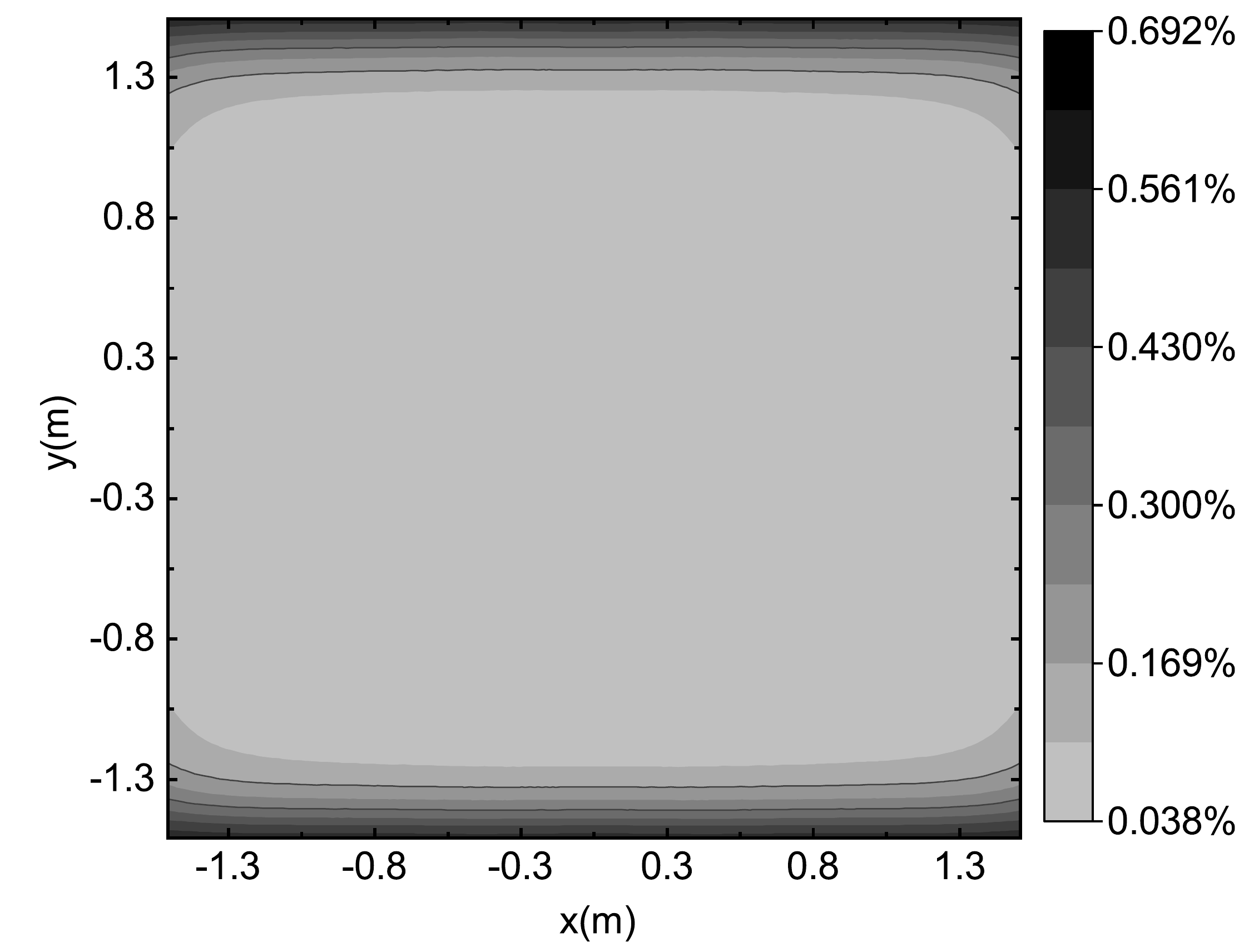}}
	\caption{The relative difference between the PD-BEM results the classical BEM results under the second loading condition with
		the Gauss kernel: ${\rm h_r} = 1/100 {\rm m}$(a); ${\rm h_r} = 1/400 {\rm m}$(b); ${\rm h_r} = 1/1600 {\rm m}$(c).}
	\label{fig16}
\end{figure}

\begin{figure}[!htb]
	\centering
	\subfigure[]{
		\label{fig27a}
		\includegraphics[scale=0.16]{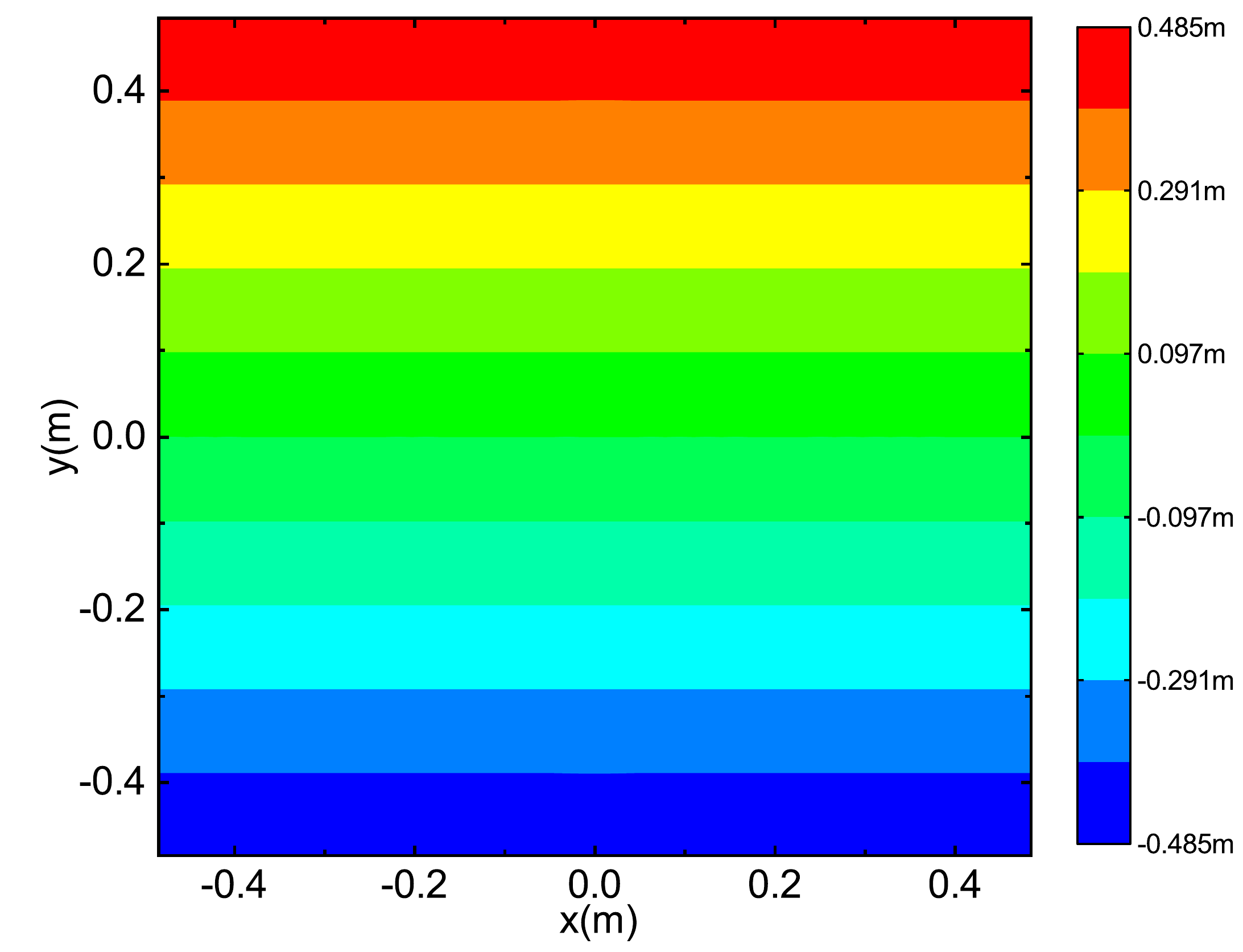}}
	\hspace{0.01in}
	\subfigure[]{
		\label{fig27b}
		\includegraphics[scale=0.16]{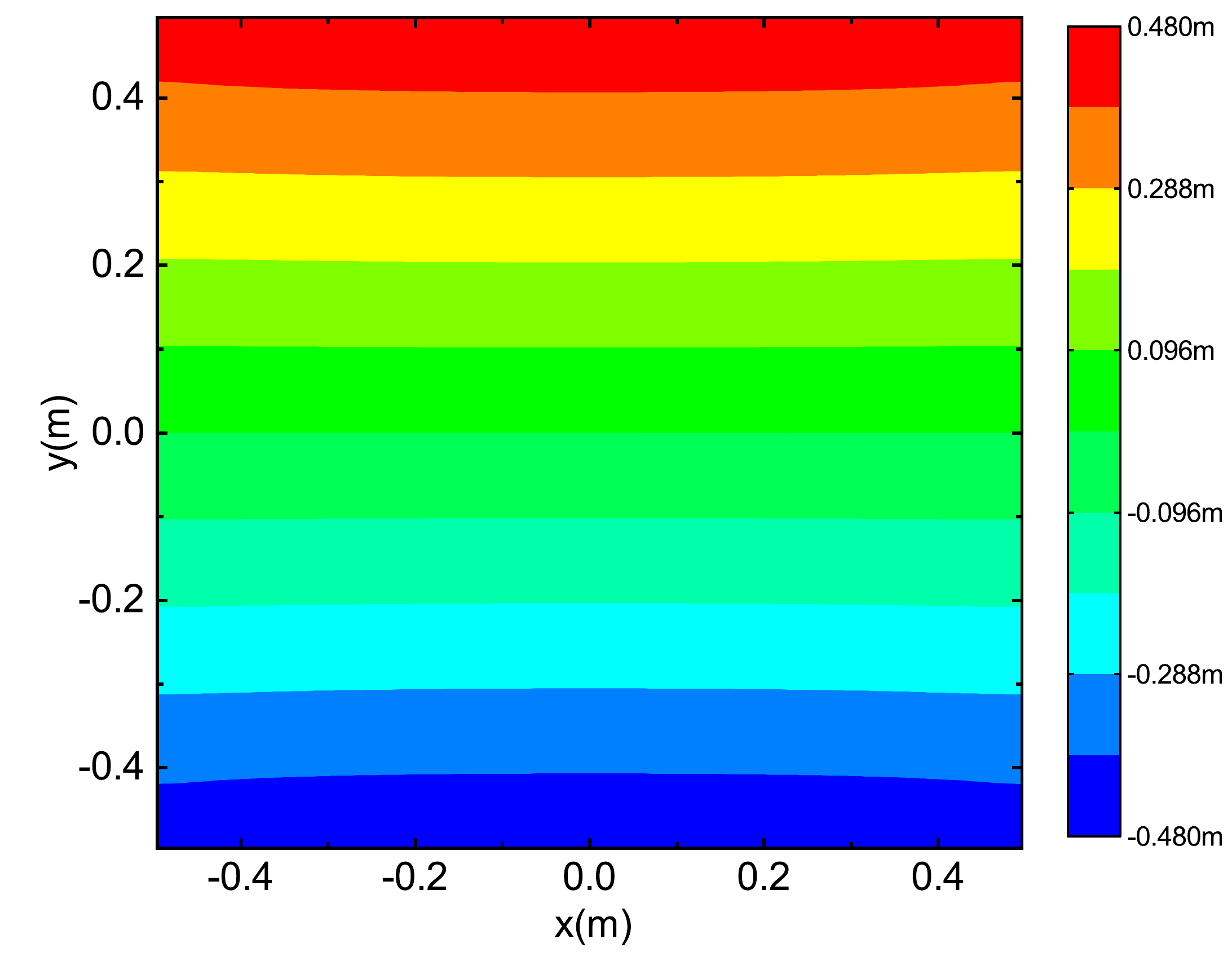}}
	\hspace{0.01in}
	\subfigure[]{
		\label{fig27c}
		\includegraphics[scale=0.16]{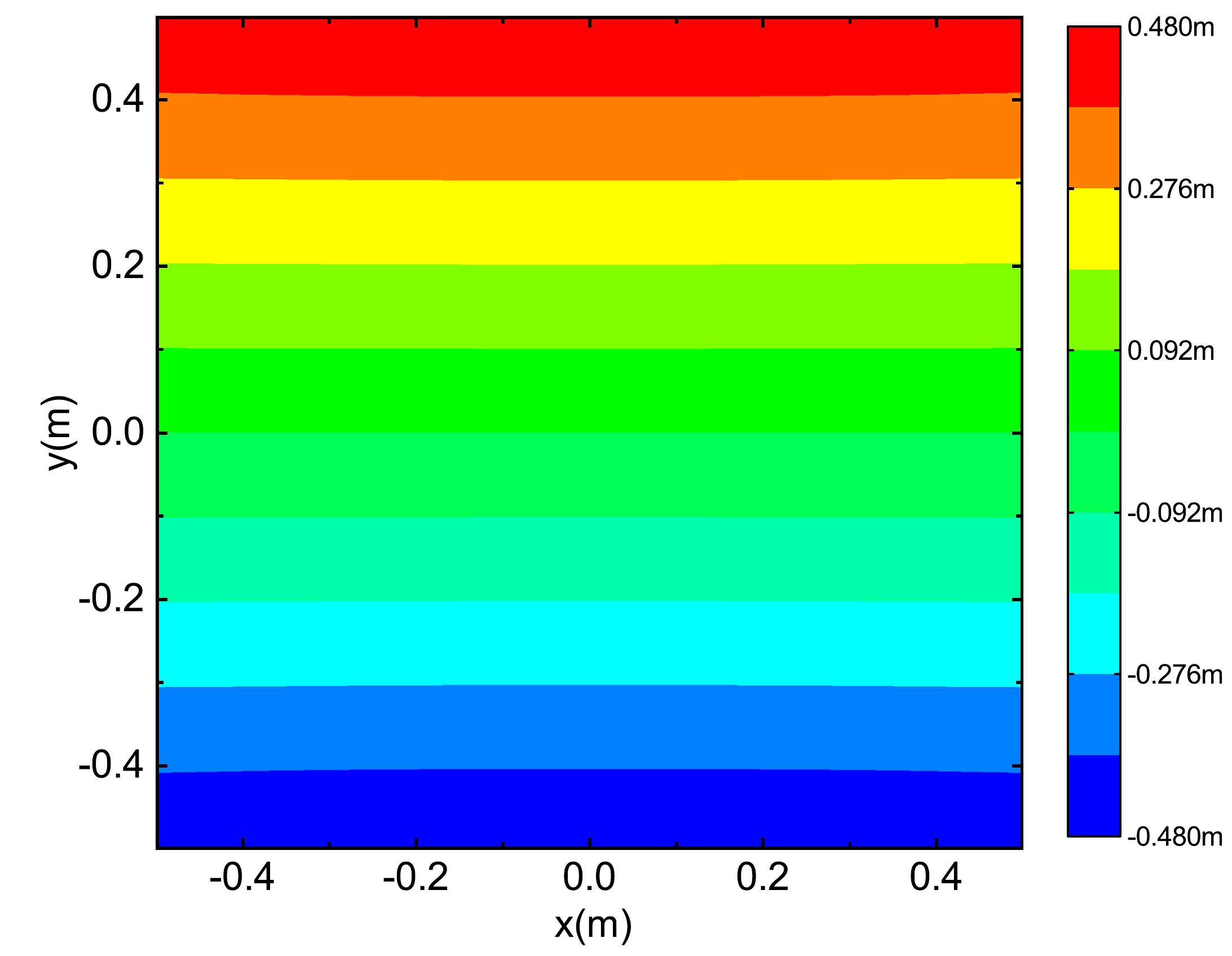}}
	\caption{The displacement under the third loading condition: PD-BEM(a); PD-MPM with $120\times120$ particles(b); PD-MPM with $360\times360$ particles(c).}
	\label{fig27}
\end{figure}

In Tables \ref{tab2} and \ref{tab3}, $E$ is the elastic modulus; $\nu$ is the Poisson ratio; $\rho$ is the density; ${\rm a}$ is the length of the sides; ${\rm t}_{n}$ and ${\rm t}_{\tau}$ are the normal and tangential surface tractions, respectively; ${\rm u}_{n}$ and ${\rm u}_{\tau}$ are the normal and tangential displacements, respectively; ${\rm b}_{x}$ and ${\rm b}_{y}$ are the body force densities in the $x$- and $y$-directions, respectively. For the first and third loading conditions, we only calculate a quarter of the domain due to the symmetry. For the second loading condition, we calculate the whole domain. For the first two loading conditions, we compute the displacements with the constant kernel and the Gauss kernel in (\ref{faaa24a}), when the characteristic length ${\rm h_r}$ are $1/100 {\rm m}$, $1/400 {\rm m}$ and $1/1600 {\rm m}$. Figure \ref{fig13} and  Figure \ref{fig14} show the relative differences between the displacements given by the PD-BEM and the classical BEM (C-BEM) in the $x$-direction, which is defined as $\left|u_x^{\rm{PD-BEM}}-u_x^{\rm{C-BEM}}\right|\times 100\%/u_x^{\rm{C-BEM}}$,  under the first loading condition for the the constant kernel and the Gauss kernel, respectively. Figure \ref{fig15} and Figure \ref{fig16} display the results under the second loading condition. For the third loading condition, we  compute the  displacement  in the $y$-direction with the constant value kernel for ${\rm h_r} = 1/60 {\rm m}$ by the PD-BEM and the PD-MPM, and the results are shown in Figure \ref{fig27}.

Figures~\ref{fig13}-\ref{fig16} show that the results obtained by using the PD-BEM can converge to the ones of the local boundary element method for the considered kernel functions, when the characteristic length tends to zero. Now we focus on the analysis of Figure \ref{fig27}. Because the results of different kernel functions are similar (Figure~\ref{fig13}-Figure~\ref{fig16}), we only calculate the displacement with the constant kernel function. For the PD-BEM, we take $262$ seconds to get a result with $32$ elements on each side, and the results do not exhibit the boundary softening phenomenon which is typically seen in the PD-MPM, as seen from the results in Figure~\ref{fig27b} which are obtained with $120\times120$ particles. The phenomenon can be alleviated with a finer grid, but the computational cost is increased. It takes $14786$ seconds to get a result with a negligible boundary softening phenomenon~\cite{tb73} with $360\times360$ particles, as seen from the results in Figure~\ref{fig27c}. In this case, the ratio of the times spent by the PD-BEM and by the PD-MPM is $1.8\%$.

\subsection{Vibration of a bar}\label{atl42}

\begin{figure}[!htb]
	\centerline{\includegraphics[scale=0.3]{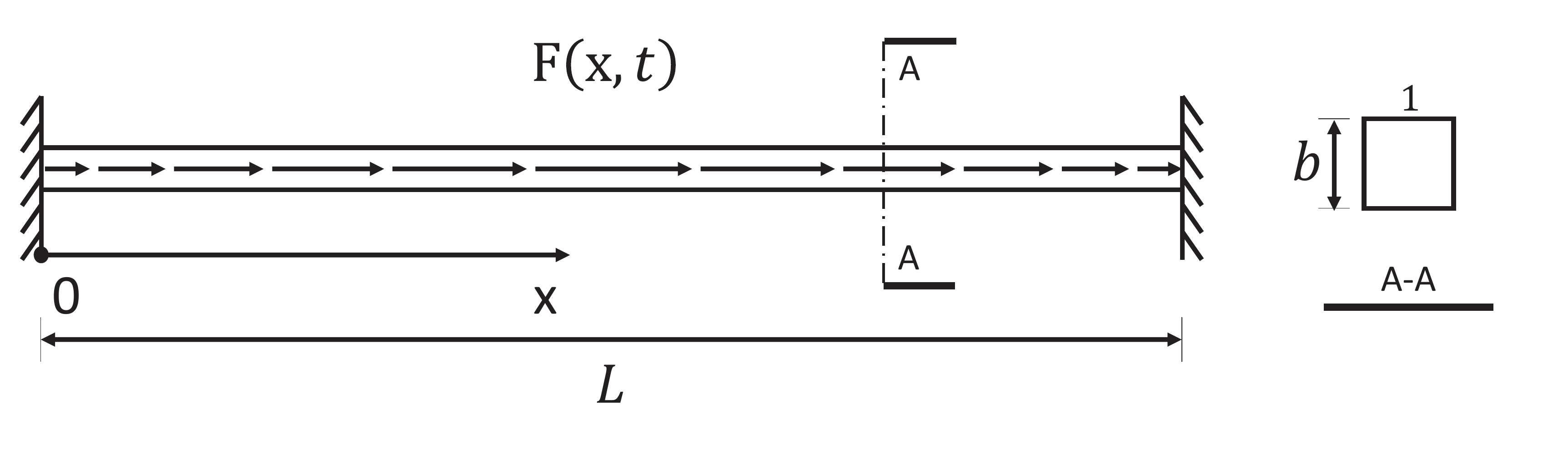}}	
	\caption{A clamped bar subjected to axial dynamic harmonic loading.\label{sc4}}
\end{figure}

We consider the vibration of a bar with two clamped ends subjected to axial dynamic harmonic loading, shown in Figure \ref{sc4}. The geometric and material parameters are listed in Table \ref{tab6}. The bar has a unit thickness, so it is treated as a two-dimensional rectangular region with a length $L$ and a height $b$.
\begin{center}
	\begin{table}[!htb]
		\centering
		\caption{Geometric and material parameters of the bar.\label{tab6}}
		\begin{tabular}{cccccc}
			\hline
			$\rho \left({\rm kg/m^2}\right)$&$E \left({\rm Pa}\right)$&$\nu$&$L \left({\rm m}\right)$&$b \left({\rm m}\right)$&${\rm F} \left(x,t\right) \left({\rm N/m}\right)$\\
			\hline
			$1.0$&$1.0$&$1/3$&$1.0$&$0.05$&${\rm sin}\left(t\right){\rm sin}\left(2x\pi\right)$\\
			\hline
		\end{tabular}
	\end{table}
\end{center}

\begin{figure}[!htb]
	\centering
	\subfigure[]{
		\label{fig23a}
		\includegraphics[scale=0.16]{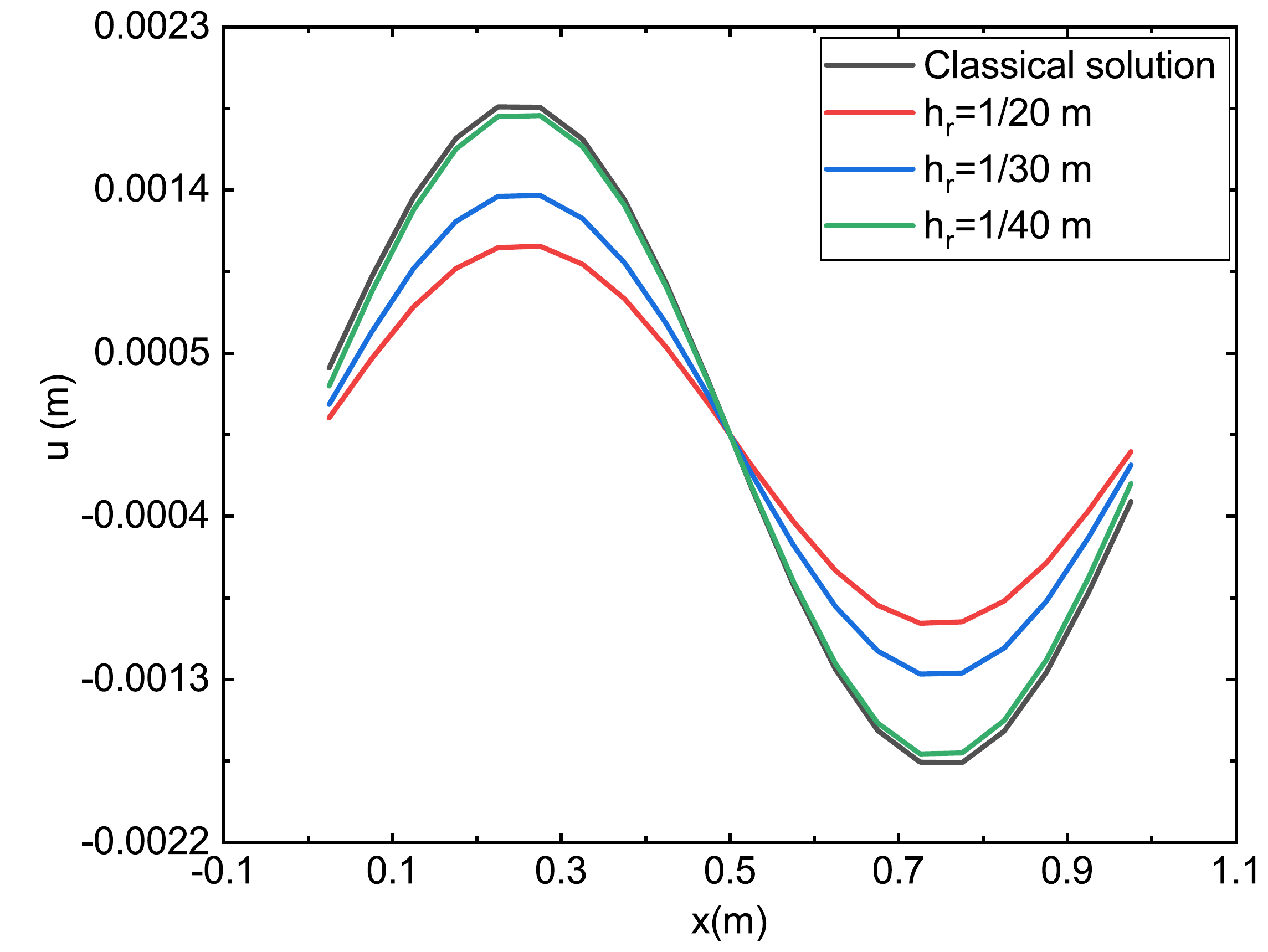}}
	\hspace{0.01in}
	\subfigure[]{
		\label{fig23b}
		\includegraphics[scale=0.16]{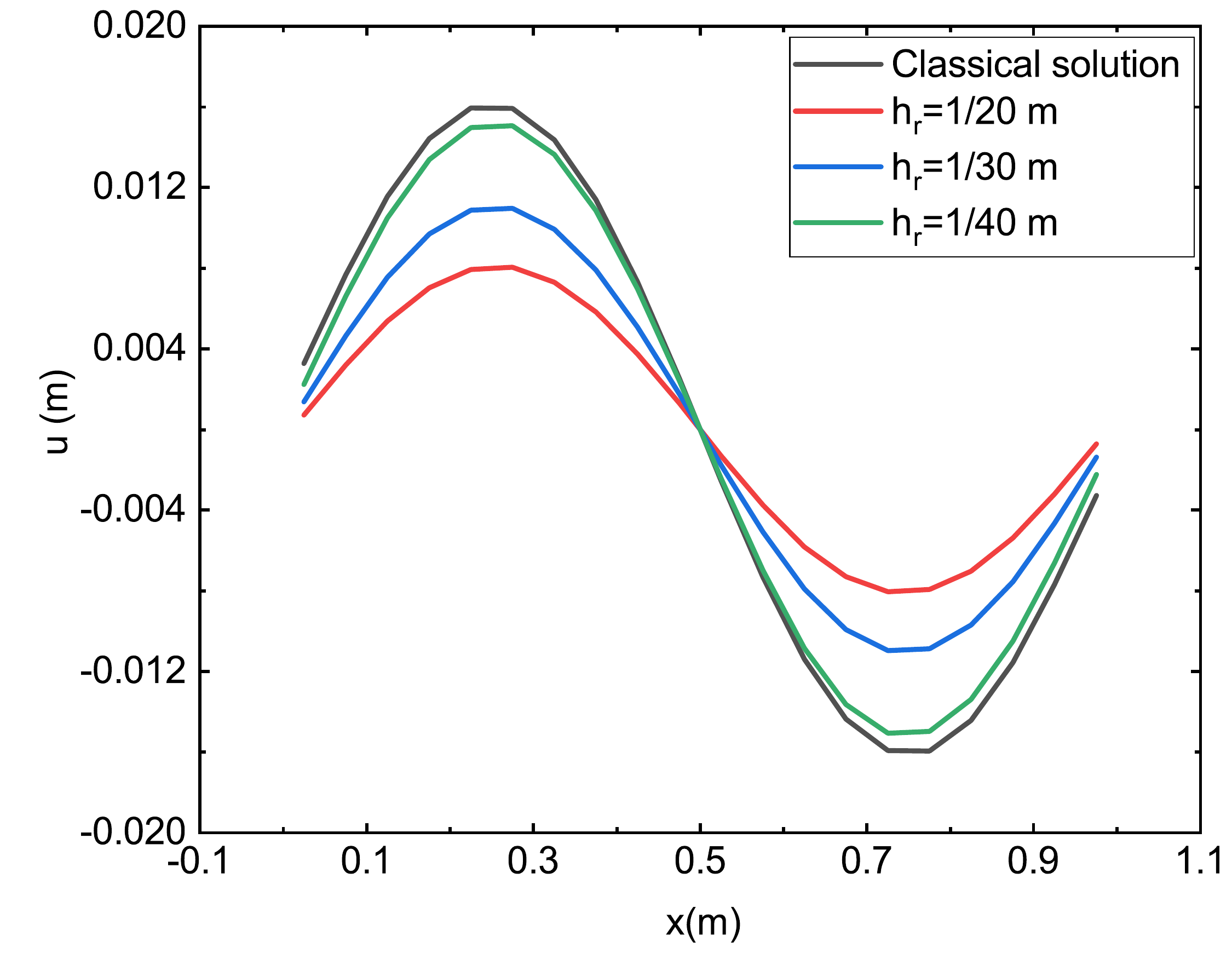}}	
	\hspace{0.01in}
	\subfigure[]{
		\label{fig23c}
		\includegraphics[scale=0.16]{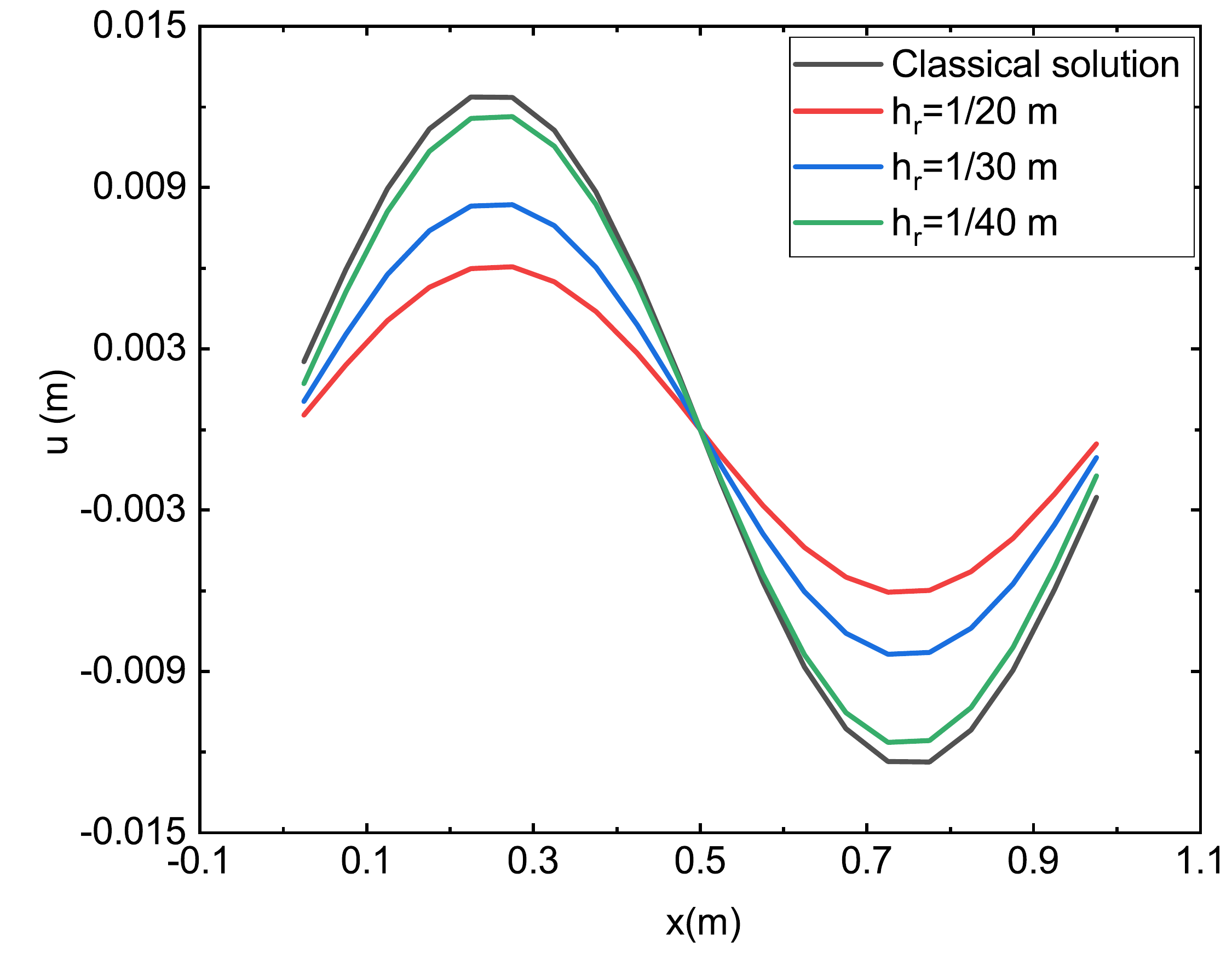}}	
	\vfill
	\subfigure[]{
		\label{fig23d}
		\includegraphics[scale=0.16]{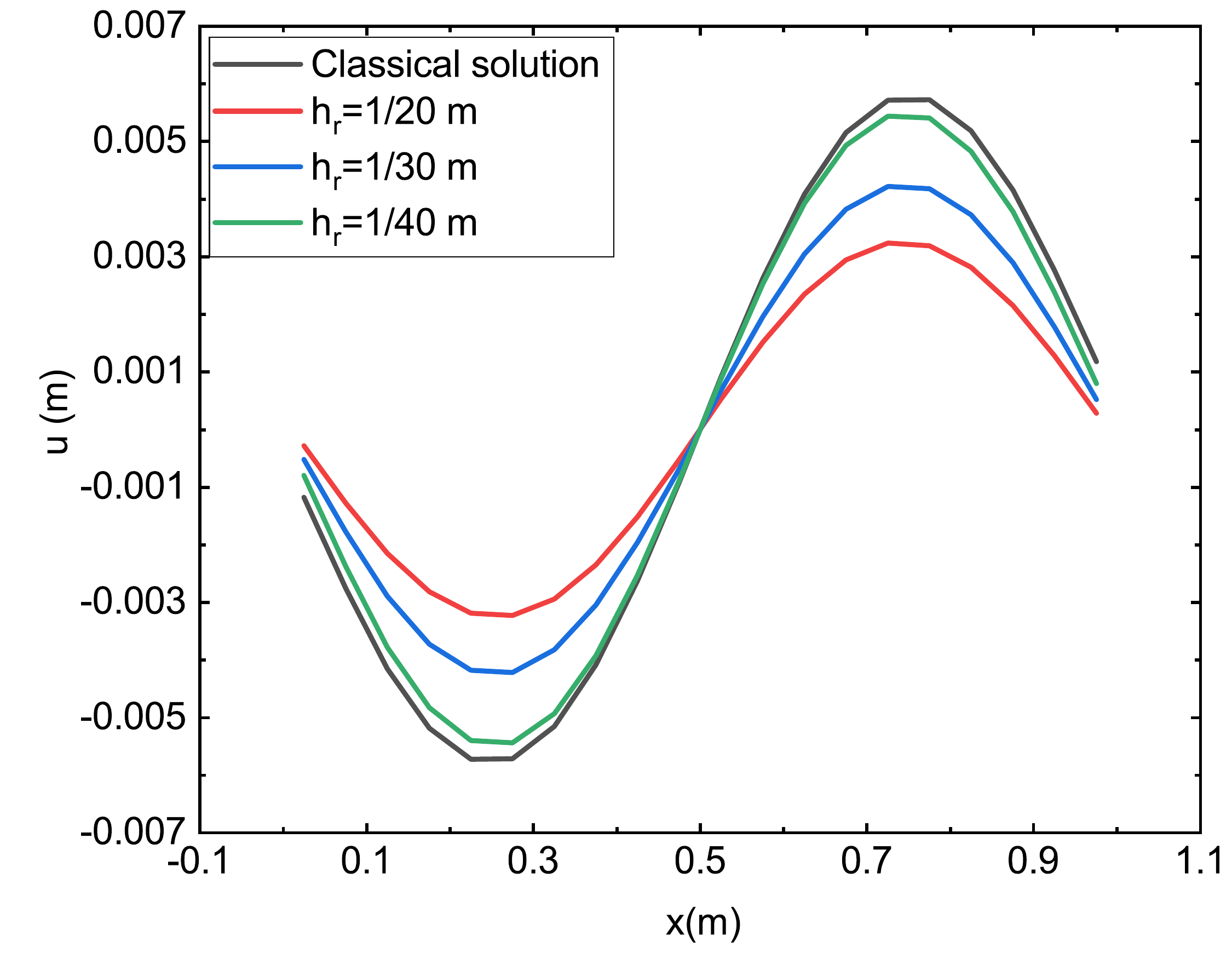}}
	\hspace{0.01in}
	\subfigure[]{
		\label{fig23e}
		\includegraphics[scale=0.16]{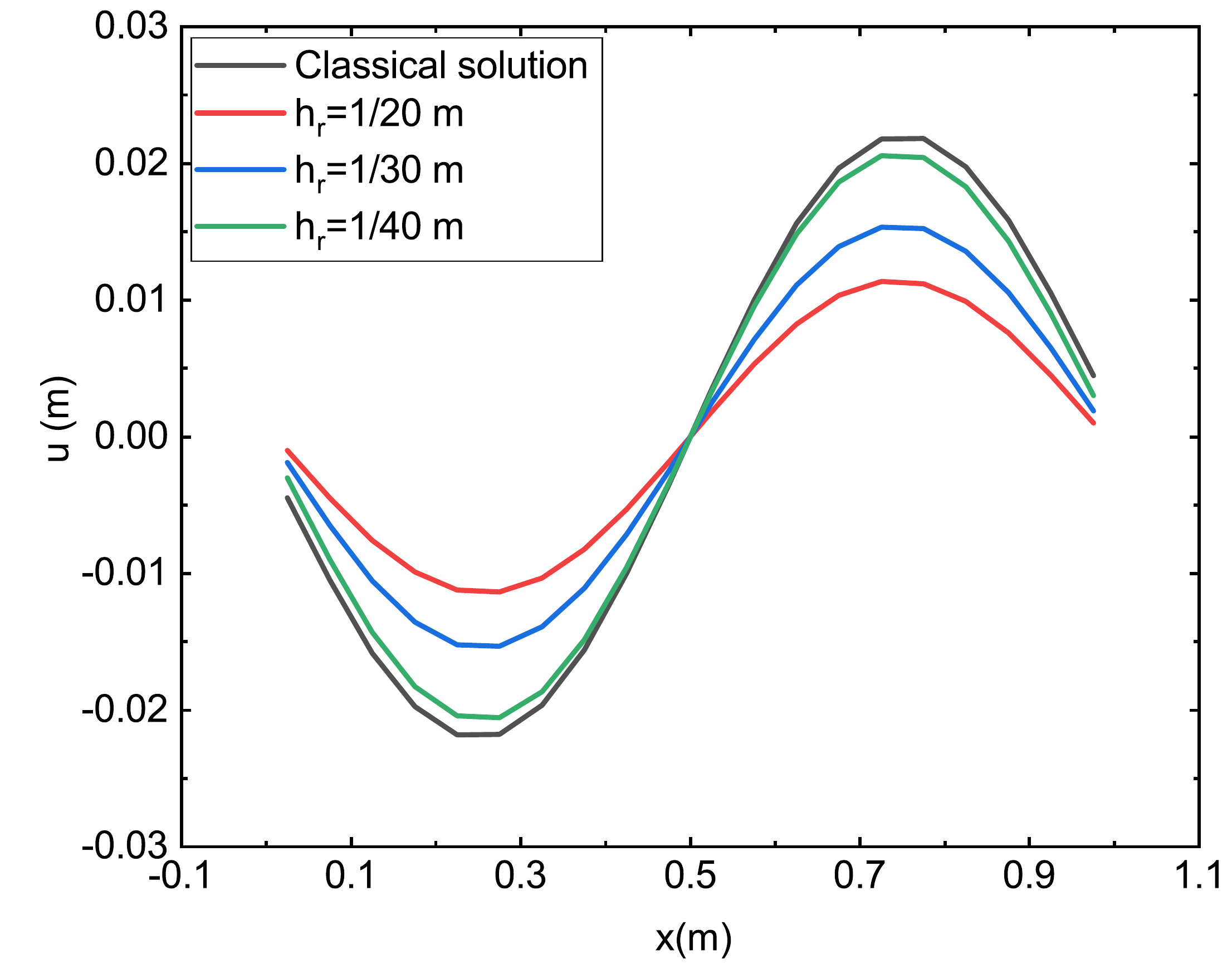}}	
	\hspace{0.01in}
	\subfigure[]{
		\label{fig23f}
		\includegraphics[scale=0.16]{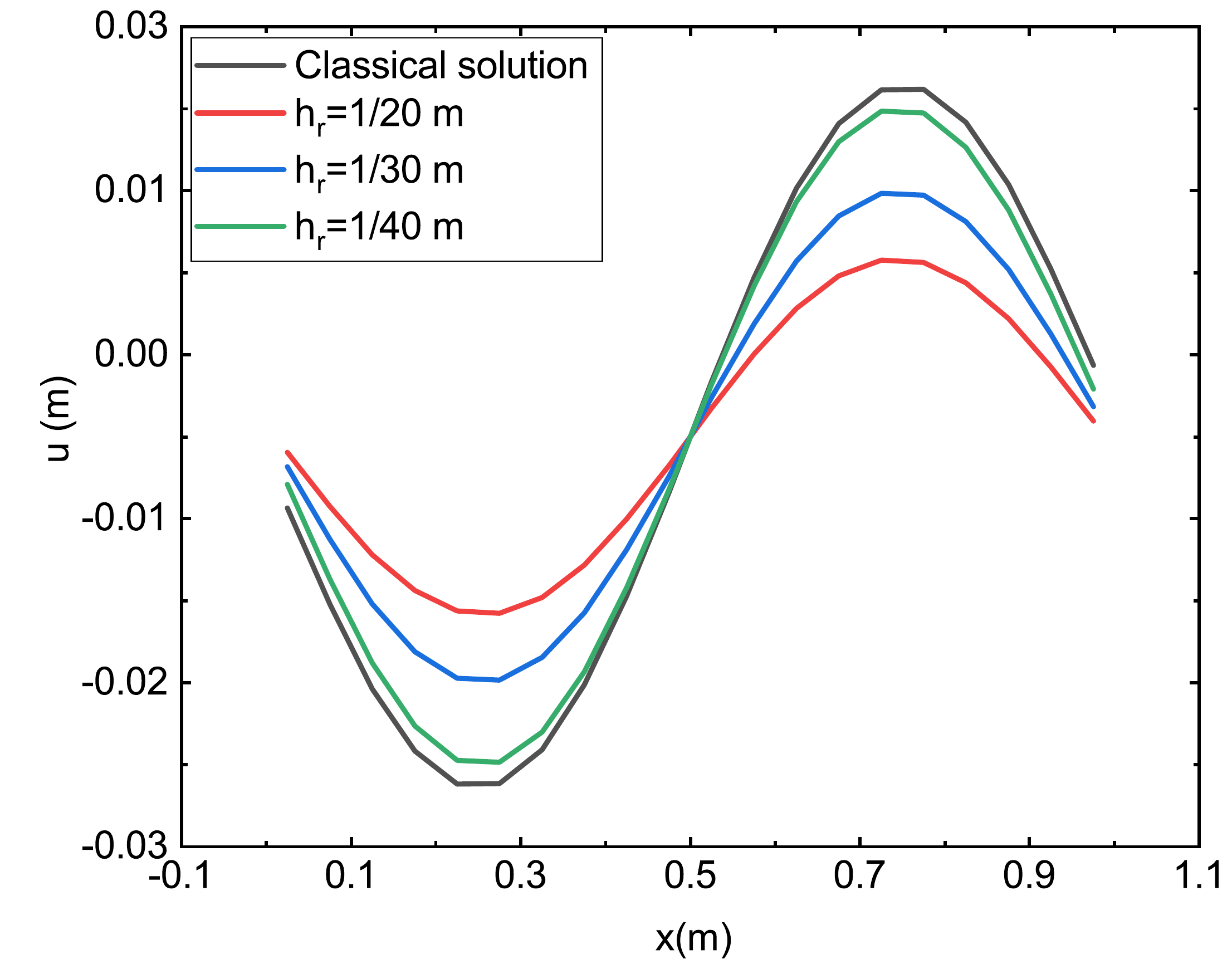}}	
	\caption{The average deflection with the constant value kernel at different moments of time: $0.25 {\rm s}$(a); $1.25 {\rm s}$(b); $2.25 {\rm s}$(c); $3.25 {\rm s}$(d); $4.25 {\rm s}$(e); $5.25 {\rm s}$(f).}
	\label{fig23}	
\end{figure}

\begin{figure}[!htb]
	\centering
	\subfigure[]{
		\label{fig24a}
		\includegraphics[scale=0.16]{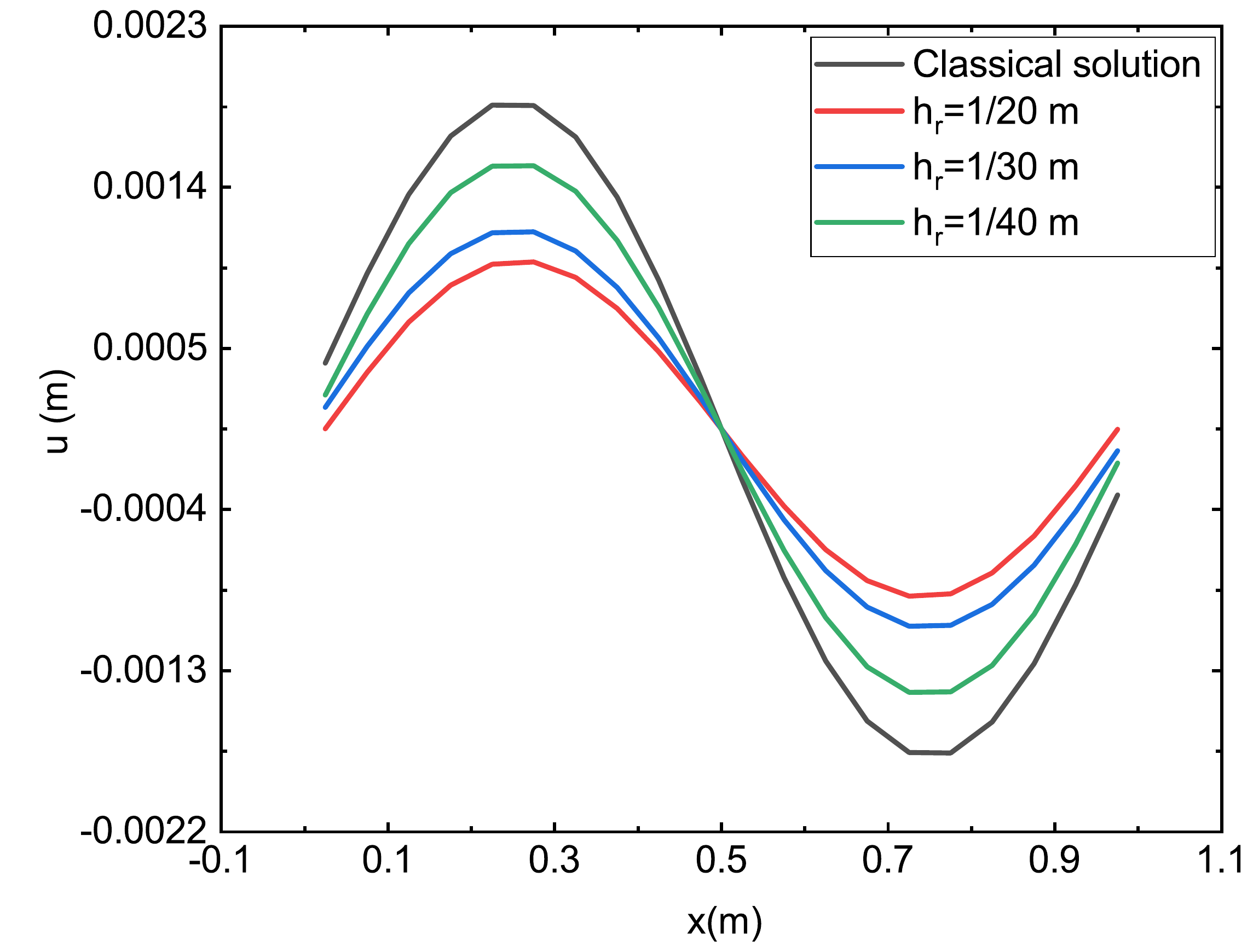}}
	\hspace{0.01in}
	\subfigure[]{
		\label{fig24b}
		\includegraphics[scale=0.16]{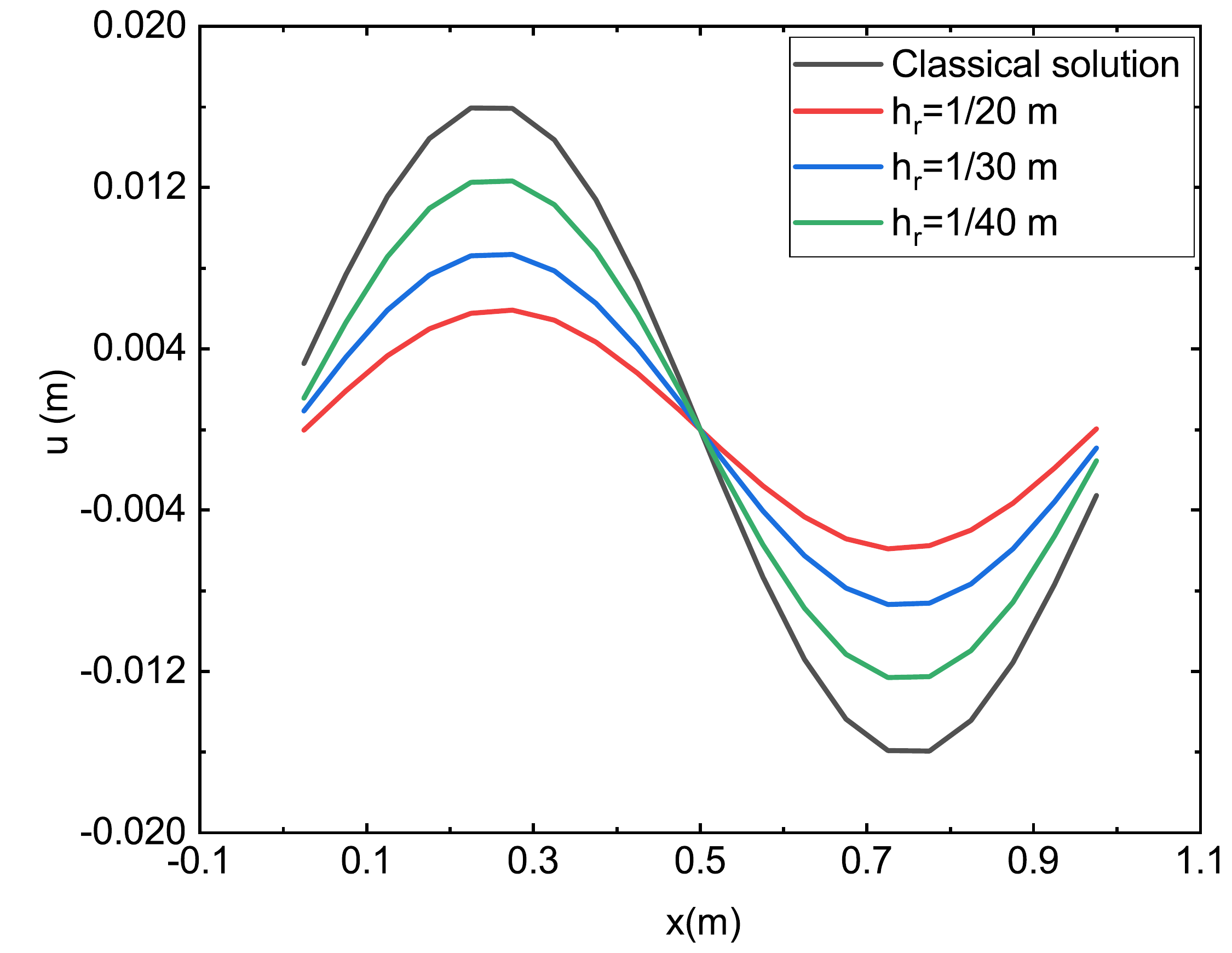}}	
	\hspace{0.01in}
	\subfigure[]{
		\label{fig24c}
		\includegraphics[scale=0.16]{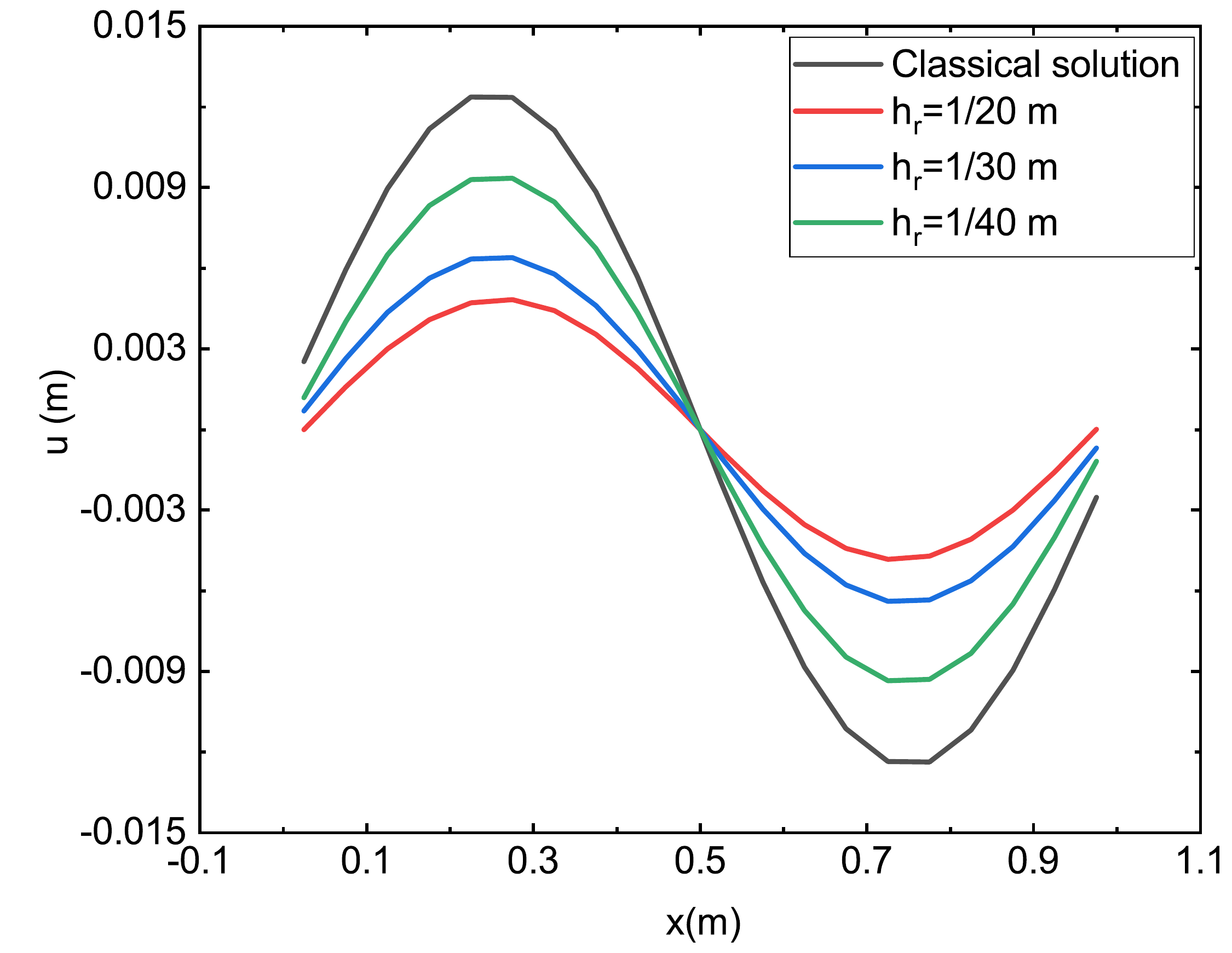}}	
	\vfill
	\subfigure[]{
		\label{fig24d}
		\includegraphics[scale=0.16]{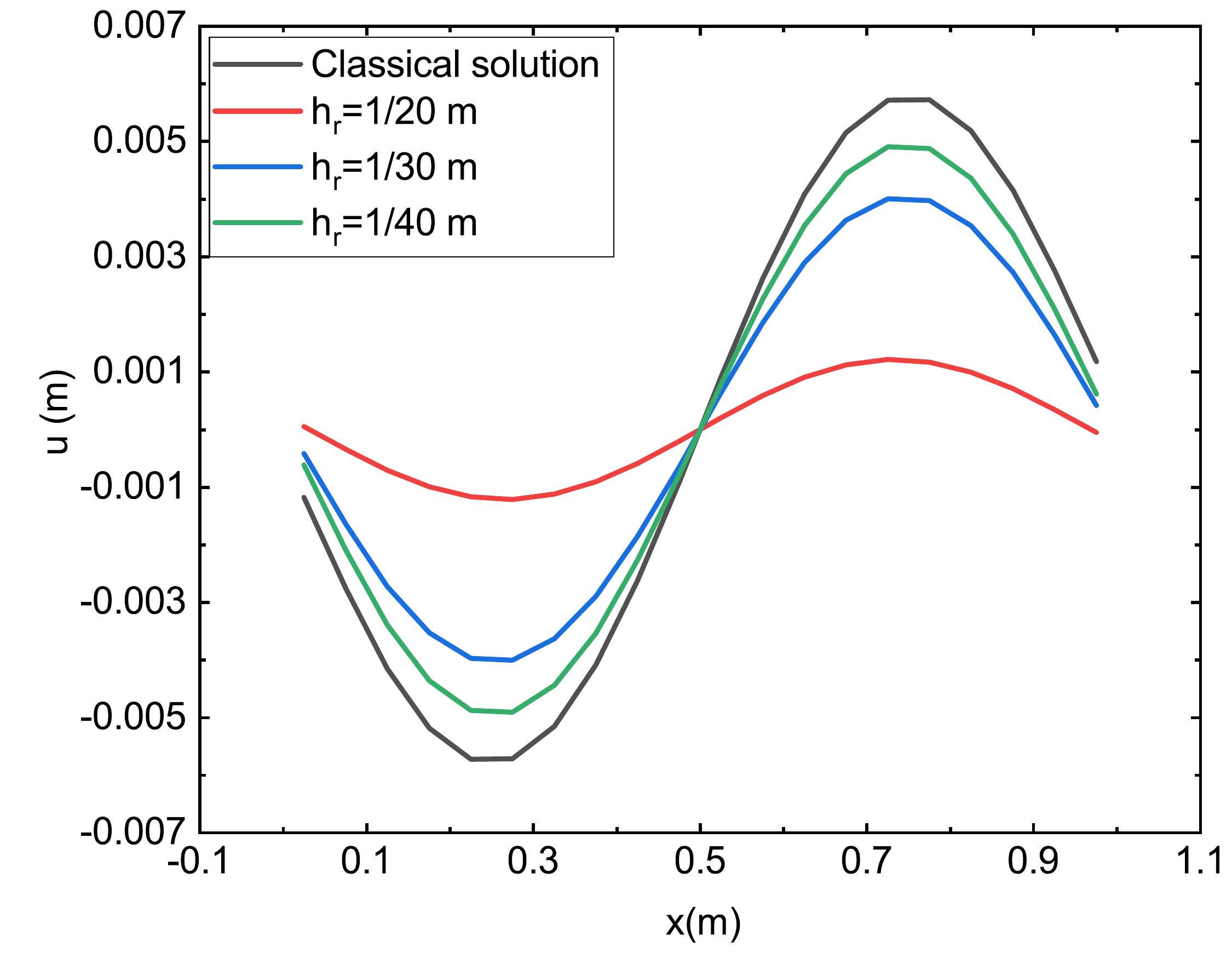}}
	\hspace{0.01in}
	\subfigure[]{
		\label{fig24e}
		\includegraphics[scale=0.16]{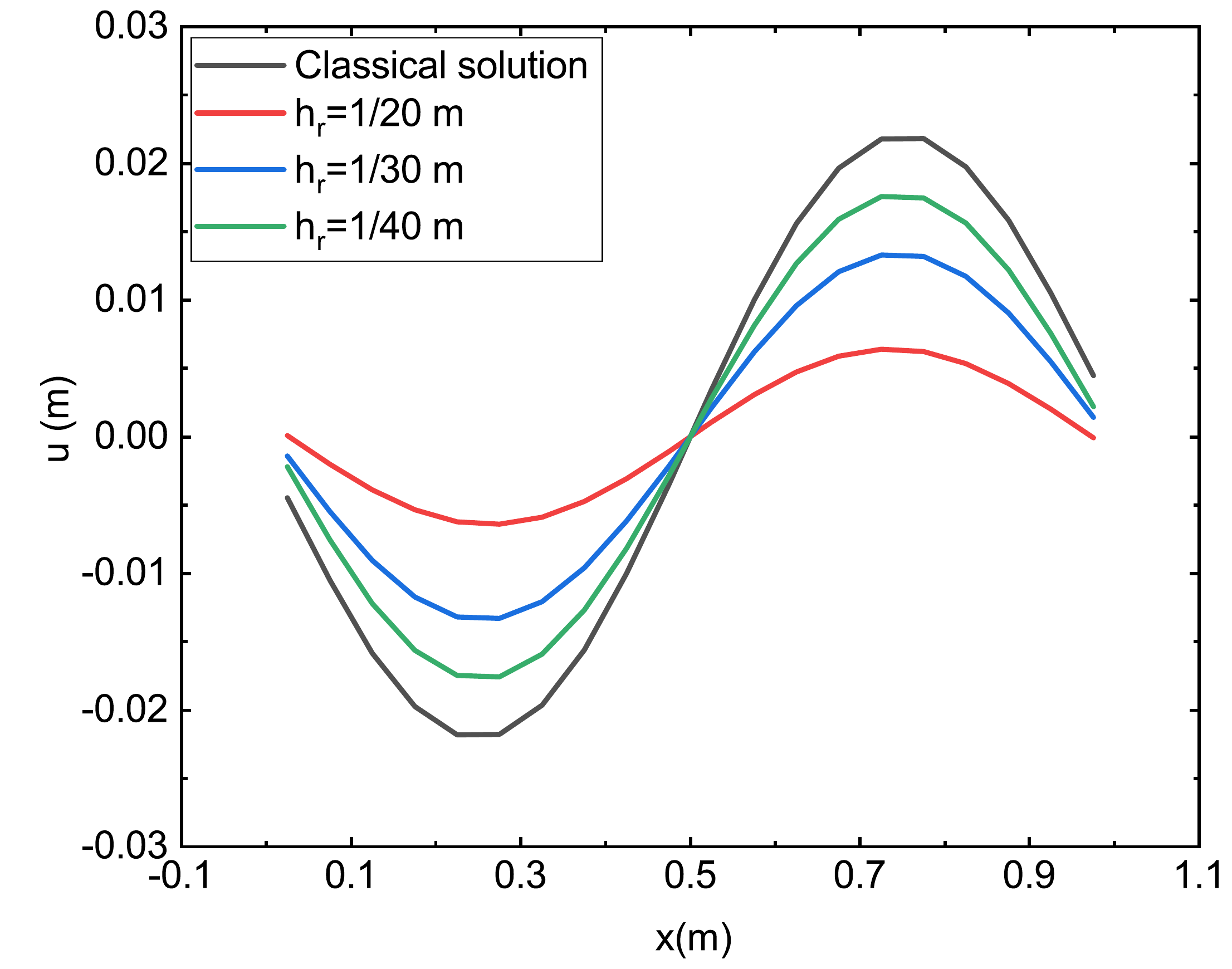}}	
	\hspace{0.01in}
	\subfigure[]{
		\label{fig24f}
		\includegraphics[scale=0.16]{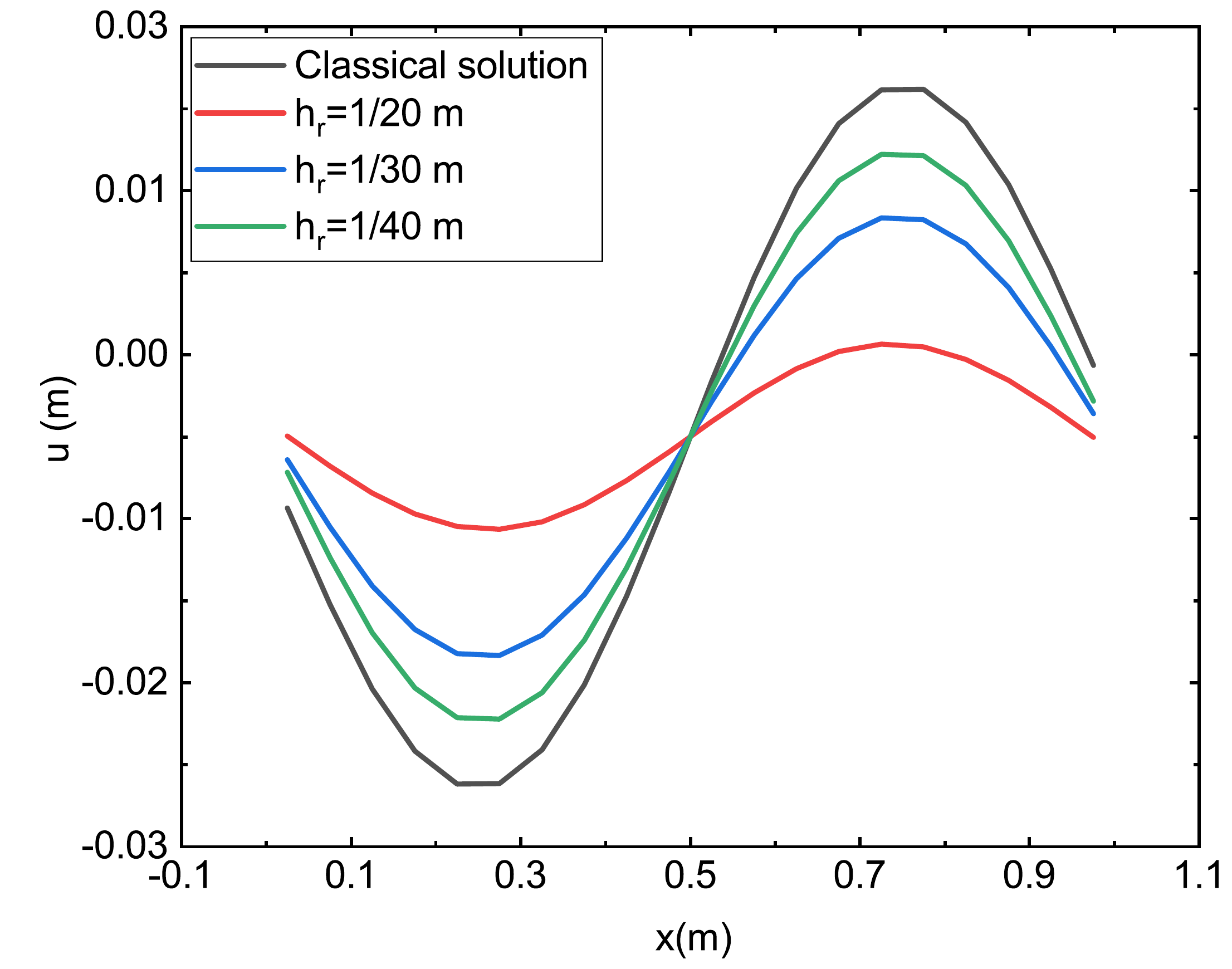}}	
	\caption{The average deflection with the Gauss kernel at different moments of time: $0.25 {\rm s}$(a); $1.25 {\rm s}$(b); $2.25 {\rm s}$(c); $3.25 {\rm s}$(d); $4.25 {\rm s}$(e); $5.25 {\rm s}$(f).}
	\label{fig24}	
\end{figure}

In Table \ref{tab6}, ${\rm F}\left(x,t\right)$ is the axial dynamic harmonic load. We examine the overall dynamic behaviour of the bar by calculating the average  axial displacement over the height. The results obtained  at different times for different characteristic lengths of $1/20 {\rm m}$, $1/30 {\rm m}$,$1/40 {\rm m}$ are shown in Figure \ref{fig23} (the constant value kernel) and Figure \ref{fig24} (the Gauss kernel). Again, the PD-BEM results approach those of the classical BEM as the characteristic length approaches zero.
\begin{figure}[!htb]
	\centerline{\includegraphics[scale=0.3]{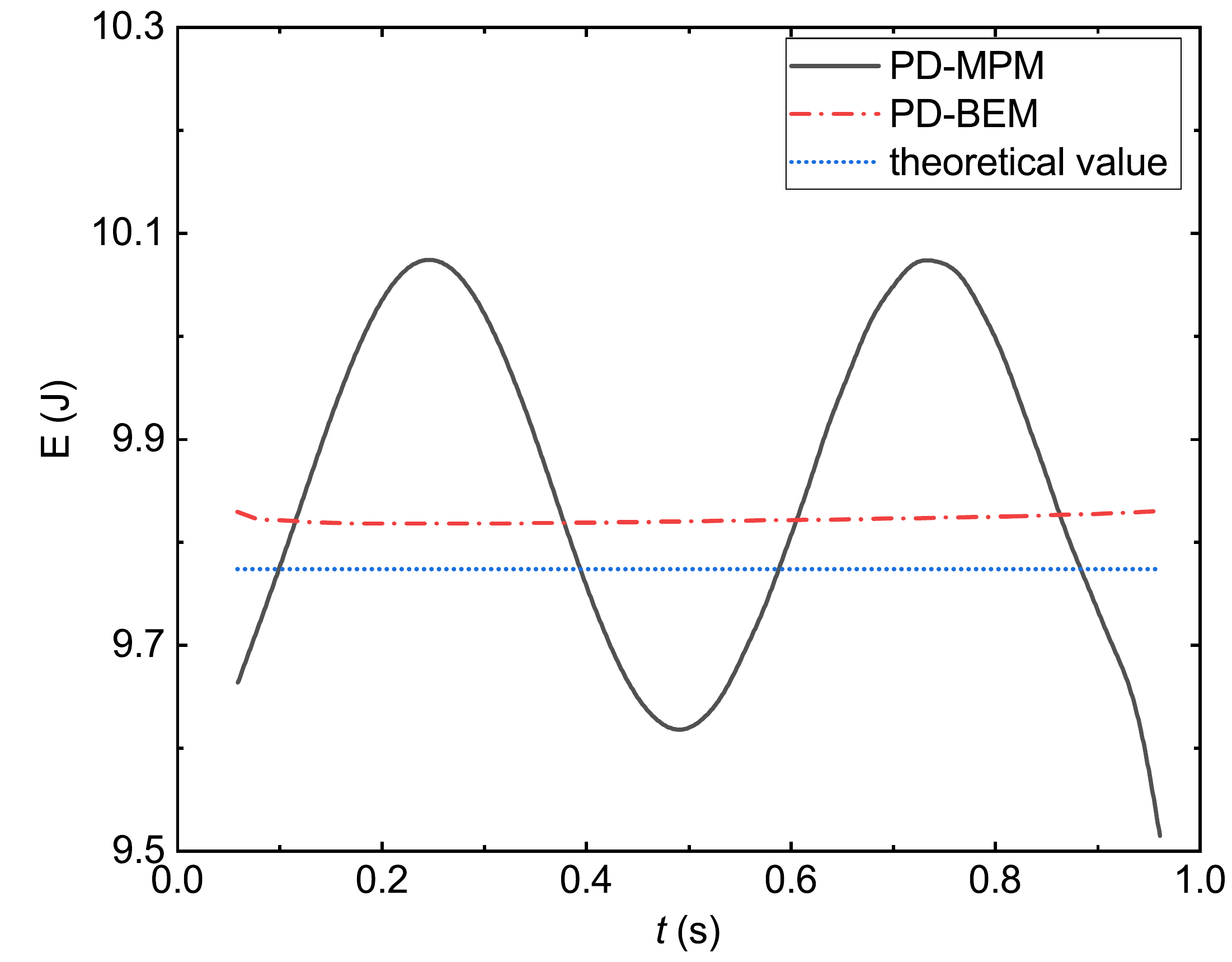}}	
	\caption{The total energy  of the bar versus the time $t$.\label{fig28}}	
\end{figure}
To highlight the advantage of the PD-BEM, we calculate the total energy of the bar by the PD-BEM and the PD-MPM, and compare them with the theoretical value in Figure \ref{fig28}. The total energy is defined as
\begin{sequation}\label{faaaa38}
\begin{aligned}
{ E_{total}} =& \dfrac{1}{4} \int_{\Omega} \int_{\Omega} \left[\mathbf{u} \left(\mathbf{x}',t\right) - \mathbf{u} \left(\mathbf{x},t\right)\right] \cdot \mathbf{C} \left(\mathbf{x}' - \mathbf{x}\right) \cdot \left[\mathbf{u} \left(\mathbf{x}',t\right) - \mathbf{u} \left(\mathbf{x},t\right)\right] {\rm d}V_{\mathbf{x}'} {\rm d}V_{\mathbf{x}} \\
& + \dfrac{1}{2} \int_{\Omega} \rho \left(\mathbf{x},t\right) \dot{\mathbf{u}} \left(\mathbf{x},t\right)\cdot \dot{\mathbf{u}} \left(\mathbf{x},t\right) {\rm d}V_{\mathbf{x}}
\end{aligned}
\end{sequation}
where $\Omega$ is the solution domain; $\mathbf{u}$ is the displacement; $\mathbf{C}$ is the micromodulus tensor; $\rho$ is the density. From Figure \ref{fig23} and Figure \ref{fig24}, we find that the results of different kernel functions are similar. Thus, we only calculate the total energy with the constant kernel function and the characteristic length ${\rm h_r}=1/100 {\rm m}$. There is no damping for this problem, so the mechanical energy must be conserved. However, the energy calculated by the PD-MPM has undulation, whereas the one by the PD-BEM is conserved. This result demonstrates that the PD-BEM does not produce time  accumulation error. Moreover, the relative difference from the theoretical value, which is calculated through the conservation of energy and equal to the initial energy, is $0.58\%$ for the PD-BEM and $2.64\%$ for the PD-MPM.

\subsection{Two-dimensional wave propagation}\label{atl43}

\begin{figure}[!htb]
	\centerline{\includegraphics[scale=0.3]{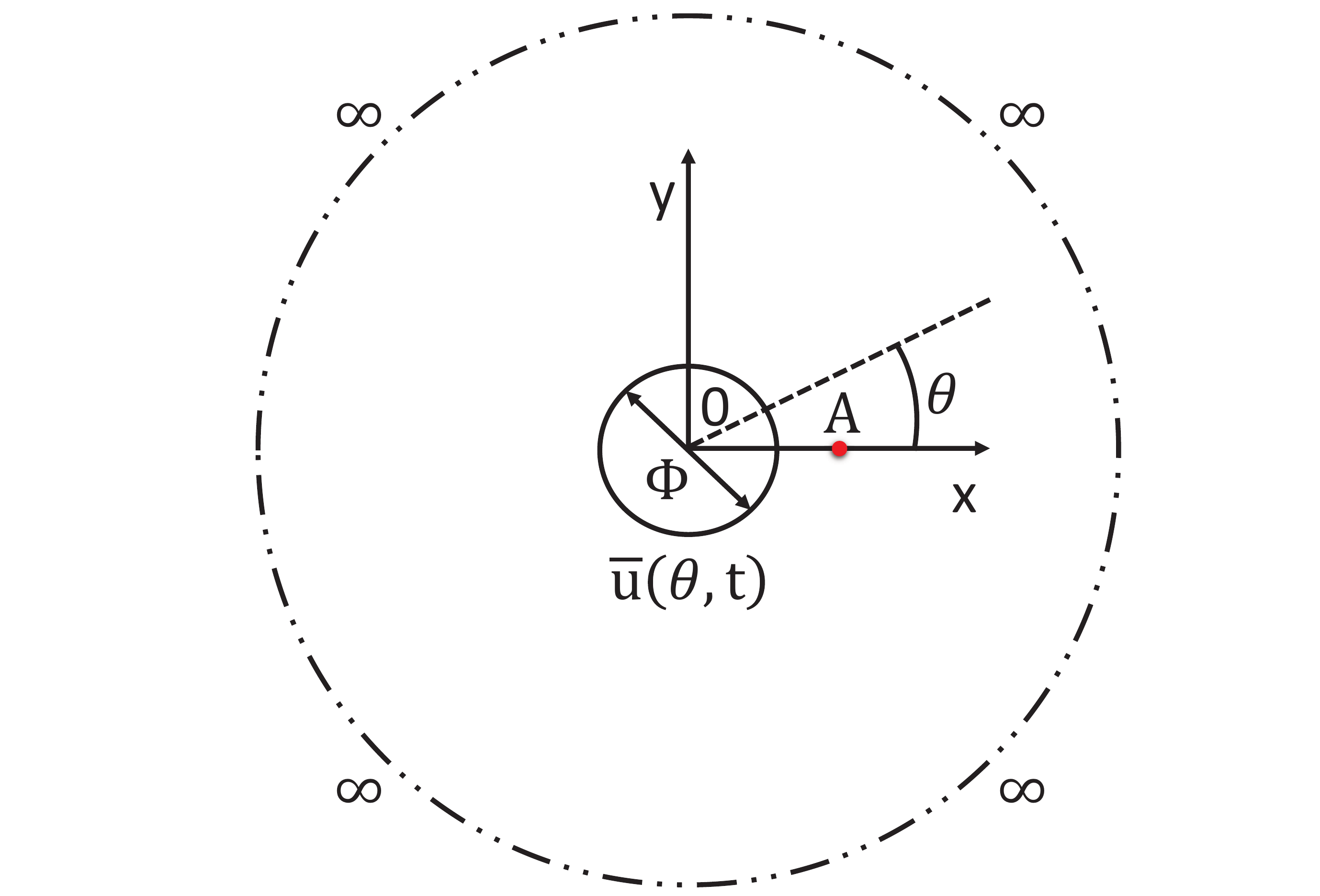}}	
	\caption{A two-dimensional infinite region with a circular void subjected to dynamic loading.\label{sc5}}
\end{figure}

Here, we calculate the two-dimensional wave propagation in an infinite domain where an excitation signal is imposed on the boundary of a circular void, as shown in Figure~\ref{sc5}. We can solve the problem in the infinite domain rather than a large finite area because of the advantage of the boundary element method, which does not need an artificial boundary. Therefore, we do not need to deal with the numerical reflection introduced by artificial boundaries~\cite{tb56,tb57}. The calculation parameters are given in Table \ref{tab7}.
\begin{center}
	\begin{table}[!htb]
		\centering
		\caption{The parameters for the two-dimensional wave propagation problem.\label{tab7}}
		\begin{tabular}{cccccc}
			\hline
			$\rho \left({\rm kg/m^2}\right)$&$E \left({\rm Pa}\right)$&$\nu$&$\Phi \left({\rm m}\right)$&$\overline{{\rm u}}_x \left(\theta,t\right) \left({\rm m}\right)$&$\overline{{\rm u}}_y \left(\theta,t\right) \left({\rm m}\right)$\\
			\hline
			$1.0$&$1.0$&$1/3$&$0.2$&${\rm sin}(20t) {\rm cos}(\theta)$&${\rm sin}(20t) {\rm sin}(\theta)$\\
			\hline
		\end{tabular}
	\end{table}
\end{center}

\begin{figure}[!htb]
	\centering
	\subfigure[]{
		\label{fig25a}
		\includegraphics[scale=0.16]{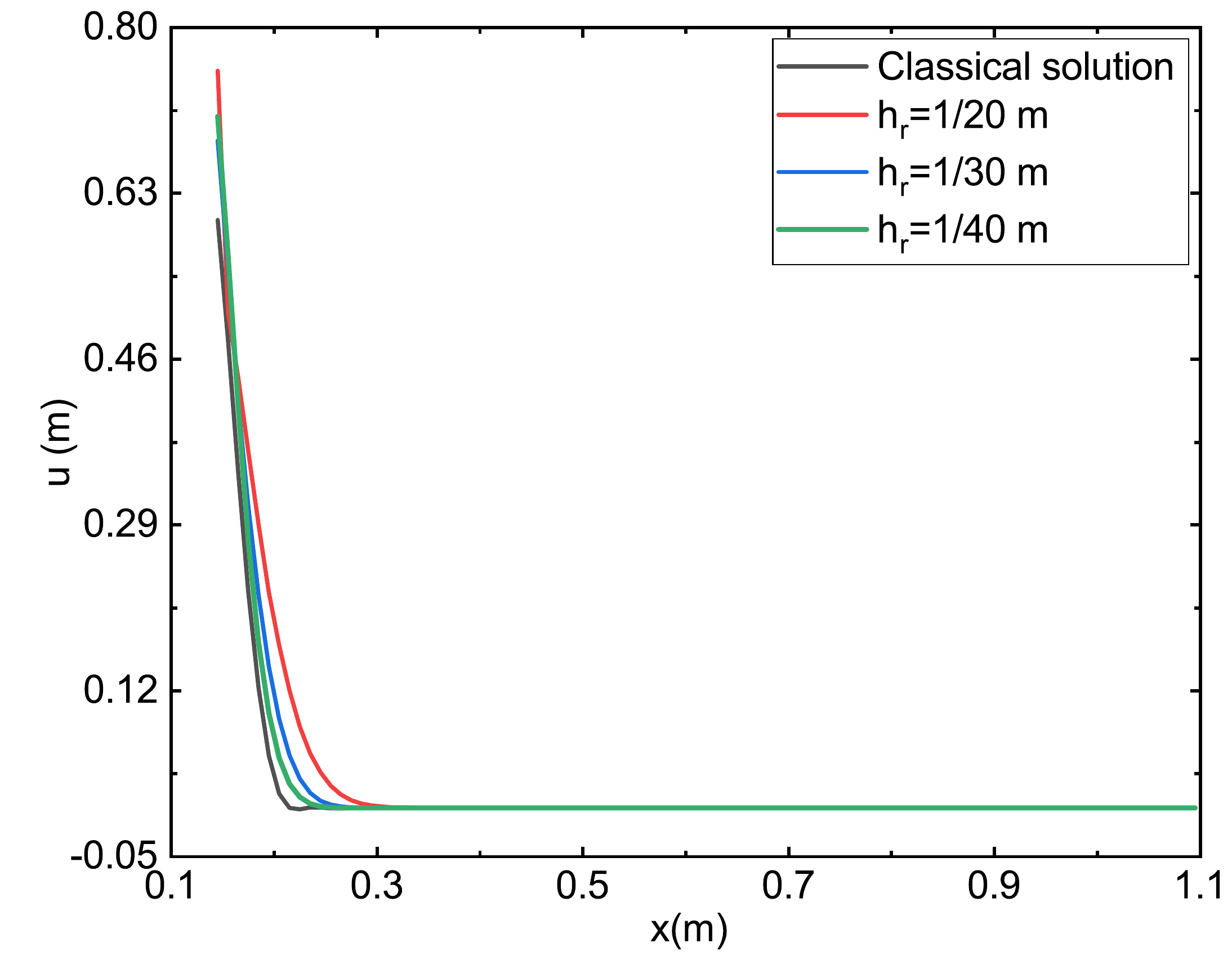}}
	\hspace{0.01in}
	\subfigure[]{
		\label{fig25b}
		\includegraphics[scale=0.16]{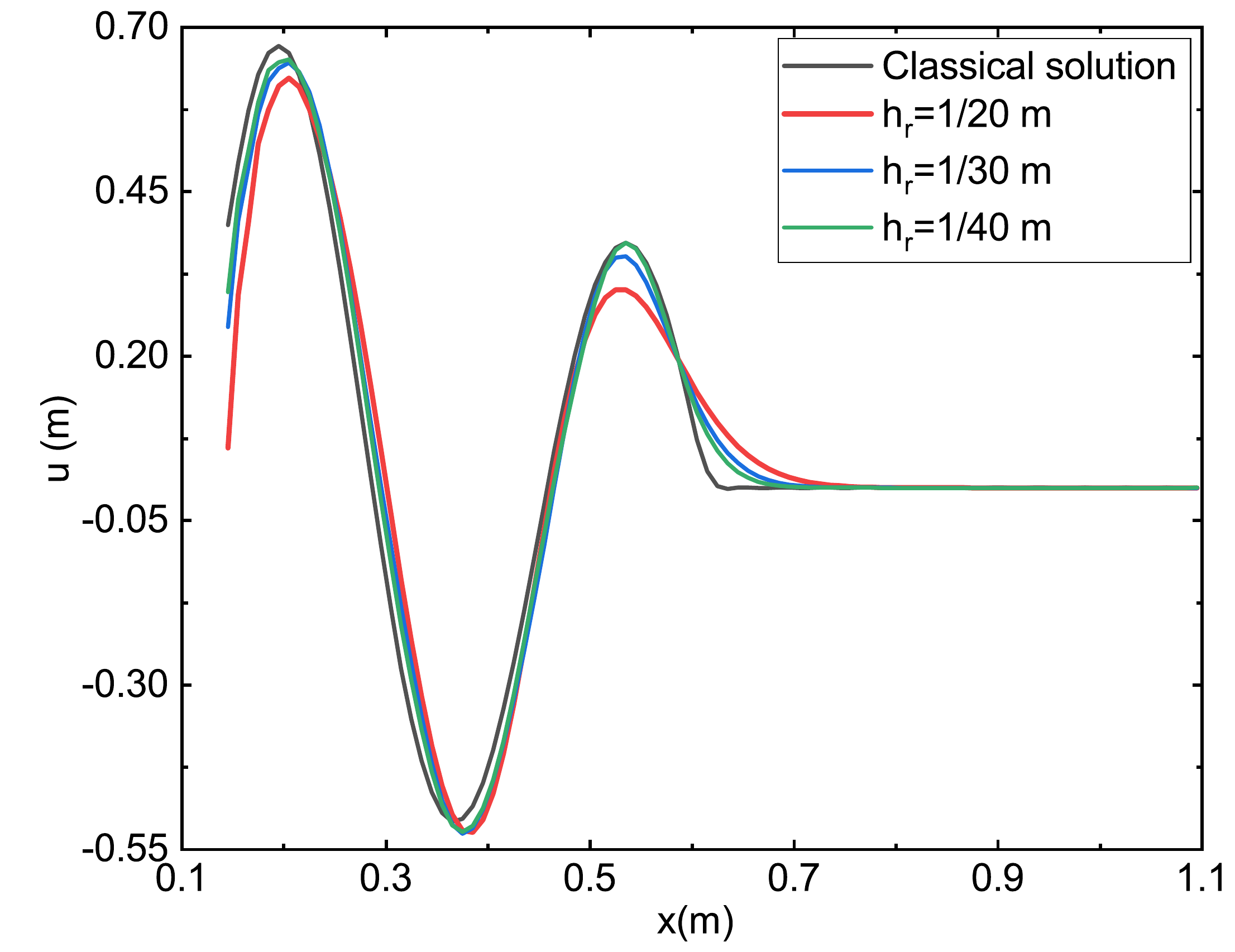}}	
	\hspace{0.01in}
	\subfigure[]{
		\label{fig25c}
		\includegraphics[scale=0.16]{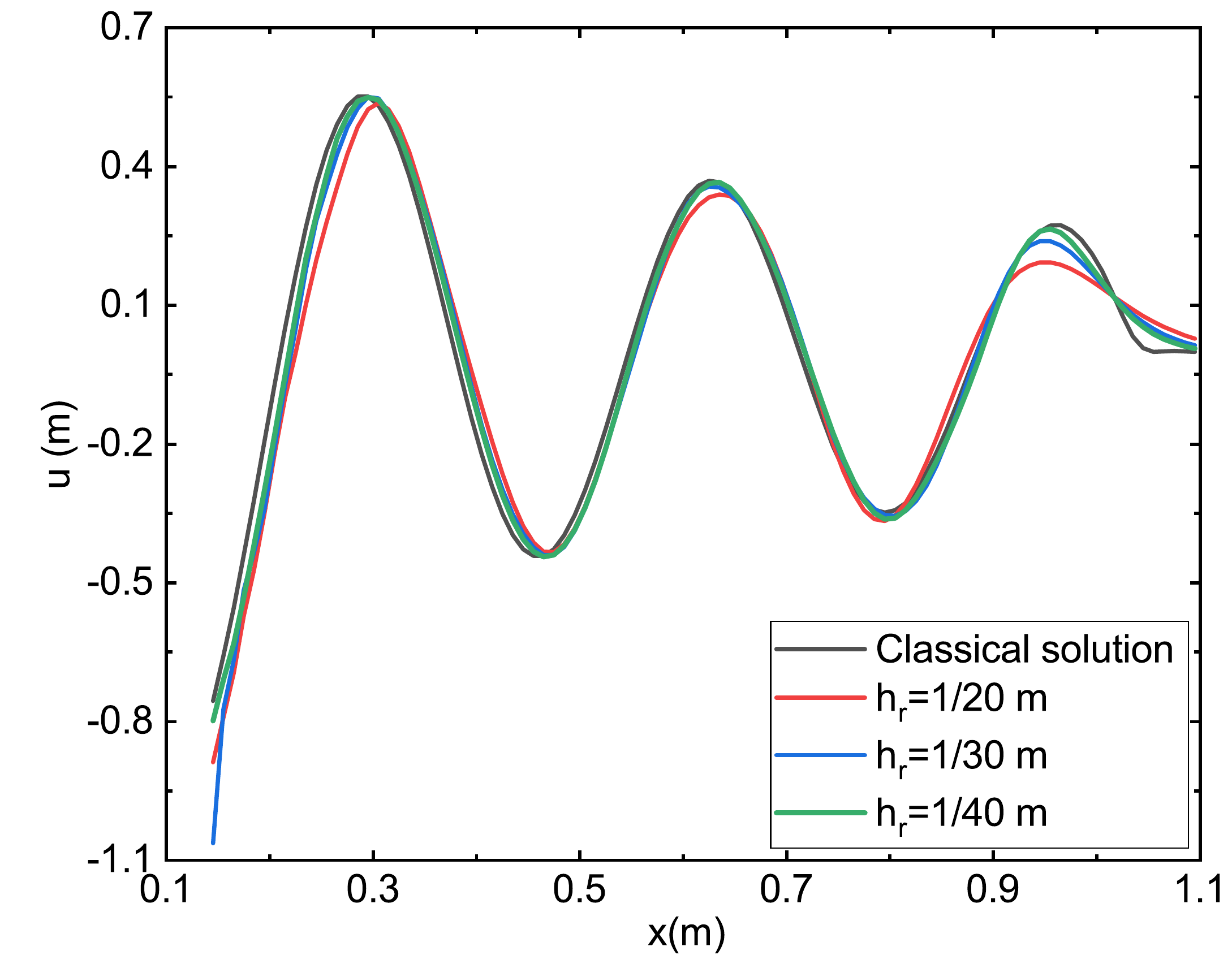}}
	\caption{The radial displacement with the constant value kernel: $t=0.085 {\rm s}$(a); $t=0.485 {\rm s}$(b); $t=0.885 {\rm s}$(c).}
	\label{fig25}	
\end{figure}

\begin{figure}[!htb]
	\centering
	\subfigure[]{
		\label{fig26a}
		\includegraphics[scale=0.16]{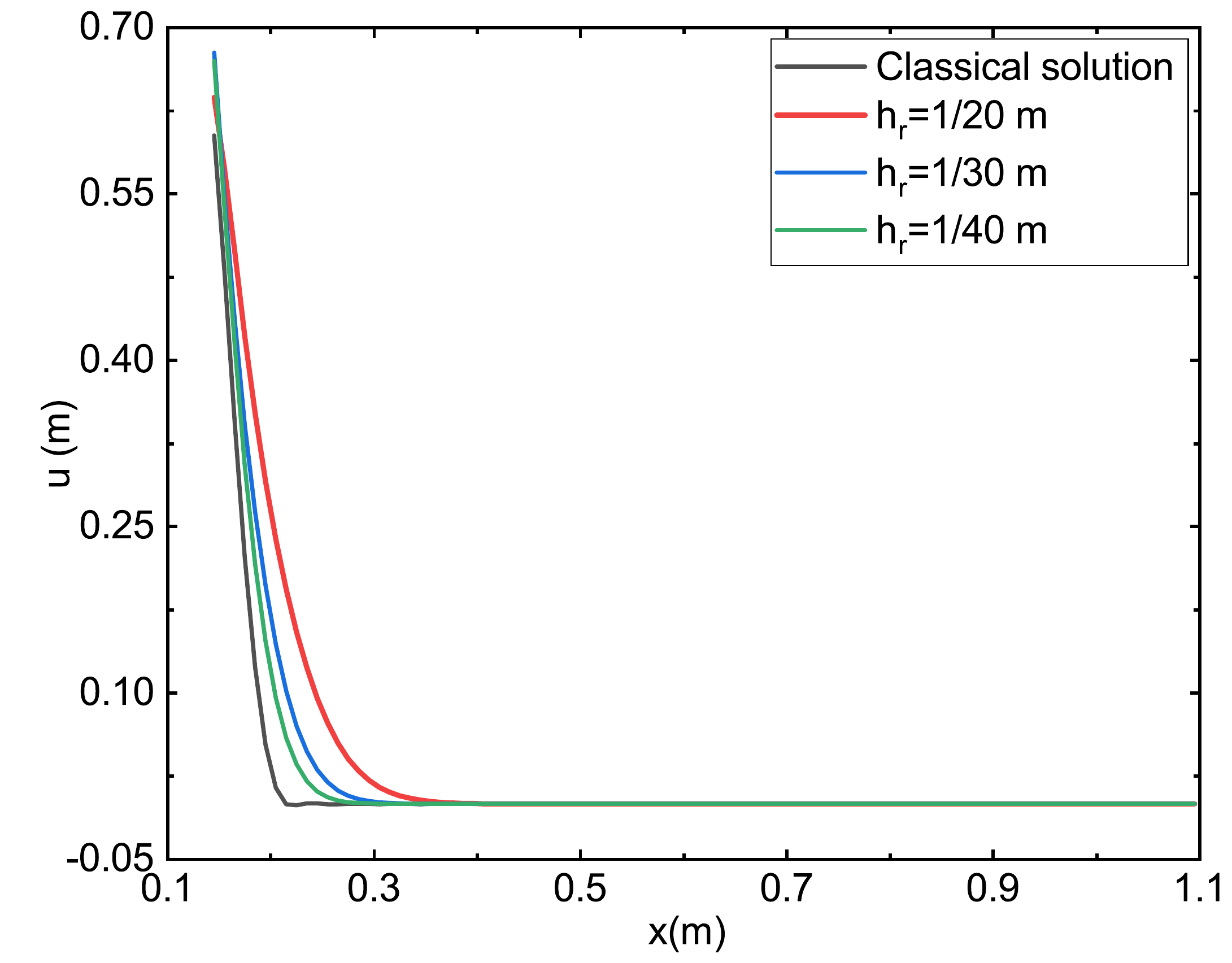}}
	\hspace{0.01in}
	\subfigure[]{
		\label{fig26b}
		\includegraphics[scale=0.16]{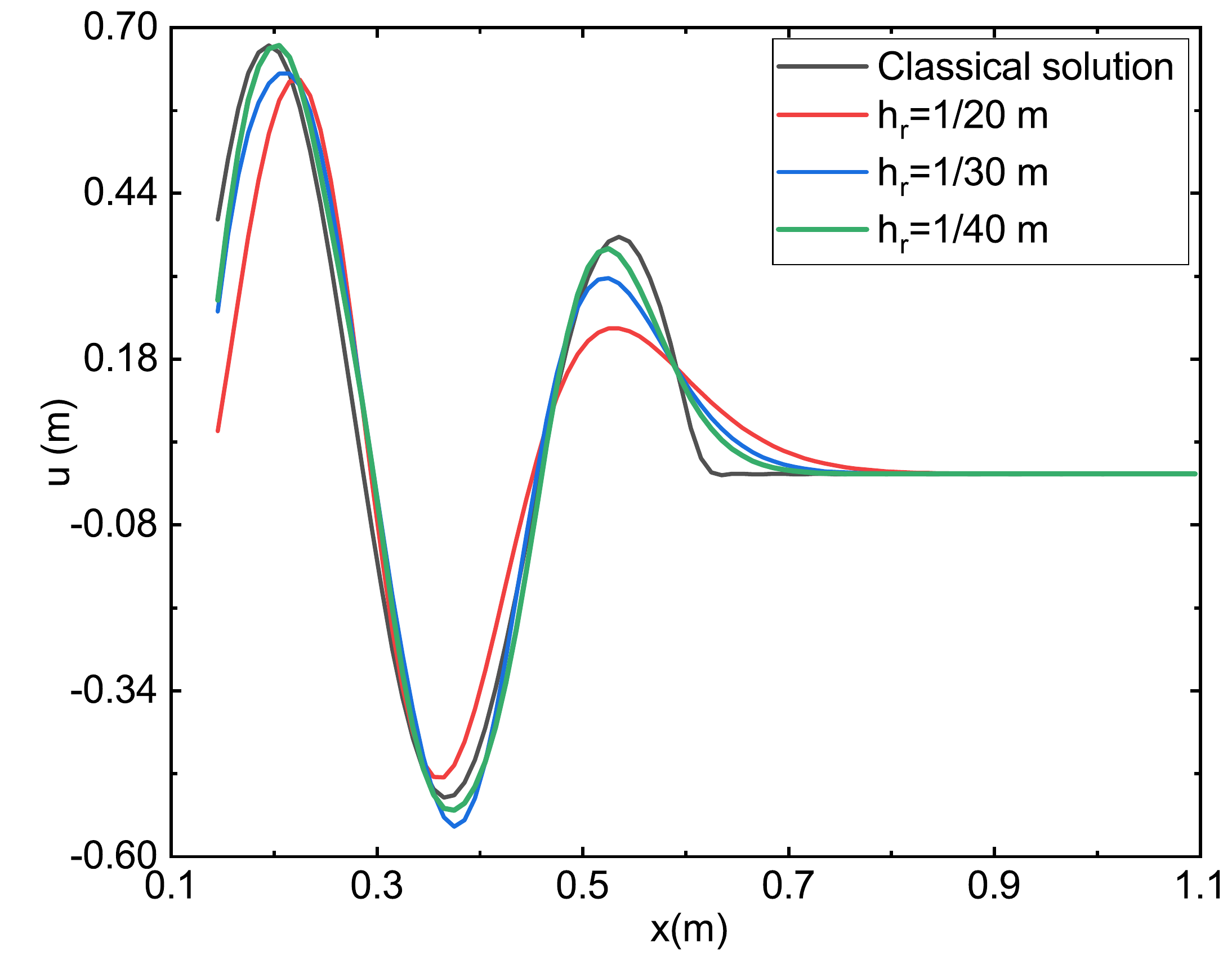}}	
	\hspace{0.01in}
	\subfigure[]{
		\label{fig26c}
		\includegraphics[scale=0.16]{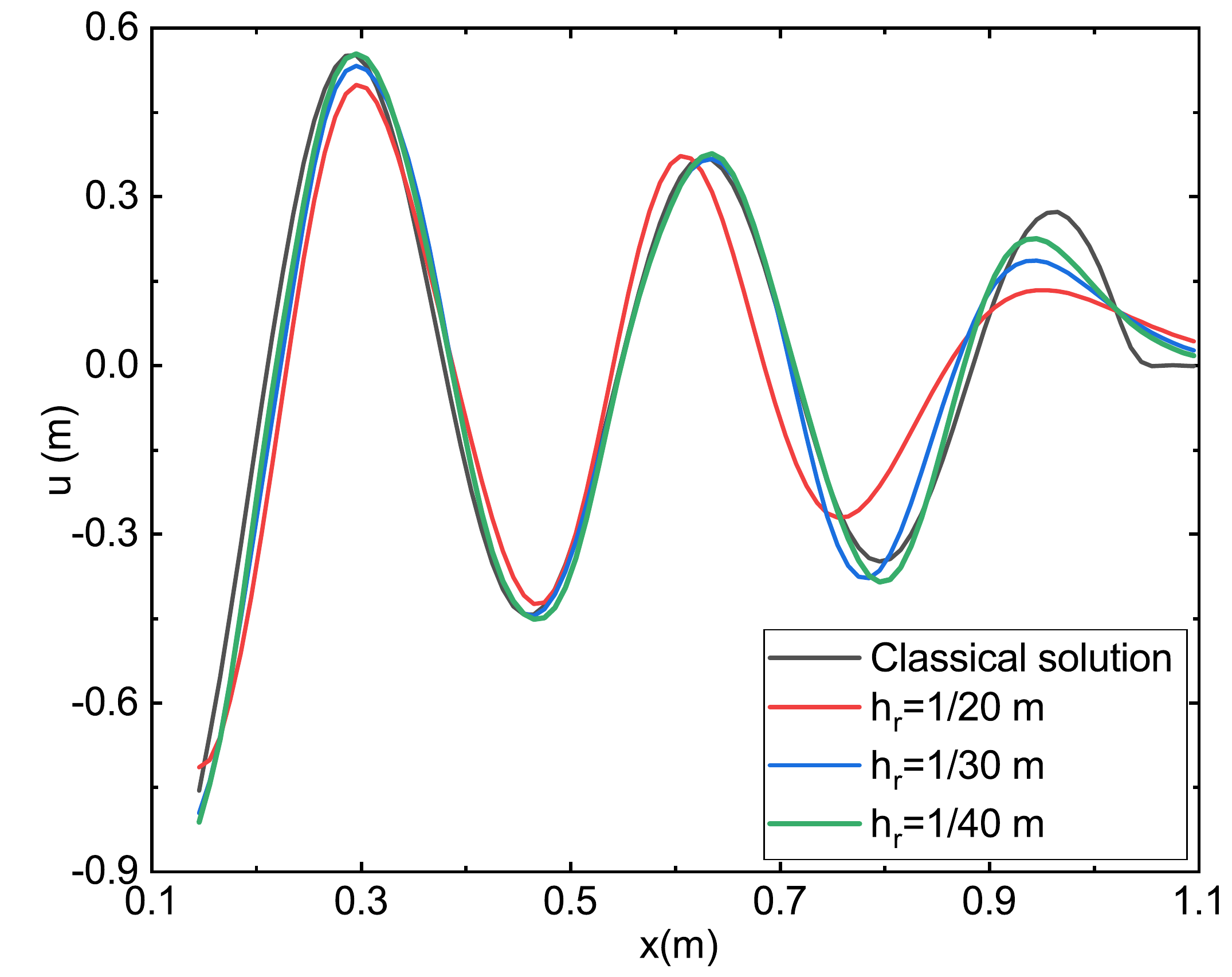}}
	\caption{The radial displacement with the Gauss kernel: $t=0.085 {\rm s}$(a); $t=0.485 {\rm s}$(b); $t=0.885 {\rm s}$(c).}
	\label{fig26}	
\end{figure}

In Table \ref{tab7}, $\overline{{\rm u}}_x\left(\theta,t\right)$ and $\overline{{\rm u}}_y\left(\theta,t\right)$ are harmonic displacements imposed on the boundary of the circular void. Under this dynamic loading condition, we calculate the radial displacement of the material at three different characteristic lengths of $1/20 {\rm m}$, $1/30 {\rm m}$, $1/40 {\rm m}$.  The variations of the radial displacement along the $x$-axis are shown in Figure \ref{fig25} and Figure \ref{fig26}.

\begin{figure}[!htb]
	\centerline{\includegraphics[scale=0.3]{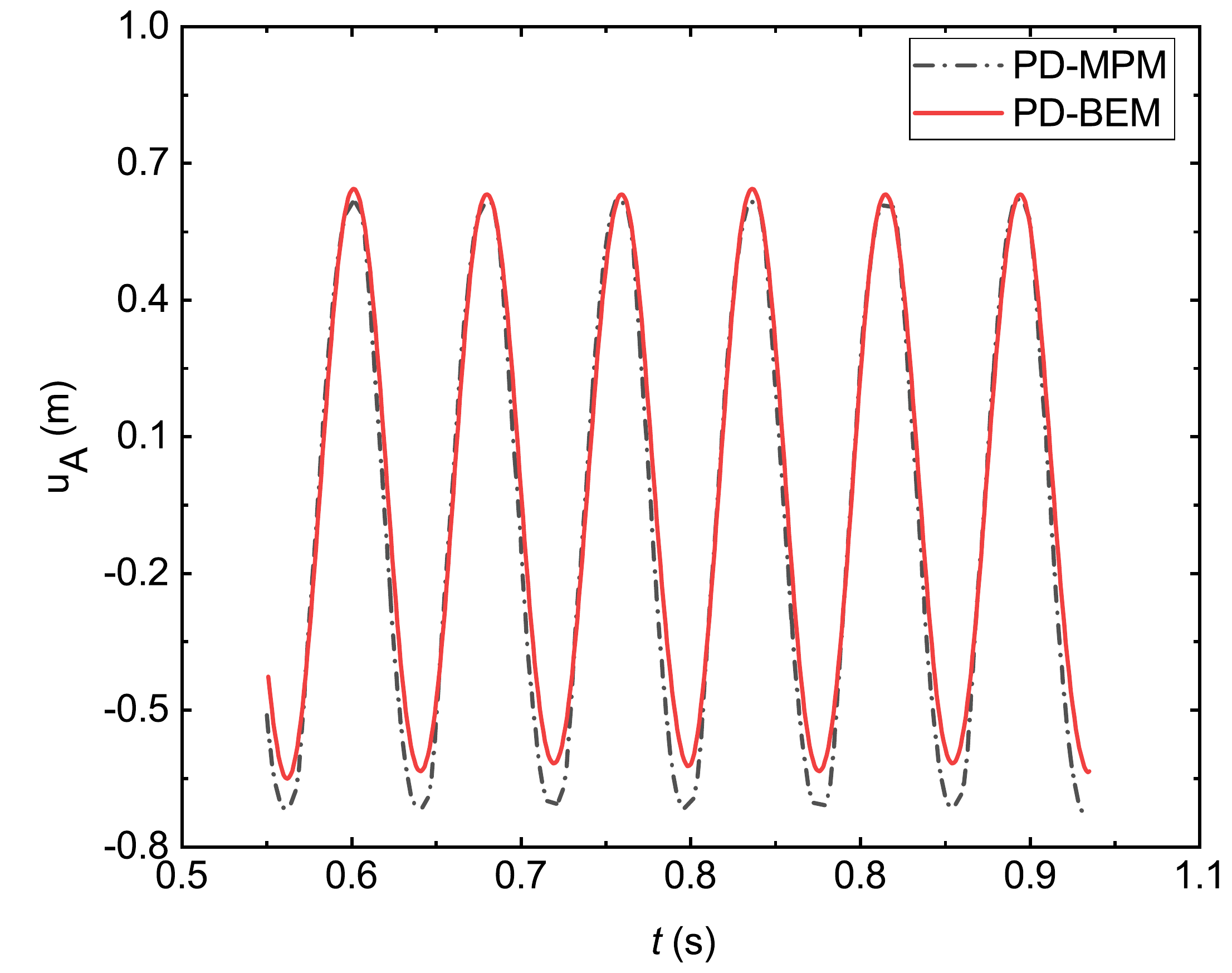}}	
	\caption{The displacement ${\rm u_A}$ in point $A$ versus the time $t$.\label{fig29}}	
\end{figure}

Generally, the larger the computed area, the more efficient the PD-BEM compared to the PD-MPM, especially, for the problem of an infinite domain. For example, we calculate the displacement at point $A$ in Figure~\ref{sc5} whose coordinates are $x=0.2,y=0.0$ from $t=0.5 {\rm s}$ to $t=1 {\rm s}$. We use the constant kernel function with the characteristic length ${\rm h_r}=1/100 {\rm m}$ here. The result is shown in Figure \ref{fig29}.
For the PD-MPM, we simulate the infinite domain by using a square region  10 times the diameter of the hole with $1000\times1000$ particles. The computation takes $46906$ seconds. For PD-BEM, we only need $36$ elements, which takes $5332$ seconds. The ratio of the computational times is about $11.4\%$.

\subsection{Static crack propagation in a double-notched specimen}\label{atl46}

\begin{figure}[!htb]
	\centerline{\includegraphics[scale=0.3]{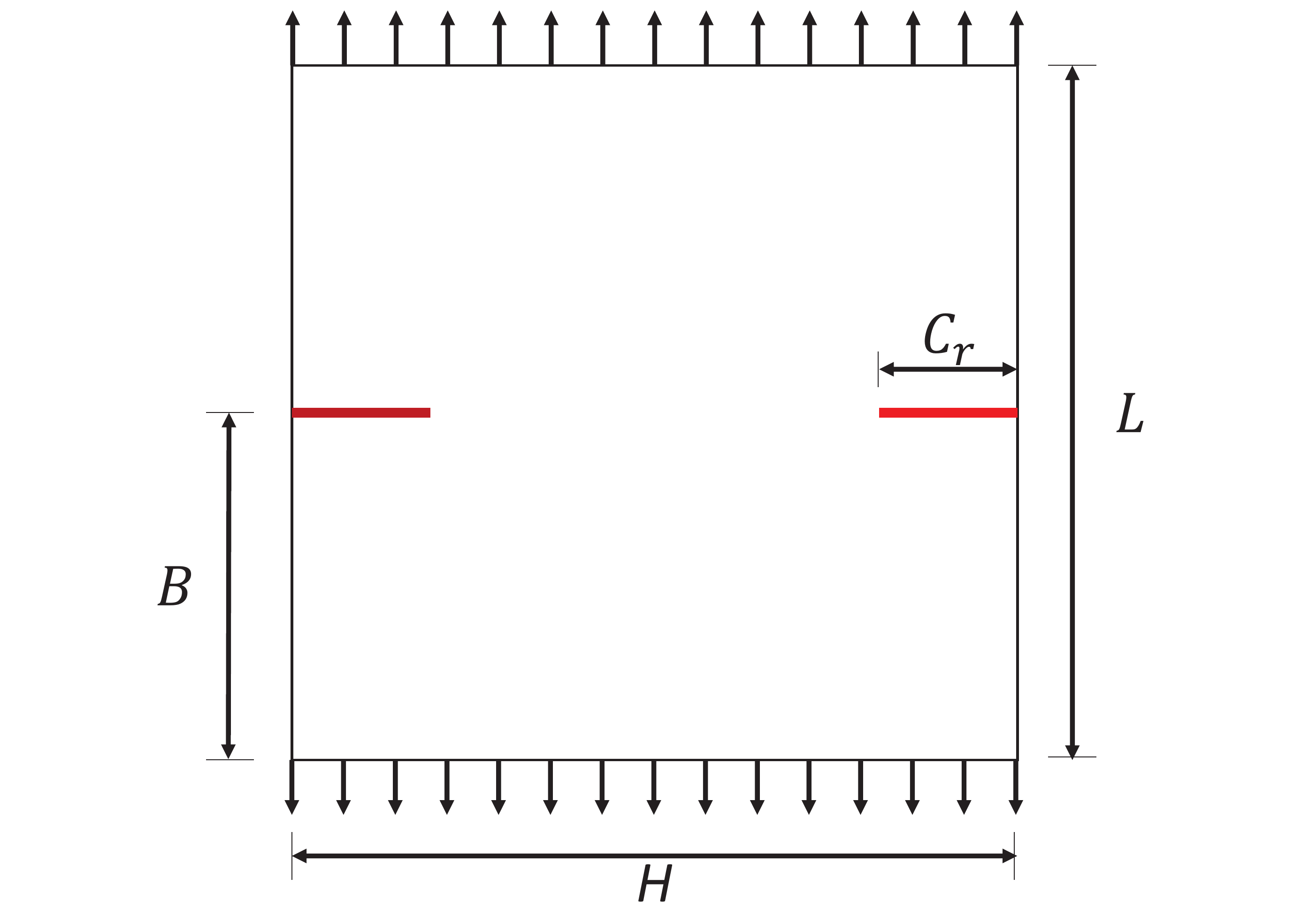}}	
	\caption{A two dimensional double-notched specimen.\label{sc7}}
\end{figure}

In this section and the next, we will apply the PD-BEM to problems of crack propagations and demonstrate its applicability and efficiency.
We only consider the constant kernel function for these examples. First, we calculate the crack propagation in a double-notched specimen under unidirectional tension shown in Figure~\ref{sc7}. For this problem, we compare the computational time and the load that corresponds to a crack propagating length of $1/20$ of the initial length. The used parameters are given in Table \ref{tab9}, where $E$ is Young's modulus; $\nu$ is Poisson's ratio; and $S_c$ is the critical stretch ratio. The meanings of other geometric parameters are shown in Figure~\ref{sc7}.
\begin{center}
	\begin{table}[!htb]
		\centering
		\caption{The material and geometric parameters for two dimensional double-notched specimen.\label{tab9}}
		\begin{tabular}{ccccccc}
			\hline
			$E \left({\rm Pa}\right)$&$\nu$&$C_{r} \left({\rm m}\right)$&$S_c$&$L \left({\rm m}\right)$&$B \left({\rm m}\right)$&$H \left({\rm m}\right)$\\
			\hline
			$1.0$&$1/3$&$0.1$&$0.1$&$1.0$&$0.5$&$1.0$\\
			\hline
		\end{tabular}
	\end{table}
\end{center}

\begin{figure}[!htb]
	\centering
	\subfigure[]{
		\label{fig33a}
		\includegraphics[scale=0.20]{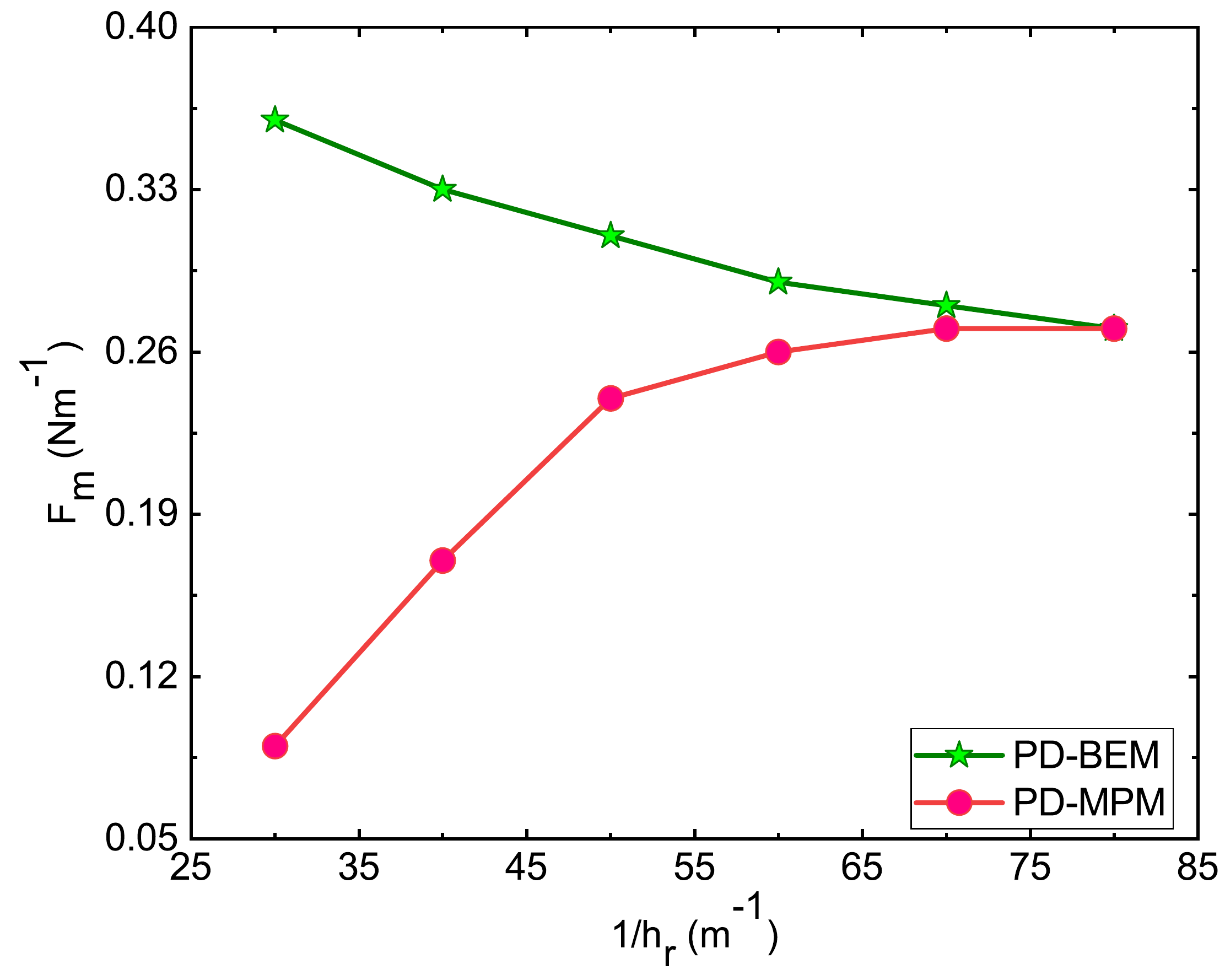}}
	\hspace{0.01in}
	\subfigure[]{
		\label{fig33b}
		\includegraphics[scale=0.20]{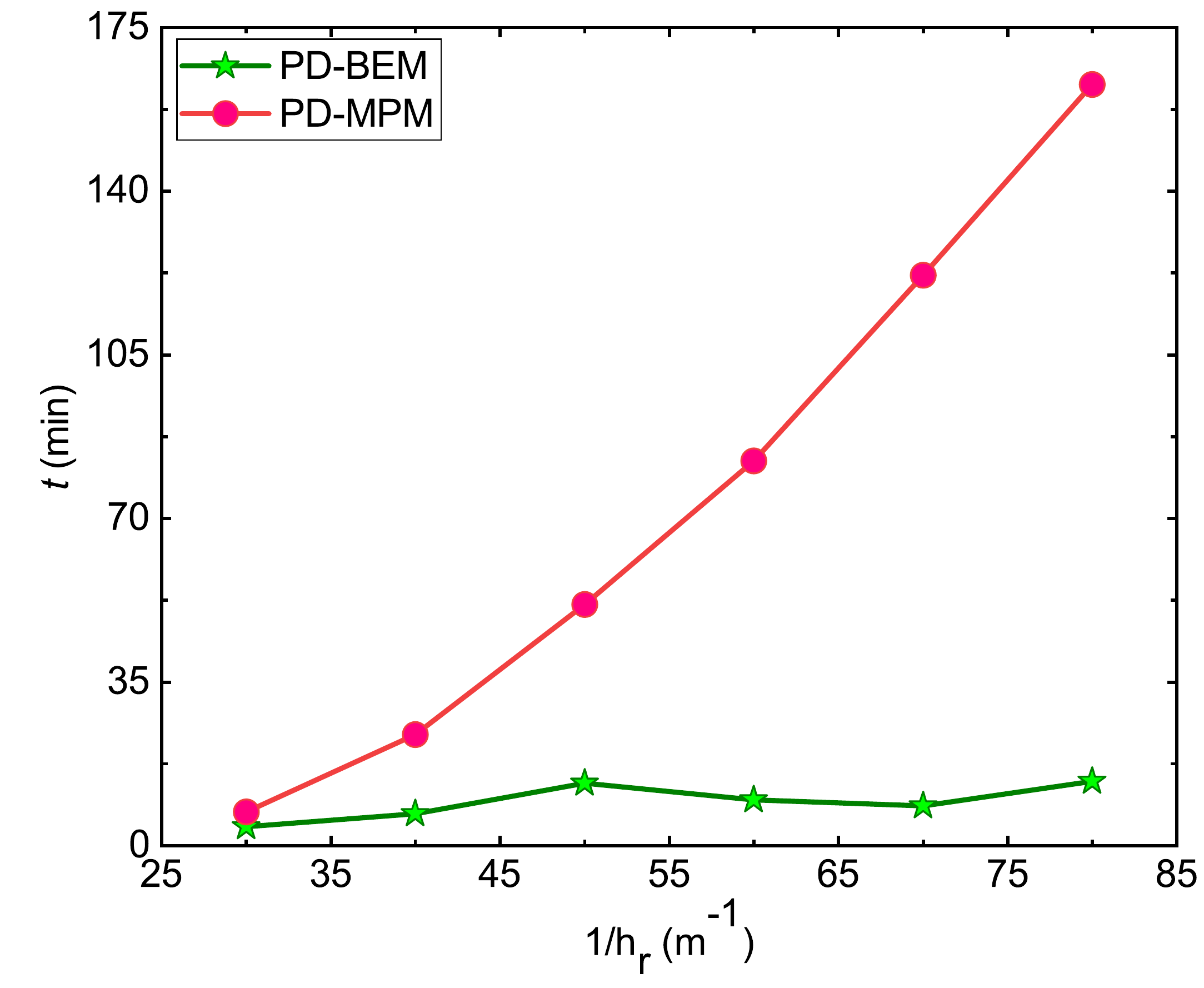}}
	\caption{The load $F_{m}$ (a) and the computational time $t$ (b) versus the reciprocal of characteristic length $1/{\rm h_r}$.}
	\label{fig33}	
\end{figure}

For the PD-BEM, we still use symmetry to simplify computation.  The critical load ($F_{m}$) corresponding to the crack propagating length of $1/20$ of the initial length is calculated by the PD-BEM and the PD-MPM. The critical load and computational time are shown in Figure \ref{fig33}.

The load predicted by the PD-BEM is higher than that by the PD-MPM in Figure \ref{fig33a}. The results of the two methods coincide with each other when the characteristic length tends to zero. There are two reasons for the lower critical load predicted by the PD-MPM. One is that the PD-BEM can give a more accurate result because of the use of the analytical solution; the other is that the surface effect~\cite{tb73} leads to the stiffness reduction near the crack surface, which results in a larger displacement and a spurius bond elongation. The larger the characteristic length, the more serious the surface effect~\cite{tb73}. In addition, the computational time of the PD-MPM rapidly increases with the decrease of  the characteristic length; however, the time of the PD-BEM does not vary much in Figure \ref{fig33b}. This demonstrates the benefits of the PD-BEM that  owns the asymptotic compatibility~\cite{tb80}.

\subsection{Fracture of pre-cracked Brazilian disk}\label{atl47}

\begin{figure}[!htb]
	\centerline{\includegraphics[scale=0.4]{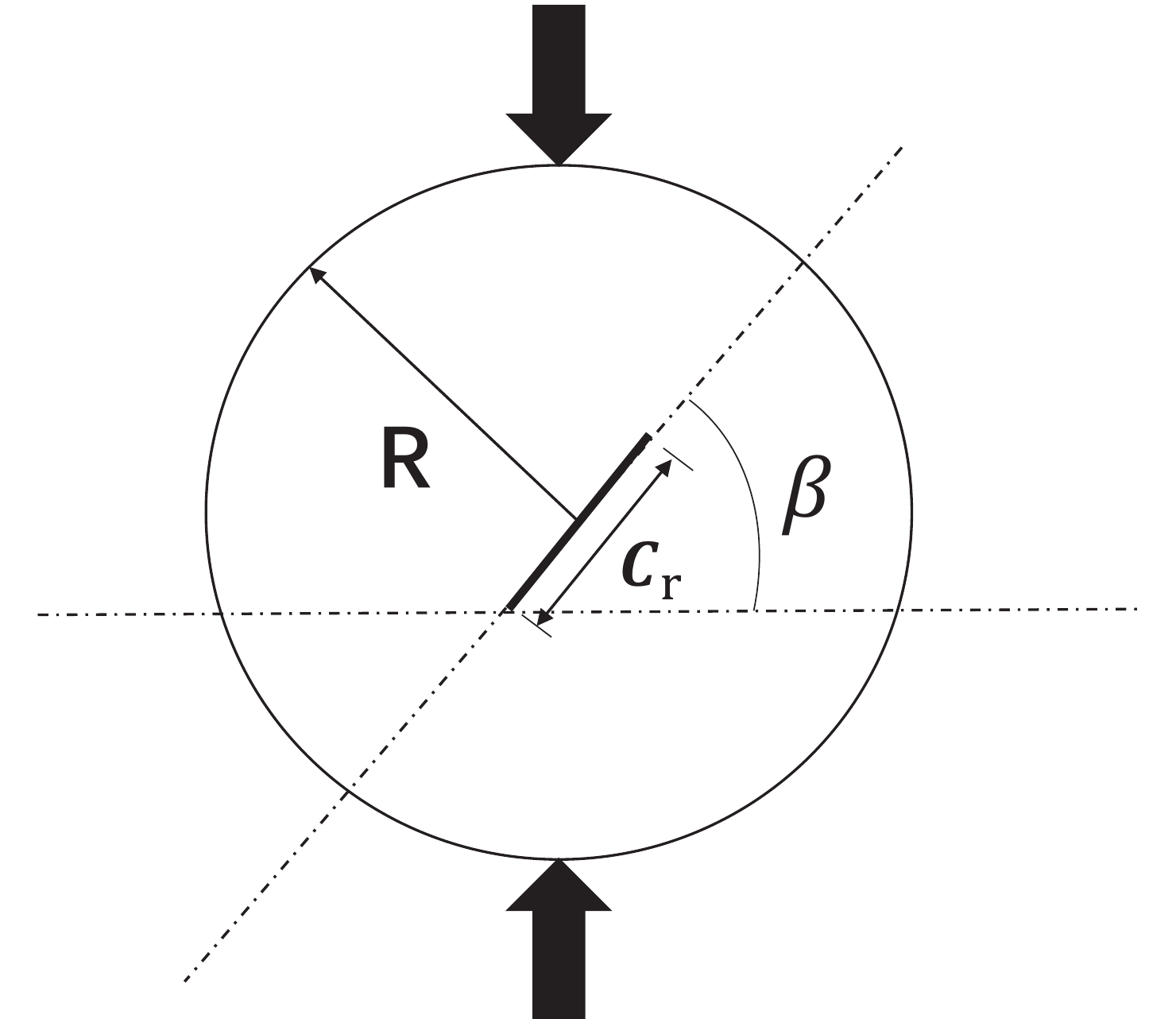}}	
	\caption{Pre-cracked Brazilian disk subjected to quasi-static displacement loading.\label{sc6}}
\end{figure}

\begin{figure}[!htb]
	\centering
	\subfigure[]{
		\label{fig34a}
		\includegraphics[height=4.5cm,width=5.5cm]{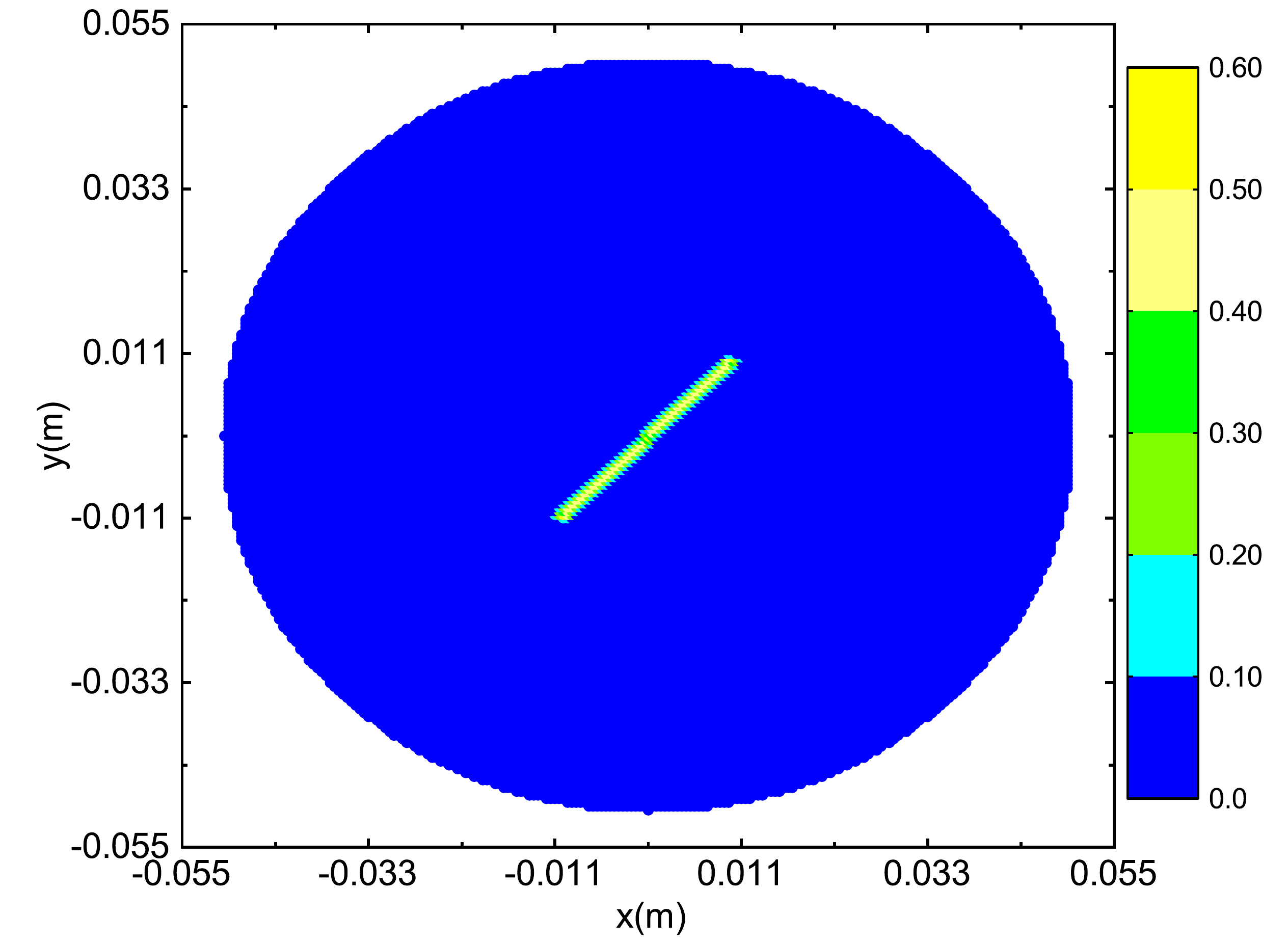}}
	\hspace{0.01in}
	\subfigure[]{
		\label{fig34b}
		\includegraphics[height=4.5cm,width=5.5cm]{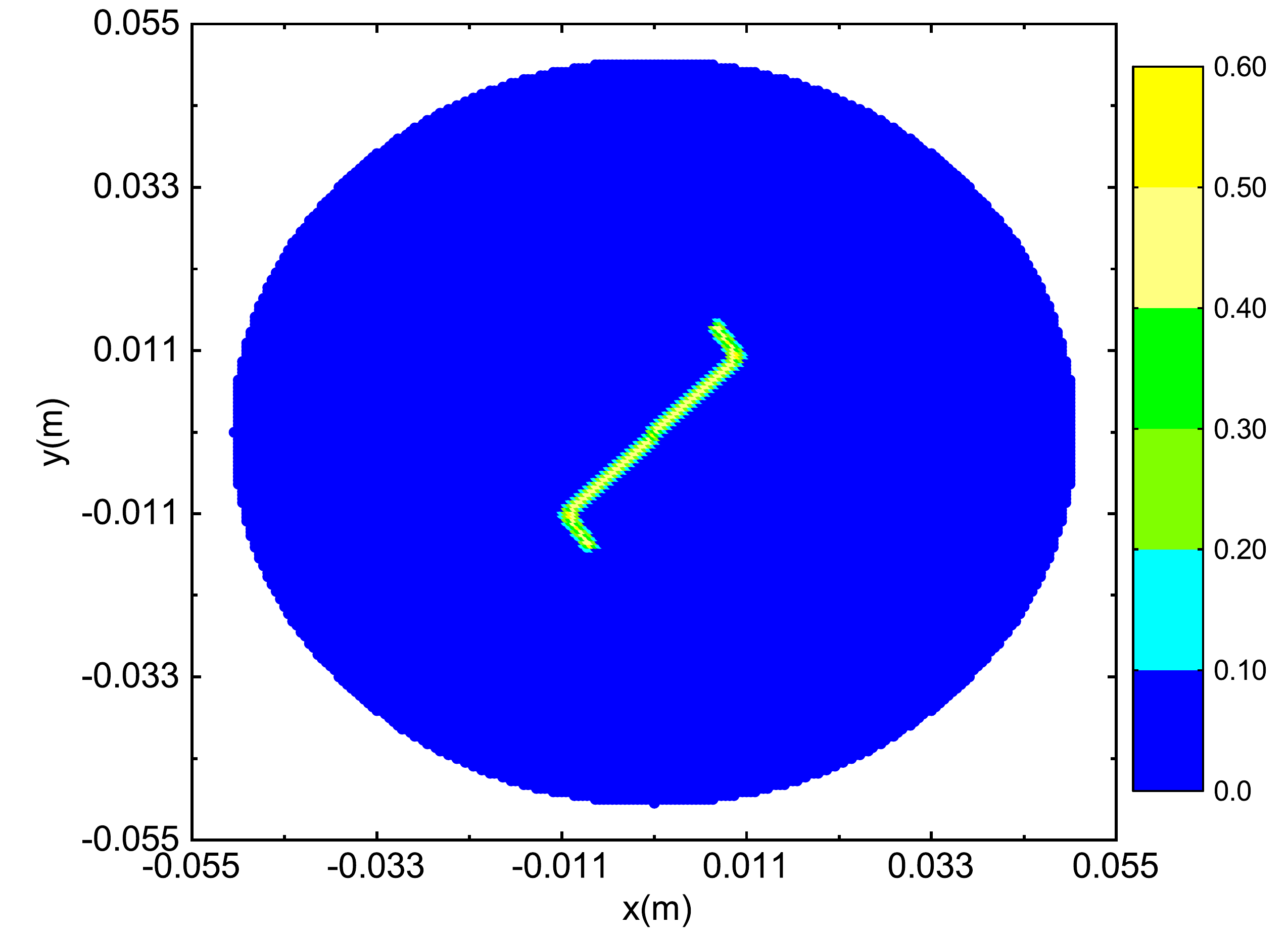}}	
	\vfill
	\subfigure[]{
		\label{fig34c}
		\includegraphics[height=4.5cm,width=5.5cm]{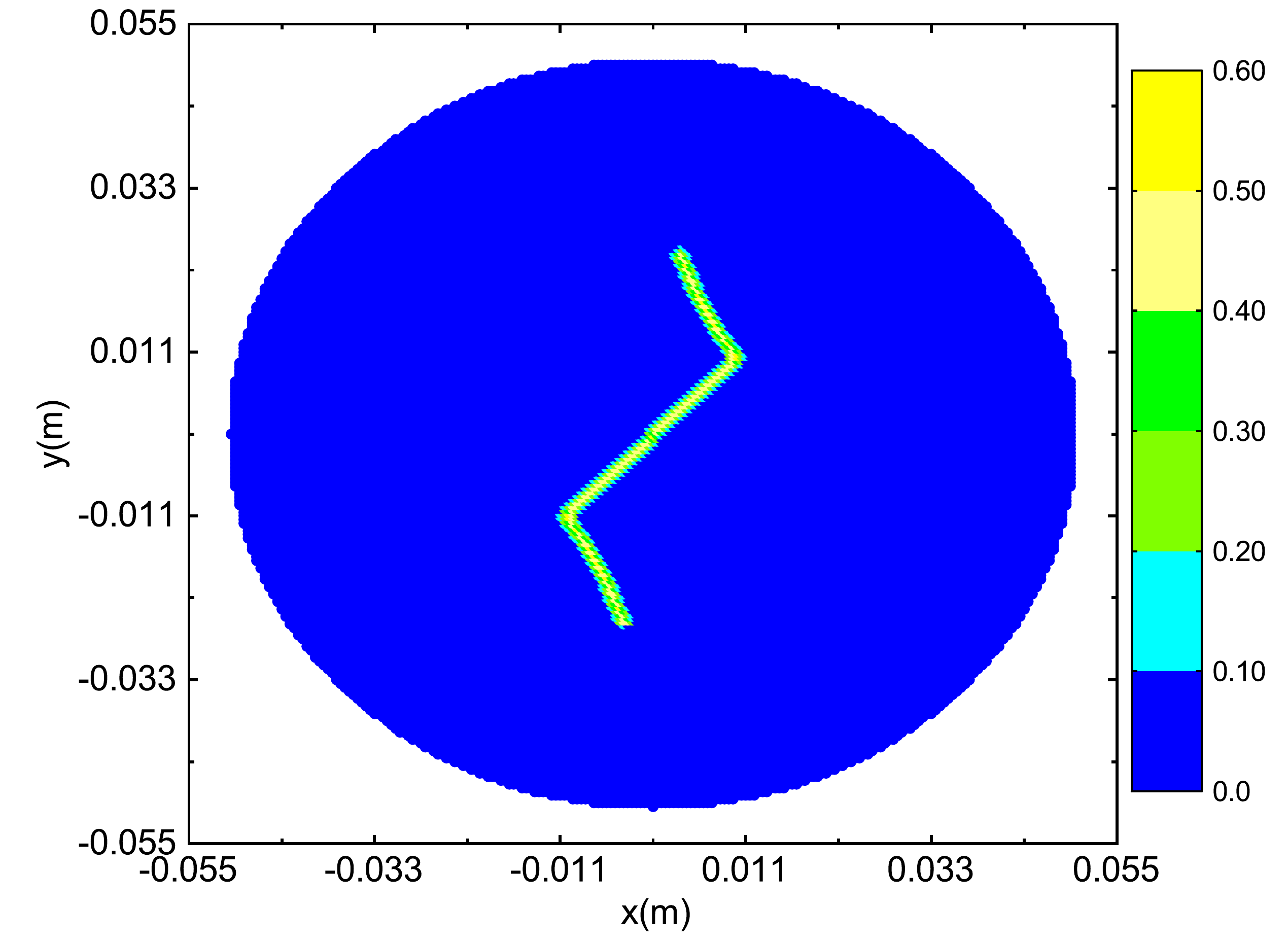}}
	\hspace{0.01in}
	\subfigure[]{
		\label{fig34d}
		\includegraphics[height=4.5cm,width=5.5cm]{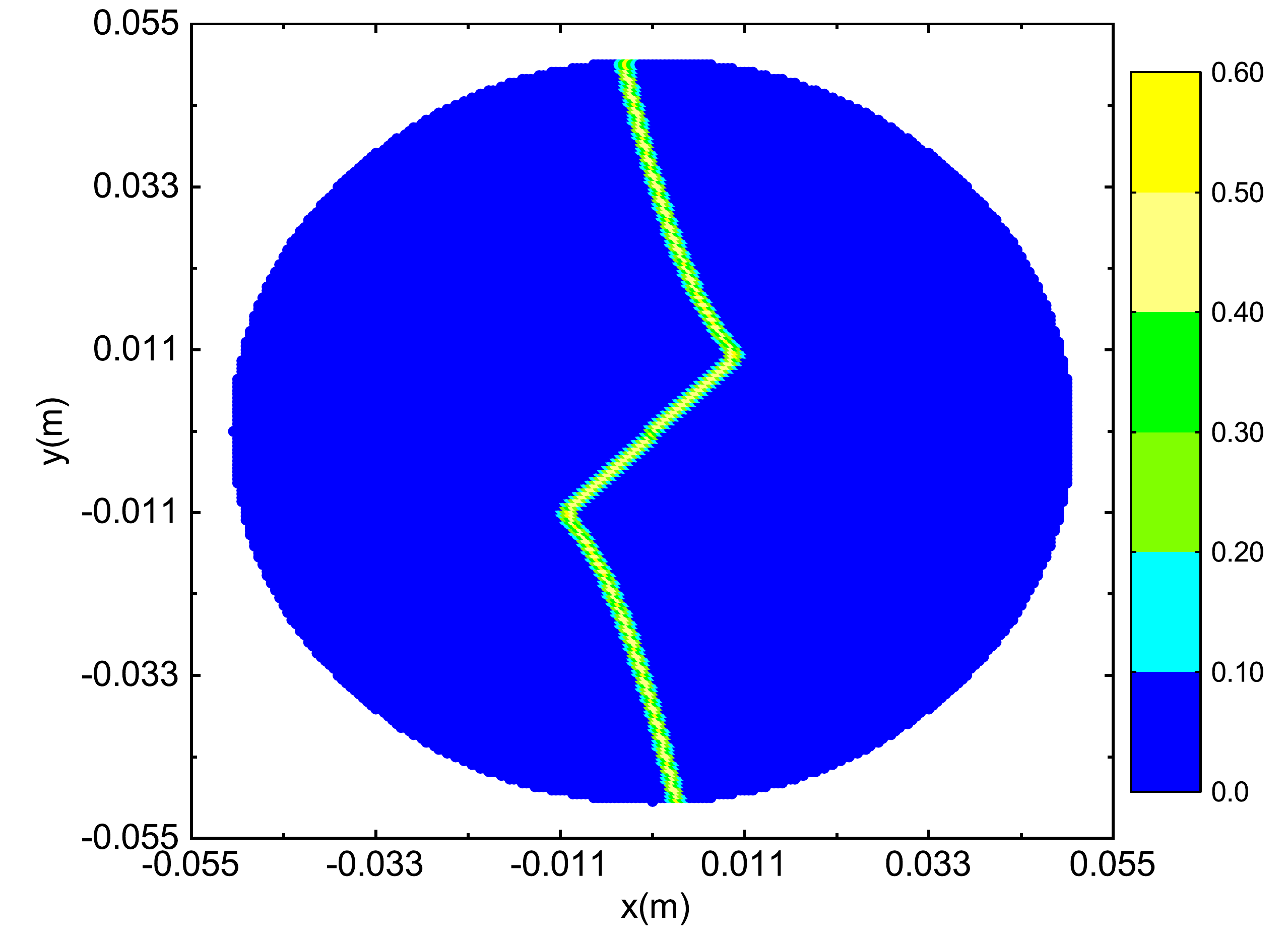}}	
	\caption{The crack propagation path at different loading displacements: $0.19 {\rm mm}$(a); $0.24 {\rm mm}$(b); $0.28 {\rm mm}$(c); $0.33 {\rm mm}$(d).}
	\label{fig34}	
\end{figure}

It is known that the PD-MPM~\cite{tb10} has the advantage of dealing with damage and crack propagation without remeshing, whereas the PD-BEM is highly efficient for problems where no new boundary emerges during the loading process, as proved by the above examples. Therefore, we propose a coupling method of the PD-MPM and PD-BEM (PD-MPM-BEM) to deal with fractures of materials. We apply this coupling method to the modelling of the fracture process of a pre-cracked Brazilian disk under quasi-static displacement loading, shown in Figure \ref{sc6}. In order to maximize the advantage of the PD-BEM in efficiency, the coupling region is adaptively selected. The inner circular area of the disk is computed by using the PD-MPM; the diameter of this PD-MPM area is the maximum distance of all the damaged material points to the center of the disk plus twice the characteristic length. The outer area is computed by using the PD-BEM. Thus, the area of the PD-MPM expands with the propagation of the crack. The loading is exerted through applying a downward displacement at the upper loading point (Figure \ref{sc6}). The geometric and material parameters are listed in Table \ref{tab8}, where $S_c$ is the critical stretch ratio which is used as a failure criterion in the PD-MPM region. The propagating stages of the pre-crack at different loading displacements are shown in Figure \ref{fig34}.
\begin{center}
	\begin{table}[!htb]
		\centering
		\caption{Geometric and material parameters of the pre-cracked Brazilian disk.\label{tab8}}
		\begin{tabular}{cccccc}
			\hline
			$R \left({\rm m}\right)$&$E \left({\rm GPa}\right)$&$\nu$&$C_r \left({\rm m}\right)$&$S_c$&$\beta$\\
			\hline
			$0.05$&$15.0$&$1/3$&$0.03$&$0.00185$&$\pi/4$\\
			\hline
		\end{tabular}
	\end{table}
\end{center}

The path of the crack propagation predicted by the coupling method PD-MPM-BEM is in good agreement with the computational result in the literature~\cite{tb82} and the experimental observation~\cite{tb81}. It takes 357 seconds to calculate the whole fracture process with this coupling method, whereas it takes 1478 seconds to calculate the whole fracture process with the PD-MPM. Therefore, the PD-BEM has the potential of computing fracture problems with high efficiency when coupled with other methods.

\section{Conclusions}\label{atl7}


In this paper, we have developed the \textit{boundary element method of peridynamics} (PD-BEM). This method has several salient features. Firstly, due to boundary discretization, this method, compared with the domain method, is not affected by breaking the symmetry of the horizon near the boundary, and it is easy to handle infinite domains. Secondly, in order to overcome the difficulty of solving dynamic problems in the time domain, we introduce the Laplace transformation with respect to time. Solving dynamic problems in the Laplace domain eliminates time history dependence and time  accumulation error. Thirdly, due to dimension reduction, it is highly efficient in terms of computational time, compared with the methods of discretizing the whole domain. Theoretically, like the classical boundary element method, it is less efficient for problems with new boundaries emerging during the process of loading because of the requirement for re-meshing the new boundaries. Notwithstanding, the PD-BEM can be coupled with other numerical methods to compute fracture process with enhanced efficiency, which may deserve further studies. Finally, we give the reciprocal theorem for the state-based peridynamic theory in Appendix \ref{ap8}. With the infinite Green function for the state-based PD~\cite{tb88}, the corresponding PD-BEM can also be established. The nonlocal boundary element method can also be developed for other physical problems such as diffusion when the the  Green function is available.

\section*{Nomenclature}

\begin{table}[!htb]
	\begin{tabular}{cc}
		PD&Peridynamics\\
		PD-BEM&Boundary element method of peridynamics\\
		PD-MPM&Meshless particle method of peridynamics\\
		$\mathcal{H}_{\mathbf{x}}$&Family of point $\mathbf{x}$\\
		$\mathcal{N}$&Nonlocal interaction operator\\
		$\mathcal{D}$&Nonlocal divergence operator\\
		$\mathcal{D}^{*}$&Conjugate operator of nonlocal divergence operator\\
		$\Omega$&Solution domain\\
		$\Omega_\tau$&Volume-constrained boundary\\
		$\Omega_{\tau_d}$&Displacement volume-constrained boundary\\
		$\Omega_{\tau_n}$&Force volume-constrained boundary\\
		$\partial\Omega$&Local boundary\\
		$\partial \Omega_d$&Local displacement boundary\\
		$\partial \Omega_n$&Local force boundary\\
		$\boldsymbol{\Theta}$&Fourth-order elasticity tensor\\
		$\overline{\mathbf{T}}$&Nonlocal surface traction operator\\
		$\underline{\omega} \left \langle \boldsymbol{\xi} \right \rangle$&Nonlocal weight function\\
		$\delta \left(\mathbf{x} - \mathbf{x}_0\right)$&Dirac function\\
		$\bbint$&Cauchy principal value of the integral\\
		$\mathcal{L}$&Laplace transformation\\
		$\mathcal{F}$&Fourier transformation\\
		$\delta_{ij}$&Kronecker delta\\
		${\rm h_r}$&Characteristic length\\
		$Re$&Real part\\
	\end{tabular}
\end{table}

\section*{Acknowledgments}

This work is supported by the National Key R\&D Program of China (2020YFA0710500), and the Natural Science Foundation of China (Grant No.~11890681). We acknowledge the High-performance Computing Platform of Peking University for providing computational resources.

\appendix

\section{Nonlocal dynamic boundary integral equation for time domain}\label{ap4}

We follow the derivation for the static problem. Two independent motion states are given as follows:
\begin{sequation}\label{dpf1}
Problem \ 5 :  \qquad
\begin{cases}
\mathcal{D} \left(\boldsymbol{\Theta} : \mathcal{D}^{*} \left(\mathbf{v}_5\right) \right) \left(\mathbf{x},t\right) + \mathbf{c} \left(\mathbf{x},t\right) = \rho \ddot{\mathbf{v}}_5 \left(\mathbf{x},t\right) \qquad & \text{in} \ \mathbf{x} \in \Omega \\
\mathbf{v}_5 \left(\mathbf{x},t\right) = \mathbf{g}_{d_5} \left(\mathbf{x},t\right) \qquad & \text{in} \ \mathbf{x} \in \Omega_{\tau_d} \\
\mathcal{N} \left(\boldsymbol{\Theta} : \mathcal{D}^{*} \left(\mathbf{v}_5\right) \right) \left(\mathbf{x},t\right) = \mathbf{g}_{n_5} \left(\mathbf{x},t\right) \qquad & \text{in} \ \mathbf{x} \in \Omega_{\tau_n} \\
\mathbf{v}_5 \left(\mathbf{x},t_0\right) = \tilde{\mathbf{v}}_5 \left(\mathbf{x}\right) \qquad \dot{\mathbf{v}}_5 \left(\mathbf{x},t_0\right) = \dot{\tilde{\mathbf{v}}}_5 \left(\mathbf{x}\right)
\end{cases}
\end{sequation}
\begin{sequation}\label{dpf2}
Problem \ 6 :  \qquad
\begin{cases}
\mathcal{D} \left(\boldsymbol{\Theta} : \mathcal{D}^{*} \left(\mathbf{v}_6\right) \right) \left(\mathbf{x},t\right) + \mathbf{d} \left(\mathbf{x},t\right) = \rho \ddot{\mathbf{v}}_6 \left(\mathbf{x},t\right) \qquad & \text{in} \ \mathbf{x} \in \Omega \\
\mathbf{v}_6 \left(\mathbf{x},t\right) = \mathbf{g}_{d_6} \left(\mathbf{x},t\right) \qquad & \text{in} \ \mathbf{x} \in \Omega_{\tau_d} \\
\mathcal{N} \left(\boldsymbol{\Theta} : \mathcal{D}^{*} \left(\mathbf{v}_6\right) \right) \left(\mathbf{x},t\right) = \mathbf{g}_{n_6} \left(\mathbf{x},t\right) \qquad & \text{in} \ \mathbf{x} \in \Omega_{\tau_n} \\
\mathbf{v}_6 \left(\mathbf{x},t_0\right) = \tilde{\mathbf{v}}_6 \left(\mathbf{x}\right) \qquad \dot{\mathbf{v}}_6 \left(\mathbf{x},t_0\right) = \dot{\tilde{\mathbf{v}}}_6 \left(\mathbf{x}\right)
\end{cases}
\end{sequation}
where $\mathbf{v}_5$ and $\mathbf{v}_6$ are the displacements of two motion states; $\mathbf{c}$ and $\mathbf{d}$ are body force densities;  $\tilde{\mathbf{v}}$ and $\dot{\tilde{\mathbf{v}}}$ are initial conditions. $t_0$ is the initial time. For these two states, we will derive the dynamic reciprocal theorem. The motion states $\mathbf{v}_5 \left(\mathbf{x},\tau\right)$ and $\mathbf{v}_6 \left(\mathbf{x},t - \tau + t_0\right)$ satisfy the following two sets of equations:
\begin{sequation}\label{dpf6}
\begin{cases}
\mathcal{D} \left(\boldsymbol{\Theta} : \mathcal{D}^{*} \left(\mathbf{v}_5\right) \right) \left(\mathbf{x},\tau\right) + \mathbf{c} \left(\mathbf{x},\tau\right) - \rho \ddot{\mathbf{v}}_5 \left(\mathbf{x},\tau\right) = \mathbf{0} \qquad & \text{in} \ \mathbf{x} \in \Omega \\
\mathbf{v}_5 \left(\mathbf{x},\tau\right) = \mathbf{g}_{d_5} \left(\mathbf{x},\tau\right) \qquad & \text{in} \ \mathbf{x} \in \Omega_{\tau_d} \\
\mathcal{N} \left(\boldsymbol{\Theta} : \mathcal{D}^{*} \left(\mathbf{v}_5\right) \right) \left(\mathbf{x},\tau\right) = \mathbf{g}_{n_5} \left(\mathbf{x},\tau\right) \qquad & \text{in} \ \mathbf{x} \in \Omega_{\tau_n}
\end{cases}
\qquad \qquad \,
\end{sequation}
\begin{sequation}\label{dpf7}
\begin{cases}
\mathcal{D} \left(\boldsymbol{\Theta} : \mathcal{D}^{*} \left(\mathbf{v}_6\right) \right) \left(\mathbf{x},t - \tau + t_0\right) + \mathbf{d} \left(\mathbf{x},t - \tau + t_0\right) - \rho \ddot{\mathbf{v}}_6 \left(\mathbf{x},t - \tau + t_0\right) = \mathbf{0} \qquad & \text{in} \ \mathbf{x} \in \Omega \\
\mathbf{v}_6 \left(\mathbf{x},t - \tau + t_0\right) = \mathbf{g}_{d_6} \left(\mathbf{x},t - \tau + t_0\right) \qquad & \text{in} \ \mathbf{x} \in \Omega_{\tau_d} \\
\mathcal{N} \left(\boldsymbol{\Theta} : \mathcal{D}^{*} \left(\mathbf{v}_6\right) \right) \left(\mathbf{x},t - \tau + t_0\right) = \mathbf{g}_{n_6} \left(\mathbf{x},t - \tau + t_0\right) \qquad & \text{in} \ \mathbf{x} \in \Omega_{\tau_n}
\end{cases}
\end{sequation}
If we fix the time and consider the inertia force as a kind of body force, then we can regard (\ref{dpf6}) and (\ref{dpf7}) as the governing equations of two independent static states. Applying the static reciprocal theorem (\ref{fa10}) for (\ref{dpf6}) and (\ref{dpf7}), we get
\begin{sequation}\label{dpf8}
\begin{split}
& \int_{\Omega} \mathbf{v}_5 \left(\mathbf{x},\tau\right) \cdot \mathcal{D} \left(\boldsymbol{\Theta} : \mathcal{D}^{*} \left(\mathbf{v}_6\right) \right) \left(\mathbf{x},t - \tau + t_0\right) - \mathbf{v}_6 \left(\mathbf{x},t - \tau + t_0\right) \cdot \mathcal{D} \left(\boldsymbol{\Theta} : \mathcal{D}^{*} \left(\mathbf{v}_5\right) \right) \left(\mathbf{x},\tau\right) {\rm d}V_{\mathbf{x}} \\
=& \int_{\Omega_\tau} \mathbf{v}_5 \left(\mathbf{x},\tau\right) \cdot \mathcal{N} \left(\boldsymbol{\Theta} : \mathcal{D}^{*} \left(\mathbf{v}_6\right) \right) \left(\mathbf{x},t - \tau + t_0\right) - \mathbf{v}_6 \left(\mathbf{x},t - \tau + t_0\right) \cdot \mathcal{N} \left(\boldsymbol{\Theta} : \mathcal{D}^{*} \left(\mathbf{v}_5\right) \right) \left(\mathbf{x},\tau\right) {\rm d}V_{\mathbf{x}}
\end{split}
\end{sequation}
Integrating (\ref{dpf8}) with respect to $\tau$ in the interval $\left(t_0,t\right)$  yields
\begin{sequation}\label{dpf5}
\begin{split}
& \int_{t_0}^{t} \int_{\Omega} \mathbf{v}_5 \left(\mathbf{x},\tau\right) \cdot \mathcal{D} \left(\boldsymbol{\Theta} : \mathcal{D}^{*} \left(\mathbf{v}_6\right) \right) \left(\mathbf{x},t - \tau + t_0\right) - \mathbf{v}_6 \left(\mathbf{x},t - \tau + t_0\right) \cdot \mathcal{D} \left(\boldsymbol{\Theta} : \mathcal{D}^{*} \left(\mathbf{v}_5\right) \right) \left(\mathbf{x},\tau\right) {\rm d}\tau {\rm d}V_{\mathbf{x}} \\
=& \int_{t_0}^{t} \int_{\Omega_\tau} \mathbf{v}_5 \left(\mathbf{x},\tau\right) \cdot \mathcal{N} \left(\boldsymbol{\Theta} : \mathcal{D}^{*} \left(\mathbf{v}_6\right) \right) \left(\mathbf{x},t - \tau + t_0\right) - \mathbf{v}_6 \left(\mathbf{x},t - \tau + t_0\right) \cdot \mathcal{N} \left(\boldsymbol{\Theta} : \mathcal{D}^{*} \left(\mathbf{v}_5\right) \right) \left(\mathbf{x},\tau\right) {\rm d}\tau {\rm d}V_{\mathbf{x}}
\end{split}
\end{sequation}
We can rewrite (\ref{dpf5}) by using convolution commutativity and exchanging the integral order as follows:
\begin{sequation}\label{dpf4}
\begin{split}
& \int_{\Omega} \int_{t_0}^{t} \mathbf{v}_5 \left(\mathbf{x},\tau\right) \cdot \mathcal{D} \left(\boldsymbol{\Theta} : \mathcal{D}^{*} \left(\mathbf{v}_6\right)\right) \left(\mathbf{x},t - \tau + t_0\right) - \mathbf{v}_6 \left(\mathbf{x},\tau\right) \cdot \mathcal{D} \left(\boldsymbol{\Theta} : \mathcal{D}^{*} \left(\mathbf{v}_5\right)\right) \left(\mathbf{x},t - \tau + t_0\right) {\rm d}\tau {\rm d}V_{\mathbf{x}} \\
=& \int_{\Omega_\tau} \int_{t_0}^{t} \mathbf{v}_5 \left(\mathbf{x},\tau\right) \cdot \mathcal{N} \left(\boldsymbol{\Theta} : \mathcal{D}^{*} \left(\mathbf{v}_6\right) \right) \left(\mathbf{x},t - \tau + t_0\right) - \mathbf{v}_6 \left(\mathbf{x},\tau\right) \cdot \mathcal{N} \left(\boldsymbol{\Theta} : \mathcal{D}^{*} \left(\mathbf{v}_5\right) \right) \left(\mathbf{x},t - \tau + t_0\right) {\rm d}\tau {\rm d}V_{\mathbf{x}}
\end{split}
\end{sequation}
Let us rewrite (\ref{dpf4}) in the form of convolution as follows:
\begin{sequation}\label{dpf3}
\begin{split}
& \int_{\Omega} \mathbf{v}_5 \left(\mathbf{x},t\right) * \mathcal{D} \left(\boldsymbol{\Theta} : \mathcal{D}^{*} \left(\mathbf{v}_6\right) \right) \left(\mathbf{x},t\right) - \mathbf{v}_6 \left(\mathbf{x},t\right) * \mathcal{D} \left(\boldsymbol{\Theta} : \mathcal{D}^{*} \left(\mathbf{v}_5\right) \right) \left(\mathbf{x},t\right) {\rm d}V_{\mathbf{x}} \\
=& \int_{\Omega_\tau} \mathbf{v}_5 \left(\mathbf{x},t\right) * \mathcal{N} \left(\boldsymbol{\Theta} : \mathcal{D}^{*} \left(\mathbf{v}_6\right) \right) \left(\mathbf{x},t\right) - \mathbf{v}_6 \left(\mathbf{x},t\right) * \mathcal{N} \left(\boldsymbol{\Theta} : \mathcal{D}^{*} \left(\mathbf{v}_5\right) \right) \left(\mathbf{x},t\right) {\rm d}V_{\mathbf{x}}
\end{split}
\end{sequation}
where $*$ is convolution about the time $t$. This is the dynamic reciprocal theorem. Using the governing equations of  $Problem\ \ 5$ and $Problem\ \ 6$, and carrying out the partial integral with respect to the time variable in (\ref{dpf3}), (\ref{dpf3}) can be further rewritten as the following form:
\begin{sequation}\label{dpf9}
\begin{split}
& \int_{\Omega} \tilde{\mathbf{v}}_5 \left(\mathbf{x}\right) \cdot \dot{\mathbf{v}}_6 \left(\mathbf{x},t\right) + \dot{\tilde{\mathbf{v}}}_5 \left(\mathbf{x}\right) \cdot \mathbf{v}_6 \left(\mathbf{x},t\right) - \mathbf{v}_5 \left(\mathbf{x},t\right) * \mathbf{d} \left(\mathbf{x},t\right) {\rm d}V_{\mathbf{x}} \\
=& \int_{\Omega_\tau} \mathbf{v}_5 \left(\mathbf{x},t\right) * \mathcal{N} \left(\boldsymbol{\Theta} : \mathcal{D}^{*} \left(\mathbf{v}_6\right) \right) \left(\mathbf{x},t\right) - \mathbf{v}_6 \left(\mathbf{x},t\right) * \mathcal{N} \left(\boldsymbol{\Theta} : \mathcal{D}^{*} \left(\mathbf{v}_5\right) \right) \left(\mathbf{x},t\right) {\rm d}V_{\mathbf{x}} \\
& + \int_{\Omega} \tilde{\mathbf{v}}_6 \left(\mathbf{x}\right) \cdot \dot{\mathbf{v}}_5 \left(\mathbf{x},t\right) + \dot{\tilde{\mathbf{v}}}_6 \left(\mathbf{x}\right) \cdot \mathbf{v}_5 \left(\mathbf{x},t\right) - \mathbf{v}_6 \left(\mathbf{x},t\right) * \mathbf{c} \left(\mathbf{x},t\right) {\rm d}V_{\mathbf{x}}
\end{split}
\end{sequation}

It is noted that the integration domain of the first integral on the right hand side of (\ref{dpf9}) is the volume constrained boundary that has a nonzero measure~\cite{tb59}. Thus, if we establish the PD-BEM directly based on (\ref{dpf9}), then we will lose the advantages of dimension reduction. Therefore, we introduce additional constraints to the volume constrained boundary~\cite{tb59}, which transform the integral in the volume constrained boundary into the one in the classical boundary for (\ref{dpf9}). To this end, as mentioned in Section \ref{atl21}, (\ref{fa13a}), (\ref{fa13ab}) and (\ref{fa14}) are the restrictions for the static problem. Similarly, the dynamic problem should also meet these restrictions. We further require them to be satisfied at any time. They are expressed as
\begin{sequation}\label{dpf10}
\begin{cases}
\mathbf{u} \left(\mathbf{x},t\right) = \overline{\mathbf{u}} \left(\mathbf{x},t\right)  & \text{in} \  \mathbf{x} \in \partial \Omega \\
\mathbf{T} \left(\mathbf{u}\right) \left(\mathbf{x},t\right) = \overline{\mathbf{T}} \left(\mathbf{v}\right) \left(\mathbf{x},t\right)  & \text{in} \  \mathbf{x} \in \partial \Omega
\end{cases}
\end{sequation}
\begin{sequation}\label{dpf11}
\int_{\Omega_\tau} \mathbf{v}_p \left(\mathbf{x},t'\right) \cdot \mathcal{N} \left(\boldsymbol{\Theta} : \mathcal{D}^{*} \left(\mathbf{v}_f\right)\right) \left(\mathbf{x},t\right) {\rm d}V_{\mathbf{x}} = \int_{\partial \Omega} \mathbf{v}_p \left(\mathbf{x},t'\right) \cdot \mathbf{T} \left(\mathbf{v}_f\right) \left(\mathbf{x},t\right) {\rm d}V_{\mathbf{x}}
\end{sequation}
where $\mathbf{T}$ is the nonlocal surface force operator; $\mathbf{v}_p$ is the possible deformed state and meets the displacement constraint. Replacing $\mathbf{v}_p$ in (\ref{dpf11}) with $\mathbf{v}_5$ in (\ref{dpf1}) and $\mathbf{v}_6$ in (\ref{dpf2}), we can get the following equations:
\begin{sequation}\label{dpf12}
\int_{\Omega_\tau} \mathbf{v}_5 \left(\mathbf{x},\tau\right) \cdot \mathcal{N} \left(\boldsymbol{\Theta} : \mathcal{D}^{*} \left(\mathbf{v}_6\right) \right) \left(\mathbf{x},t - \tau + t_0\right) {\rm d}V_{\mathbf{x}} = \int_{\partial \Omega} \mathbf{v}_5 \left(\mathbf{x},\tau\right) \cdot \mathbf{T} \left(\mathbf{v}_6\right) \left(\mathbf{x},t - \tau + t_0\right) {\rm d}V_{\mathbf{x}}
\end{sequation}
\begin{sequation}\label{dpf13}
\int_{\Omega_\tau} \mathbf{v}_6 \left(\mathbf{x},\tau\right) \cdot \mathcal{N} \left(\boldsymbol{\Theta} : \mathcal{D}^{*} \left(\mathbf{v}_5\right) \right) \left(\mathbf{x},t - \tau + t_0\right) {\rm d}V_{\mathbf{x}} = \int_{\partial \Omega} \mathbf{v}_6 \left(\mathbf{x},\tau\right) \cdot \mathbf{T} \left(\mathbf{v}_5\right) \left(\mathbf{x},t - \tau + t_0\right) {\rm d}V_{\mathbf{x}}
\end{sequation}
Substituting (\ref{dpf12}) and (\ref{dpf13}) into (\ref{dpf9}), we can construct the dynamic reciprocal theorem as follows:
\begin{sequation}\label{dpf14}
\begin{split}
& \int_{\Omega} \tilde{\mathbf{v}}_5 \left(\mathbf{x}\right) \cdot \dot{\mathbf{v}}_6 \left(\mathbf{x},t\right) + \dot{\tilde{\mathbf{v}}}_5 \left(\mathbf{x}\right) \cdot \mathbf{v}_6 \left(\mathbf{x},t\right) - \mathbf{v}_5 \left(\mathbf{x},t\right) * \mathbf{d} \left(\mathbf{x},t\right) {\rm d}V_{\mathbf{x}} \\
=& \int_{\Omega} \tilde{\mathbf{v}}_6 \left(\mathbf{x}\right) \cdot \dot{\mathbf{v}}_5 \left(\mathbf{x},t\right) + \dot{\tilde{\mathbf{v}}}_6 \left(\mathbf{x}\right) \cdot \mathbf{v}_5 \left(\mathbf{x},t\right) - \mathbf{v}_6 \left(\mathbf{x},t\right) * \mathbf{c} \left(\mathbf{x},t\right) {\rm d}V_{\mathbf{x}}\\
& + \int_{\partial \Omega} \mathbf{v}_5 \left(\mathbf{x},t\right) * \mathbf{T} \left(\mathbf{v}_6\right) \left(\mathbf{x},t\right) - \mathbf{v}_6 \left(\mathbf{x},t\right) * \mathbf{T} \left(\mathbf{v}_5\right) \left(\mathbf{x},t\right) {\rm d}V_{\mathbf{x}}
\end{split}
\end{sequation}
Similar to the static case, we can construct the fundamental solution in the finite domain and apply the reciprocal theorem to the fundamental solution and the real state. At last, we can get the dynamic boundary integral equation as follows:
\begin{sequation}\label{dpf17}
\begin{split}
{\rm u}_k \left(\mathbf{x},t\right) =& \int_{\Omega} \tilde{\mathbf{u}} \left(\mathbf{x}'\right) \cdot \dot{\mathbf{v}}_k \left(\mathbf{x}' - \mathbf{x},t\right) + \dot{\tilde{\mathbf{u}}} \left(\mathbf{x}'\right) \cdot \mathbf{v}_k \left(\mathbf{x}' - \mathbf{x},t\right) {\rm d}V_{\mathbf{x}'} + \int_{\Omega} \mathbf{v}_k \left(\mathbf{x}' - \mathbf{x},t\right) * \mathbf{f} \left(\mathbf{x}',t\right) {\rm d}V_{\mathbf{x}'} \\
& - \bbint_{\partial \Omega} \mathbf{u} \left(\mathbf{x}',t\right) * \mathbf{T} \left(\mathbf{v}_k\right) \left(\mathbf{x}' - \mathbf{x},t\right) - \mathbf{v}_k \left(\mathbf{x}' - \mathbf{x},t\right) * \mathbf{T} \left(\mathbf{u}\right) \left(\mathbf{x}',t\right) {\rm d}V_{\mathbf{x}'}
\end{split}
\end{sequation}
The governing equation and boundary conditions of the fundamental solution $\mathbf{v}_k$ and the real state ${\mathbf{u}}$ are described as follows:
\begin{sequation}\label{dpf15}
Problem \ 7 :
\begin{cases}
\mathcal{D} \left(\boldsymbol{\Theta} : \mathcal{D}^{*} \left(\mathbf{u}\right) \right) \left(\mathbf{x},t\right) + \mathbf{f} \left(\mathbf{x},t\right) = \rho \ddot{\mathbf{u}} \left(\mathbf{x},t\right)   \qquad \qquad \, \, \,  & \text{in} \ \mathbf{x} \in \Omega \\
\mathbf{u} \left(\mathbf{x},t\right) = \mathbf{g}_{d} \left(\mathbf{x},t\right)   \qquad \qquad \, \, \,  & \text{in} \ \mathbf{x} \in \Omega_{\tau_d} \\
\mathcal{N} \left(\boldsymbol{\Theta} : \mathcal{D}^{*} \left(\mathbf{u}\right) \right) \left(\mathbf{x},t\right) = \mathbf{g}_{n} \left(\mathbf{x},t\right)   \qquad \qquad \, \, \,  & \text{in} \ \mathbf{x} \in \Omega_{\tau_n} \\
\mathbf{u} \left(\mathbf{x},t_0\right) = \tilde{\mathbf{u}} \left(\mathbf{x}\right) \qquad \dot{\mathbf{u}} \left(\mathbf{x},t_0\right) = \dot{\tilde{\mathbf{u}}} \left(\mathbf{x}\right)
\end{cases}
\end{sequation}
\begin{sequation}\label{dpf16}
Problem \ 8 :
\begin{cases}
\mathcal{D} \left(\boldsymbol{\Theta} : \mathcal{D}^{*} \left(\mathbf{v}_k\right) \right) \left(\mathbf{x} - \mathbf{x}_0,t\right) = \rho \ddot{\mathbf{v}}_k \left(\mathbf{x} - \mathbf{x}_0,t\right)  & \text{in} \ \mathbf{x} \in \Omega \\
\mathbf{v}_k \left(\mathbf{x} - \mathbf{x}_0,t\right) = \overline{\mathbf{v}}_k \left(\mathbf{x} - \mathbf{x}_0,t\right)  & \text{in} \ \mathbf{x} \in \Omega_{\tau_d} \\
\mathcal{N} \left(\boldsymbol{\Theta} : \mathcal{D}^{*} \left(\mathbf{v}_k\right) \right) \left(\mathbf{x} - \mathbf{x}_0,t\right) = \mathcal{N} \left(\boldsymbol{\Theta} : \mathcal{D}^{*} \left(\overline{\mathbf{v}}_k\right) \right) \left(\mathbf{x} - \mathbf{x}_0,t\right)  & \text{in} \ \mathbf{x} \in \Omega_{\tau_n} \\
\mathbf{v}_k \left(\mathbf{x} - \mathbf{x}_0,t_0\right) = \mathbf{0} \qquad \dot{\mathbf{v}}_k \left(\mathbf{x} - \mathbf{x}_0,t_0\right) = \delta \left(\mathbf{x} - \mathbf{x}_0\right) \mathbf{e}_k
\end{cases}
\end{sequation}
where $\delta \left(\mathbf{x} - \mathbf{x}_0\right)$ is the Dirac function. $\mathbf{x}_0 \in \Omega$. $\tilde{\mathbf{u}}$ and $\dot{\tilde{\mathbf{u}}}$ are the initial conditions. $\mathbf{e}_k$ is a coordinate base vector. $\overline{\mathbf{v}}_k$ is the infinite fundamental solution which satisfies the following equation:
\begin{sequation}\label{dpfa16}
\begin{cases}
\mathcal{D} \left(\boldsymbol{\Theta} : \mathcal{D}^{*} \left(\overline{\mathbf{v}}_k\right) \right) \left(\mathbf{x} - \mathbf{x}_0,t\right) = \rho \ddot{\overline{\mathbf{v}}}_k \left(\mathbf{x} - \mathbf{x}_0,t\right) \qquad  \text{in} \ \mathbf{x} \in \Omega \\
\overline{\mathbf{v}}_k \left(\mathbf{x} - \mathbf{x}_0,t_0\right) = \mathbf{0} \\
\dot{\overline{\mathbf{v}}}_k \left(\mathbf{x} - \mathbf{x}_0,t_0\right) = \delta \left(\mathbf{x} - \mathbf{x}_0\right) \mathbf{e}_k
\end{cases}
\end{sequation}
where $\mathbf{x}_0, \mathbf{x}\in R^n$.

\section{Numerical Laplace transformation inversion based on Fourier cosine transformation}\label{ap9}

Here, we recapitulate the main steps of the method of Laplace transformation inversion in the literature~\cite{tb69}. The Laplace transformation inversion and the Laplace transformation for a real function are
\begin{sequation}\label{hpg1}
f\left(t\right) = \dfrac{2{\rm e}^{at}}{\pi} \int_{0}^{\infty} Re\left\{F\left(s\right)\right\} {\rm cos} \left(\omega t\right) {\rm d}\omega
\end{sequation}
\begin{sequation}\label{hpg2}
Re\left\{F\left(s\right)\right\} = \int_{0}^{\infty} {\rm e}^{-at} f\left(t\right) {\rm cos} \left(\omega t\right) {\rm d}t
\end{sequation}
where $s = a + {\rm i}\omega$. $f\left(t\right)$ is a function of $t$ with its Laplace transformation $F\left(s\right)$. It is noted that (\ref{hpg2}) can also be regard as the Fourier cosine transformation for the function ${\rm e}^{-at} f\left(t\right)$, which is the core of this algorithm. Then, the main steps for the inversion are~\cite{tb69}
\begin{enumerate}
	\item Regarding the Laplace transformation inversion as the Fourier cosine transformation by reshaping $f \left(t\right)$ as a new function $g \left(t\right)$.	
	\item Carrying out the periodic extension for $g \left(t\right)$ in the time interval.
	\item Implementing the Fourier expansion for the periodic function that is obtained in step 2.
	\item Transforming the function $g \left(t\right)$ in step 3 back to the function $f \left(t\right)$.
\end{enumerate}
The operations according to the above four steps are described as follows.  First, two new functions are introduced
\begin{sequation}\label{hpg3}
h\left(t\right) = {\rm e}^{-at} f\left(t\right)
\end{sequation}
\begin{sequation}\label{hpg4}
g\left(t\right) = {\rm e}^{-at} fp\left(t\right)
\end{sequation}
where $fp\left(t\right)$ is a function formed by periodic extension for the interval in which the value of the function $f\left(t\right)$ is to be calculated.  Executing the Fourier cosine expansion for $g\left(t\right)$, we can get
\begin{sequation}\label{hpg5}
g\left(t\right) = \dfrac{2}{T} \left[\dfrac{A\left(\omega_0\right)}{2} + \sum_{k = 1}^{\infty} A\left(\omega_k\right) {\rm cos} \dfrac{k\pi t}{T} \right]
\end{sequation}
In (\ref{hpg5}), $2T$ is the length of the interval in which the value of the function $f\left(t\right)$ is to be calculated. $A\left(\omega_k\right)$ is expressed as follows:
\begin{sequation}\label{hpg6}
A\left(\omega_k\right) = \int_{0}^{\infty} h\left(t\right) {\rm cos} \left(\dfrac{k\pi t}{T}\right) {\rm d}t
\end{sequation}
From (\ref{hpg2}), we know
\begin{sequation}\label{hpg7}
A\left(\omega_k\right) = Re\left\{F\left(s\right)\right\}
\end{sequation}
Substituting (\ref{hpg7}) and (\ref{hpg4}) into (\ref{hpg5}), we can get
\begin{sequation}\label{hpg8}
fp\left(t\right) = \dfrac{2{\rm e}^{at}}{T} \left[\dfrac{Re\left\{F\left(a\right)\right\}}{2} + \sum_{k = 1}^{\infty} Re\left\{F\left(a + \dfrac{k\pi {\rm i}}{T}\right)\right\}{\rm cos} \dfrac{k\pi t}{T} \right]
\end{sequation}
By the summation of the truncated series, the Laplace transformation inversion formula is obtained
\begin{sequation}\label{hpg9}
f\left(t\right) \approx \dfrac{2{\rm e}^{at}}{T} \left[\dfrac{Re\left\{F\left(a\right)\right\}}{2} + \sum_{k = 1}^{N} Re\left\{F\left(a + \dfrac{k\pi {\rm i}}{T}\right)\right\}{\rm cos} \dfrac{k\pi t}{T} \right]
\end{sequation}

\section{The reciprocal theorem for state-based PD}\label{ap8}

Firstly, we give the Navier equation for the linearized state-based PD as follows~\cite{tb74}:
\begin{sequation}\label{hpf1}
\begin{cases}
\mathcal{D} \left(\eta \varpi \left(\mathcal{D}^{*} \left(\mathbf{u}\right)\right)^{T} + \omega\sigma Tr\left(\mathcal{D}^{*}_{\omega} \left(\mathbf{u}\right)\right)\mathbf{I} \right) \left(\mathbf{x}\right) + \mathbf{b} \left(\mathbf{x}\right) = \mathbf{0} \qquad & \text{in} \ \mathbf{x} \in \Omega \\
\mathbf{u} \left(\mathbf{x}\right) = \mathbf{h} \left(\mathbf{x}\right) \qquad & \text{in} \ \mathbf{x} \in \Omega_{\tau_d} \\
\mathcal{N} \left(\eta \varpi \left(\mathcal{D}^{*} \left(\mathbf{u}\right)\right)^{T} + \omega\sigma Tr\left(\mathcal{D}^{*}_{\omega} \left(\mathbf{u}\right)\right)\mathbf{I} \right) \left(\mathbf{x}\right) = \mathbf{g} \left(\mathbf{x}\right) \qquad & \text{in} \ \mathbf{x} \in \Omega_{\tau_n}
\end{cases}
\end{sequation}
where $\varpi$ is the influence function; $\mathbf{I}$ is the metric tensor; $Tr$ is the trace operator; $T$ is the transposed symbol;  $\mathbf{b}\left(\mathbf{x}\right)$ is the body force density; $\mathbf{g} \left(\mathbf{x}\right)$ is the force constraint on the boundary; $\mathbf{h} \left(\mathbf{x}\right)$ is the displacement constraint on the boundary. The definition of the nonlocal operator is shown in the literature~\cite{tb59}. $\mathcal{D}$, $\mathcal{D}^{*}$ and $\mathcal{N}$ have been given in Section \ref{atl21}. We give the definition of $\mathcal{D}^{*}_{\omega}$ as follows:
\begin{sequation}\label{hpf2}
\mathcal{D}^{*}_{\omega} \left(\mathbf{x}\right) = \int_{\Omega_\tau \cup \Omega} \mathcal{D}^{*} \left(\mathbf{u}\right) \left(\mathbf{x},\mathbf{x}'\right) \omega \left(\mathbf{x},\mathbf{x}'\right) {\rm d}V_{\mathbf{x}'}
\end{sequation}
where
\begin{sequation}\label{hpf3}
\omega \left(\mathbf{x},\mathbf{x}'\right) = \dfrac{\left| \mathbf{x}' - \mathbf{x} \right| \varpi \left(\mathbf{x},\mathbf{x}'\right)}{\phi \left(\mathbf{x}\right)} \qquad \qquad \qquad \boldsymbol{\alpha} \left(\mathbf{x},\mathbf{x}'\right) = \dfrac{\mathbf{x} - \mathbf{x}'}{\left| \mathbf{x}' - \mathbf{x} \right|}
\end{sequation}
For the three-dimensional case, the parameters in the Navier equation are
\begin{sequation}\label{hpf4}
\sigma = \kappa - \dfrac{\eta\phi}{3} \qquad \qquad \phi \left(\mathbf{x}\right) = \dfrac{1}{3} m \left(\mathbf{x}\right) \qquad \qquad \eta = \dfrac{15\mu}{m \left(\mathbf{x}\right)}
\end{sequation}
where $\kappa$ is the bulk modulus; $\mu$ is the shear modulus; and $m$ is
\begin{sequation}\label{hpf5}
m \left(\mathbf{x}\right) = \int_{\Omega_\tau \cup \Omega} \left| \boldsymbol{\xi} \right|^2 \varpi \left(\mathbf{x},\boldsymbol{\xi}\right) {\rm d}V_{\boldsymbol{\xi}}
\end{sequation}
For the two-dimensional case, the parameters in the Navier equation are
\begin{sequation}\label{hpf6}
\sigma = \kappa - \dfrac{\eta\phi}{2} \qquad \qquad \phi \left(\mathbf{x}\right) = \dfrac{1}{2} m \left(\mathbf{x}\right) \qquad \qquad \eta = \dfrac{8\mu}{m \left(\mathbf{x}\right)}
\end{sequation}
Next, we derive the reciprocal theorem. Firstly, we give the following two states:
\begin{sequation}\label{hpf7}
Case  1 :
\begin{cases}
\mathcal{D} \left(\eta \varpi \left(\mathcal{D}^{*} \left(\mathbf{v}_1\right)\right)^{T} + \omega\sigma Tr\left(\mathcal{D}^{*}_{\omega} \left(\mathbf{v}_1\right)\right)\mathbf{I} \right) \left(\mathbf{x}\right) + \mathbf{b}_1 \left(\mathbf{x}\right) = \mathbf{0} \qquad & \text{in} \ \mathbf{x} \in \Omega \\
\mathbf{v}_1 \left(\mathbf{x}\right) = \mathbf{h}_1 \left(\mathbf{x}\right) \qquad & \text{in} \ \mathbf{x} \in \Omega_{\tau_d} \\
\mathcal{N} \left(\eta \varpi \left(\mathcal{D}^{*} \left(\mathbf{v}_1\right)\right)^{T} + \omega\sigma Tr\left(\mathcal{D}^{*}_{\omega} \left(\mathbf{v}_1\right)\right)\mathbf{I} \right) \left(\mathbf{x}\right) = \mathbf{g}_1 \left(\mathbf{x}\right) \qquad & \text{in} \ \mathbf{x} \in \Omega_{\tau_n}
\end{cases}
\end{sequation}
\begin{sequation}\label{hpf8}
Case  2 :
\begin{cases}
\mathcal{D} \left(\eta \varpi \left(\mathcal{D}^{*} \left(\mathbf{v}_2\right)\right)^{T} + \omega\sigma Tr\left(\mathcal{D}^{*}_{\omega} \left(\mathbf{v}_2\right)\right)\mathbf{I} \right) \left(\mathbf{x}\right) + \mathbf{b}_2 \left(\mathbf{x}\right) = \mathbf{0} \qquad & \text{in} \ \mathbf{x} \in \Omega \\
\mathbf{v}_2 \left(\mathbf{x}\right) = \mathbf{h}_2 \left(\mathbf{x}\right) \qquad & \text{in} \ \mathbf{x} \in \Omega_{\tau_d} \\
\mathcal{N} \left(\eta \varpi \left(\mathcal{D}^{*} \left(\mathbf{v}_2\right)\right)^{T} + \omega\sigma Tr\left(\mathcal{D}^{*}_{\omega} \left(\mathbf{v}_2\right)\right)\mathbf{I} \right) \left(\mathbf{x}\right) = \mathbf{g}_2 \left(\mathbf{x}\right) \qquad & \text{in} \ \mathbf{x} \in \Omega_{\tau_n}
\end{cases}
\end{sequation}
The virtual work produced by $Case\ \ 1$ and $Case\ \ 2$ is
\begin{sequation}\label{hpf9}
\begin{split}
& \int_{\Omega} \mathbf{v}_1 \left(\mathbf{x}\right) \cdot \mathcal{D} \left(\eta \varpi \left(\mathcal{D}^{*} \left(\mathbf{v}_2\right)\right)^{T} + \omega\sigma Tr\left(\mathcal{D}^{*}_{\omega} \left(\mathbf{v}_2\right)\right)\mathbf{I} \right) \left(\mathbf{x}\right) {\rm d}V_{\mathbf{x}} \\
=& \int_{\Omega_\tau \cup \Omega} \int_{\Omega_\tau \cup \Omega} \mathcal{D}^{*} \left(\mathbf{v}_1\right) \left(\mathbf{x},\mathbf{x}'\right) : \left(\eta \varpi \left(\mathcal{D}^{*} \left(\mathbf{v}_2\right)\right)^{T} + \omega\sigma Tr\left(\mathcal{D}^{*}_{\omega} \left(\mathbf{v}_2\right)\right)\mathbf{I} \right) \left(\mathbf{x},\mathbf{x}'\right) {\rm d}V_{\mathbf{x}'} {\rm d}V_{\mathbf{x}} \\
& + \int_{\Omega_\tau} \mathbf{v}_1\left(\mathbf{x}\right) \cdot \mathcal{N} \left(\eta \varpi \left(\mathcal{D}^{*} \left(\mathbf{v}_2\right)\right)^{T} + \omega\sigma Tr\left(\mathcal{D}^{*}_{\omega} \left(\mathbf{v}_2\right)\right)\mathbf{I} \right) \left(\mathbf{x}\right) {\rm d}V_{\mathbf{x}} \\
=& \int_{\Omega_\tau \cup \Omega} \sigma Tr\left(\mathcal{D}^{*}_{\omega} \left(\mathbf{v}_2\right)\right)Tr\left(\mathcal{D}^{*}_{\omega} \left(\mathbf{v}_1\right)\right) {\rm d}V_{\mathbf{x}} + \int_{\Omega_\tau} \mathbf{v}_1\left(\mathbf{x}\right) \cdot \mathcal{N} \left(\eta \varpi \left(\mathcal{D}^{*} \left(\mathbf{v}_2\right)\right)^{T} + \omega\sigma Tr\left(\mathcal{D}^{*}_{\omega} \left(\mathbf{v}_2\right)\right)\mathbf{I} \right) \left(\mathbf{x}\right) {\rm d}V_{\mathbf{x}} \\
& + \int_{\Omega_\tau \cup \Omega} \int_{\Omega_\tau \cup \Omega} \eta \varpi \left(\mathbf{x},\mathbf{x}'\right) \left(\left( \mathbf{v}_1 \left(\mathbf{x}'\right) - \mathbf{v}_1 \left(\mathbf{x}\right) \right) \cdot \boldsymbol{\alpha} \left(\mathbf{x},\mathbf{x}'\right) \right) \left(\left( \mathbf{v}_2 \left(\mathbf{x}'\right) - \mathbf{v}_2 \left(\mathbf{x}\right) \right) \cdot \boldsymbol{\alpha} \left(\mathbf{x},\mathbf{x}'\right) \right) {\rm d}V_{\mathbf{x}'} {\rm d}V_{\mathbf{x}}
\end{split}
\end{sequation}
\begin{sequation}\label{hpf10}
\begin{split}
& \int_{\Omega} \mathbf{v}_2 \left(\mathbf{x}\right) \cdot \mathcal{D} \left(\eta \varpi \left(\mathcal{D}^{*} \left(\mathbf{v}_1\right)\right)^{T} + \omega\sigma Tr\left(\mathcal{D}^{*}_{\omega} \left(\mathbf{v}_1\right)\right)\mathbf{I} \right) \left(\mathbf{x}\right) {\rm d}V_{\mathbf{x}} \\
=& \int_{\Omega_\tau \cup \Omega} \sigma Tr\left(\mathcal{D}^{*}_{\omega} \left(\mathbf{v}_1\right)\right)Tr\left(\mathcal{D}^{*}_{\omega} \left(\mathbf{v}_2\right)\right) {\rm d}V_{\mathbf{x}} + \int_{\Omega_\tau} \mathbf{v}_2\left(\mathbf{x}\right) \cdot \mathcal{N} \left(\eta \varpi \left(\mathcal{D}^{*} \left(\mathbf{v}_1\right)\right)^{T} + \omega\sigma Tr\left(\mathcal{D}^{*}_{\omega} \left(\mathbf{v}_1\right)\right)\mathbf{I} \right) \left(\mathbf{x}\right) {\rm d}V_{\mathbf{x}} \\
& + \int_{\Omega_\tau \cup \Omega} \int_{\Omega_\tau \cup \Omega} \eta \varpi \left(\mathbf{x},\mathbf{x}'\right) \left(\left( \mathbf{v}_2 \left(\mathbf{x}'\right) - \mathbf{v}_2 \left(\mathbf{x}\right) \right) \cdot \boldsymbol{\alpha} \left(\mathbf{x},\mathbf{x}'\right) \right) \left(\left( \mathbf{v}_1 \left(\mathbf{x}'\right) - \mathbf{v}_1 \left(\mathbf{x}\right) \right) \cdot \boldsymbol{\alpha} \left(\mathbf{x},\mathbf{x}'\right) \right) {\rm d}V_{\mathbf{x}'} {\rm d}V_{\mathbf{x}}
\end{split}
\end{sequation}
(\ref{hpf9}) corresponds to the virtual work generated by the load in $Case\ \ 2$ on the displacement in $Case\ \ 1$, and (\ref{hpf10}) corresponds to the virtual work generated by the load in $Case\ \ 1$ on the displacement in $Case\ \ 2$.  Subtracting (\ref{hpf9}) from (\ref{hpf10}), we can get the reciprocal theorem as follows:
\begin{sequation}\label{hpf11}
\begin{split}
& \int_{\Omega} \mathbf{v}_2 \left(\mathbf{x}\right) \cdot \mathcal{D} \left(\eta \varpi \left(\mathcal{D}^{*} \left(\mathbf{v}_1\right)\right)^{T} + \omega\sigma Tr\left(\mathcal{D}^{*}_{\omega} \left(\mathbf{v}_1\right)\right)\mathbf{I} \right) \left(\mathbf{x}\right) {\rm d}V_{\mathbf{x}} \\
& - \int_{\Omega} \mathbf{v}_1 \left(\mathbf{x}\right) \cdot \mathcal{D} \left(\eta \varpi \left(\mathcal{D}^{*} \left(\mathbf{v}_2\right)\right)^{T} + \omega\sigma Tr\left(\mathcal{D}^{*}_{\omega} \left(\mathbf{v}_2\right)\right)\mathbf{I} \right) \left(\mathbf{x}\right) {\rm d}V_{\mathbf{x}} \\
=& \int_{\Omega_\tau} \mathbf{v}_2\left(\mathbf{x}\right) \cdot \mathcal{N} \left(\eta \varpi \left(\mathcal{D}^{*} \left(\mathbf{v}_1\right)\right)^{T} + \omega\sigma Tr\left(\mathcal{D}^{*}_{\omega} \left(\mathbf{v}_1\right)\right)\mathbf{I} \right) \left(\mathbf{x}\right) {\rm d}V_{\mathbf{x}} \\
& - \int_{\Omega_\tau} \mathbf{v}_1\left(\mathbf{x}\right) \cdot \mathcal{N} \left(\eta \varpi \left(\mathcal{D}^{*} \left(\mathbf{v}_2\right)\right)^{T} + \omega\sigma Tr\left(\mathcal{D}^{*}_{\omega} \left(\mathbf{v}_2\right)\right)\mathbf{I} \right) \left(\mathbf{x}\right) {\rm d}V_{\mathbf{x}}
\end{split}
\end{sequation}
The derivation and detailed form of the boundary integral equation are the same as those of the bond-based PD. The proof of the equivalence between the linearized state-based PD~\cite{tb75} and the Navier equation is shown in the literature~\cite{tb74}.

\end{document}